\documentclass[12pt,twoside]{amsart}
\usepackage{amssymb}

\numberwithin{equation}{section}
\font\tengothic=eufm10 scaled\magstep 1
\font\sevengothic=eufm7 scaled\magstep 1
\newfam\gothicfam
       \textfont\gothicfam=\tengothic
       \scriptfont\gothicfam=\sevengothic

\newtheorem{theorem}{Theorem}[section]
\newtheorem{lemma}[theorem]{Lemma}
\newtheorem{proposition}[theorem]{Proposition}
\newtheorem{corollary}[theorem]{Corollary}
\newtheorem{conjecture}[theorem]{Conjecture}

\theoremstyle{definition}
\newtheorem{definition}[theorem]{Definition} 
\newtheorem{remark}[theorem]{Remark}
\newtheorem{example}[theorem]{Example}

\newtheorem{notation}[theorem]{Notation}
\newtheorem{question}[theorem]{Question}
\newtheorem{mquestion}[theorem]{Main question}

\newcommand{\codim}{\operatorname{codim}}
\newcommand{\coker}{\operatorname{coker}}
\newcommand{\Ann}{\operatorname{Ann}}

\newcommand{\Ass}{\operatorname{Ass}}
\newcommand{\Hom}{\operatorname{Hom}}
\newcommand{\Ext}{\operatorname{Ext}}
\newcommand{\reg}{\operatorname{reg}}

\newcommand{\rank}{\operatorname{rank}}
\newcommand{\rankk}{\operatorname{rank}_K}
\newcommand{\depth}{\operatorname{depth}}

\newcommand{\im}{\operatorname{im}}

\newcommand{\HH}{H_{\mathfrak m}}

\newcommand{\ER}{\operatorname{Ext}_R}

\newcommand{\Proj}{\operatorname{Proj}}
\newcommand{\proj}[1]
{ \mathchoice
            { {\mathbb P}^{#1} }
            { {\mathbb P}^{#1} }
            { {\mathbb P}^{#1} }
            { {\mathbb P}^{#1} }
          }

\newcommand{\s}{\; | \;}
\newcommand{\mif}{\mbox{if} ~}
\newcommand{\fora}{\quad \mbox{for all}\;\;}

\newcommand{\ffi}{\varphi}
\newcommand{\veps}{\varepsilon}

\newcommand{\Rad}{\operatorname{Rad}}

\newcommand{\cI}{{\mathcal I}}

\newcommand{\cF}{{\mathcal F}}

\newcommand{\cL}{{\mathcal L}}

\newcommand{\cM}{{\mathcal M}}
\newcommand{\fm}{{\mathfrak m}}
\newcommand{\fa}{{\mathfrak a}}
\newcommand{\fb}{{\mathfrak b}}
\newcommand{\fc}{{\mathfrak c}}
\newcommand{\fq}{{\mathfrak q}}
\newcommand{\fp}{{\mathfrak p}}
\newcommand{\fd}{{\mathfrak d}}

\newcommand{\acm}{arithmetically Cohen-Macaulay}
\newcommand{\aCM}{arithmetically Cohen-Macaulay}
\newcommand{\ACM}{arithmetically Cohen-Macaulay}
\newcommand{\lCM}{locally Cohen-Macaulay}
\newcommand{\CM}{Cohen-Macaulay}
\newcommand{\aG}{arithmetically Gorenstein}
\newcommand{\aBM}{arithmetically Buchsbaum}
\newcommand{\aci}{almost complete intersection}
\newcommand{\HR}{Hartshorne-Rao module}

\newcommand {\ZZ}{\mathbb{Z}}
\newcommand {\NN}{\mathbb{N}}
\newcommand {\PP}{\mathbb{P}}

\begin{document}
\title[Liaison and Related Topics: Notes from the Torino
Workshop/School]{Liaison and Related Topics: Notes from the Torino
Workshop/School}

\author[J.\ Migliore, U.\ Nagel]{J.\ Migliore, U.\ Nagel}
\address{Department of Mathematics,
         University of Notre Dame,
         Notre Dame, IN 46556,
         USA}
\email{migliore.1@nd.edu}
\address{Fachbereich Mathematik und Informatik, Universit\"at
Paderborn, D--33095 Paderborn, Germany}
\email{uwen@uni-paderborn.de
}







\maketitle

\tableofcontents


  \section{Introduction} \label{intro}

These are the expanded and detailed notes of the lectures given by the authors
during the school and workshop entitled ``Liaison and Related Topics," held at
the Politecnico di Torino during the period October 1-5, 2001.

The authors each gave five lectures of length 1.5 hours each.  We attempted to
cover liaison theory from first principles, through the main developments
(especially in codimension two) and the standard applications, to the recent
developments in Gorenstein liaison and a discussion of open problems.  Given
the extensiveness of the subject, it was not possible to go into great detail
in every proof.  Still, it is hoped that the material that we chose will be
beneficial and illuminating for the participants, and for the reader.

We believe that these notes will be a valuable addition to the literature, and
give details and points of view that cannot be found in other expository works
on this subject.  Still, we would like to point out that a number of such works
do exist.  In particular, the interested reader should also consult
\cite{hartshorne-exp}, \cite{migbook}, \cite{kyoto}, \cite{rosa1},
\cite{rosa2}.
\smallskip

We are going to describe the contents of these notes. In the expository
Section 2 we discuss the origins of liaison theory, its scope and several
results and problems which are more carefully treated in later sections.

Sections 3 and 4 have preparatory character. We recall several results
which are used later on. In Section 3 we discuss in particular the relation
between local and sheaf  cohomology, and modules and sheaves. Sections 4 is
devoted to Gorenstein ideals where among other things we describe various
constructions of such ideals.

  The discussions of liaison theory begins in Section 5. Besides giving the
basic definitions we state the first results justifying the name, i.e. showing
that indeed the properties of directly linked schemes can be related to each
other.

  Two key results of Gorenstein liaison are presented in Section 6: the somewhat
  surprisingly general version of basic double linkage and the fact that
linearly equivalent divisors on ``nice'' \aCM\ subschemes are Gorenstein linked
in two steps.

  The equivalence classes generated by the various concepts of linkage are
discussed in Sections 7 - 10. Rao's correspondence is explained in Section 7.
It is a relation between  even liaison classes  and certain reflexive
modules/sheaves which gives necessary conditions on two subschemes for being
linked in an even number of steps. In Section 8 it is shown that these
conditions are also sufficient for subschemes of codimension two. It is the
main open problem of Gorenstein liaison to decide if this is also true for
subschemes of higher codimension.  Several results are mentioned which provide
evidence for an affirmative answer. Examples show that the answer is negative
if one links by complete intersections only. In Section 9 we consider the
structure of an even liaison class. For subschemes of codimension
two it is described by the Lazarsfeld-Rao property. Moreover, we discuss the
  possibility of extending it to subschemes of higher codimension. In
Section 10 we compare the equivalence relations generated by the different
concepts of linkage. In particular, we explain how invariants for complete
intersection liaison can be used to distinguish complete intersection liaison
classes within one Gorenstein liaison class.

  Section 11 gives a flavour of the various applications of liaison theory.

  Throughout these notes we mention various open problems. Some of them and
  further problems related to liaison theory are stated in Section 12.

  Although most of the results are true more generally for subschemes of an \aG\
  subscheme, for simplicity we restrict ourselves to subschemes of $\PP^n$.
\smallskip

Both authors were honored and delighted to be invited to give the lectures for
this workshop.  We are grateful to the main organizers, Gianfranco Casnati,
Nadia Chiarli and Silvio Greco, for their kind hospitality.  We are also
grateful to the participants, especially Roberto Notari and Maria Luisa
Spreafico, for their hospitality and mathematical discussions, and for their
hard work in preparing this volume.  Finally, we are grateful to Robin
Hartshorne and Rosa Mir\'o-Roig for helpful comments about the contents of
these notes, and especially to Hartshorne for his Example \ref{hartshorne ex}.


\section{Overview and History} \label{sec-overv}

This section will give an expository overview of the subject of liaison theory,
and the subsequent sections will provide extensive detail.  Liaison theory has
its roots dating to more than a century ago.  The greatest activity, however,
has been in the last quarter century, beginning with the work of Peskine and
Szpir\'o
\cite{PS} in 1974.  There are at least three perspectives on liaison that we
hope to stress in these notes:

\begin{itemize}

\item Liaison is a very interesting subject in its own right.  There are many
hard open problems, and recently there is hope for a broad theory in arbitrary
codimension that neatly encompasses the codimension two case, where a
fairly complete picture has been understood for many years.

\item Liaison is a powerful tool for constructing examples.  Sometimes a
hypothetical situation arises but it is not known if a concrete example exists
to fit the theoretical constraints.  Liaison is often used to find such an
example.

\item Liaison is a useful method of proof.  It often happens that one can study
an object by linking to something which is intrinsically easier to study.  It
is also a useful method of proving that an object does not exist, because if it
did then a link would exist to something which can be proved to be
non-existent.

\end{itemize}

Let $R = K[x_0,\dots,x_n]$ where $K$ is a field.  For a sheaf ${\mathcal F}$ of
${\mathcal O}_{\proj{n}}$-modules, we set
\[
H^i_* ({\mathcal F}) = \bigoplus_{t \in {\mathbb Z}} H^i(\proj{n}, {\mathcal
F}(t))
\]
This is a graded $R$-module.  One use of this module comes in the following
notion.

\begin{definition}
A subscheme $X \subset \proj{n}$ is {\em arithmetically
Cohen-Macaulay} if $R/I_X$ is a Cohen-Macaulay ring, i.e.\ $\dim R/I = \depth
R/I$, where $\dim$ is the Krull-dimension.
\end{definition}

These notions will be discussed in greater detail in coming sections.  We will
see in Section 3 that $X$ is \ACM\ if and only if $H^i_*({\mathcal I}_X) = 0$
for $1 \leq i \leq \dim X$.  When $X$ is \ACM\ of codimension $c$, say, the
minimal free resolution of $I_X$ is as short as possible:
\[
0 \rightarrow F_c \rightarrow F_{c-1} \rightarrow \cdots \rightarrow F_1
\rightarrow I_X \rightarrow 0.
\]
(This follows from the Auslander-Buchsbaum theorem and the definition of a
Cohen-Macaulay ring.)  The {\em Cohen-Macaulay type} of $X$, or of $R/I_X$, is
the rank of $F_c$.  We will take as our definition that $X$ is {\em
arithmetically Gorenstein} if $X$ is \ACM\ of Cohen-Macaulay type 1, although
in Section 4 we will see equivalent formulations (Proposition \ref{equiv to
gor}).  For example, thanks to the Koszul resolution we know that a complete
intersection is always arithmetically Gorenstein.  The converse holds only in
codimension two.  We will discuss these notions again later, but we assume
these
basic ideas for the current discussion.

Liaison is, roughly, the study of unions of subschemes, and in particular what
can be determined if one knows that the union is ``nice.'' Let us begin with a
very simple situation.   Let $C_1$ and $C_2$ be equidimensional subschemes in
$\proj{n}$  with saturated ideals $I_{C_1}, I_{C_2} \subset R$ (i.e.\ $I_{C_1}$
and $I_{C_2}$ are unmixed homogeneous ideals in $R$).  We assume that $C_1$ and
$C_2$ have no common component.  We can study the union $X = C_1 \cup C_2$,
with
saturated ideal $I_X = I_{C_1} \cap I_{C_2}$, and the intersection $Z = C_1
\cap
C_2$, defined by the ideal $I_{C_1} + I_{C_2}$.  Note that this latter ideal is
not necessarily saturated, so $I_Z = (I_{C_1} + I_{C_2})^{sat}$.  These are
related by the exact sequence
\begin{equation} \label{int dirsum sum}
0 \rightarrow I_{C_1} \cap I_{C_2} \rightarrow I_{C_1} \oplus I_{C_2}
\rightarrow I_{C_1} + I_{C_2} \rightarrow 0.
\end{equation}
Sheafifying gives
\[
0 \rightarrow {\mathcal I}_X \rightarrow {\mathcal I}_{C_1} \oplus {\mathcal
I}_{C_2} \rightarrow {\mathcal I}_Z \rightarrow 0.
\]
Taking cohomology and forming a direct sum over all twists, we get
\[
\begin{array}{rclccccccccccccccccccccccc}
0 \ \rightarrow \ I_X \ \rightarrow \ I_{C_1} \oplus I_{C_2} & \hskip -.4cm
\longrightarrow & \hskip -.4cm I_Z \  \rightarrow \ H^1_*({\mathcal I}_X) \
\rightarrow \ H^1_*({\mathcal I}_{C_1}) \oplus H^1_*({\mathcal I}_{C_2}) \
\rightarrow
\cdots
\\
& \hskip -1.3cm \searrow & \hskip -.9cm \nearrow \\
& \hskip 0cm I_{C_1} + I_{C_2} \\
& \hskip -1.3cm \nearrow & \hskip -.9cm \searrow \\
& \hskip -.4cm  0 \hfill & \hskip -.3cm 0
\end{array}
\]
So one can see immediately that somehow $H^1_*({\mathcal I}_X)$
(or really a submodule) measures the failure of $I_{C_1} +
I_{C_2}$ to be saturated, and that if this cohomology is zero then
the ideal is saturated.  More observations about how submodules of
$H^1_*({\mathcal I}_X)$ measure various deficiencies can be found
in \cite{migbook}.

\begin{remark}\label{observations}
We can make the following observations about our union $X = C_1 \cup C_2$:

\begin{enumerate}
\item  If $H^1_*({\mathcal I_X}) = 0$ (in particular if $X$ is \ACM) then
$I_{C_1} + I_{C_2} = I_Z$ is saturated.

\item $I_X \subset I_{C_1}$ and $I_X \subset I_{C_2}$.

\item \label{ideal quotient} $[I_X : I_{C_1}] = I_{C_2}$ and $[I_X : I_{C_2}] =
I_{C_1}$ since $C_1$ and $C_2$ have no common component (cf.\ \cite{CLO} page
192).

\item It is not hard to see that we have an exact sequence
\[
0 \rightarrow R/I_X \rightarrow R/I_{C_1} \oplus R/I_{C_2} \rightarrow
R/(I_{C_1} + I_{C_2}) \rightarrow 0.
\]
Hence we get the relations
\[
\begin{array}{rcl}
\deg C_1 + \deg C_2 & = & \deg X \\
p_a C_1 + p_a C_2 & = & p_a X + 1 - \deg Z \ \ \hbox{(if $C_1$ and $C_2$ are
curves)}
\end{array}
\]
where $p_a$ represents the arithmetic genus.

\item \label{acm not pres} Even if $X$ is \ACM, it is possible that $C_1$ is
\ACM\ but $C_2$ is not \ACM.  For instance, consider the case where $C_2$ is
the disjoint union of two lines in $\proj{3}$ and $C_1$ is a proper secant
line of
$C_2$.  The union is an \ACM\ curve of degree 3.

\item \label{ACM probs} If $C_1$ and $C_2$ are allowed to have common
components
then observations 3 and 4 above fail.  In particular, even if $X$ is \ACM,
knowing something about $C_1$ and something about $X$ does not allow us to say
anything helpful about $C_2$.  See Example~\ref{bad acm}.

\end{enumerate}
\end{remark}

The amazing fact, which is the starting point of liaison theory, is that when
we restrict $X$ further by assuming that it is arithmetically Gorenstein, then
these problems can be overcome.  The following definition will be re-stated
in
more algebraic language later (Definition \ref{def of alg link}).

\begin{definition}
Let $C_1, C_2$ be  equidimensional subschemes of $\proj{n}$ having no common
component.  Assume that $X := C_1 \cup C_2$ is arithmetically Gorenstein.
Then
$C_1$ and $C_2$ are said to be {\em (directly) geometrically G-linked by $X$},
and we say that $C_2$ is {\em residual} to $C_1$ in $X$.  If $X$ is a complete
intersection, we say that $C_1$ and $C_2$ are {\em (directly)
geometrically CI-linked}.
\end{definition}

\begin{example}

If $X$ is the complete intersection in $\proj{3}$ of a surface consisting
of the
union of two planes with a surface consisting of one plane then $X$ links a
line $C_1$ to a different line $C_2$.

\begin{figure}[ht]
  \begin{picture}(500,80)
\put (255, 47){$C_1$}
\put (300,47){$C_2$}
  \put (70,20){\line (0,1){40}}
\put (70,60){\line (1,0){60}}
\put (130,20){\line (0,1){40}}
\put (70,20){\line (1,0){30}}
\put (115,5){\line(0,1){40}}
\put (100,20){\line (0,1){40}}
\put (85,75){\line (1,-1){30}}
\put (115,20){\line (1,0){15}}
\put (85,60){\line (0,1){15}}
\put (100,20){\line (1,-1){15}}
\put (140,30){$\cap$}
\put (145,10){\line (1,0){60}}
\put (145,10){\line (1,1){30}}
\put (175,40){\line (1,0){60}}
\put (205,10){\line (1,1){30}}
\put (245,30){$=$}
\put (270,10){\line (1,1){30}}
\put (270,40){\line (1,-1){30}}
  \end{picture}
\caption{Geometric Link} \label{geom-link}
\end{figure}
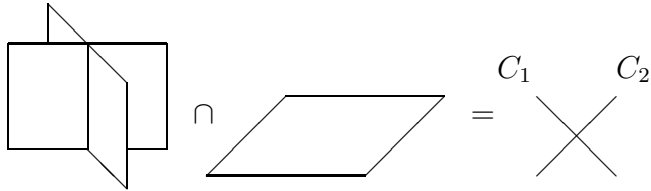

\end{example}

\begin{remark}
\begin{enumerate}
\item Given a scheme $C_1$, it is relatively easy (theoretically or on a
computer) to find a complete intersection $X$ containing $C_1$.  It is much
less easy to find one which gives a geometric link (see Example \ref{dble
line}).  In any case, $X$ is \ACM, and if one knows the degrees of the
generators
of
$I_X$ then one knows the degree and arithmetic genus of $X$ and even the
minimal free resolution of $I_X$, thanks to the Koszul resolution.

\item We will see that when $X$ is a complete intersection, a great deal of
information is passed from $C_1$ to $C_2$.  For example, $C_1$ is \ACM\ if and
only if $C_2$ is \ACM.  We saw above that this is not true when $X$ is merely
\ACM.  In fact, much stronger results hold, as we shall see.  An important
problem in general is to find {\em liaison invariants}.

\item While the notion of direct links has generated a theory, {\em liaison
theory}, that has become an active and fruitful area of study, it began as an
idea that did not quite work.  Originally, it was hoped that starting with {\em
any} curve $C_1$ in $\proj{3}$ one could always find a way to link it to a
``simpler'' curve $C_2$ (e.g.\ one of smaller degree), and use information
about $C_2$ to study $C_1$.  Based on a suggestion of Harris, Lazarsfeld
and Rao
\cite{LR} showed that this idea is fatally flawed: for a general curve $C
\subset \proj{3}$ of large degree, there is no simpler curve that can be
obtained from $C$ in {\em any} number of steps.

However, this actually led to a structure theorem for codimension two even
liaison classes \cite{BBM}, \cite{MP}, \cite{N-gorliaison}, \cite{nollet},
often called the {\em Lazarsfeld-Rao property}, which is one of the main
results
of liaison theory.
\end{enumerate}
\end{remark}

We now return to the question of how easy it is to find a complete intersection
containing a given scheme $C_1$ and providing a geometric link.  Since our
schemes are only assumed to be equidimensional, we will consider a non-reduced
example.

\begin{example}\label{dble line}
Let $C_1$ be a non-reduced scheme of degree two in $\proj{3}$, a so-called {\em
double line}.  It turns out (see e.g.\ \cite{dble}, \cite{montreal}) that the
homogeneous ideal of $C_1$ is of the form
\[
I_{C_1} = (x_0^2, x_0x_1, x_1^2, x_0 F(x_2,x_3) - x_1 G(x_2,x_3))
\]
where $F,G$ are homogeneous of the same degree, with no common factor.  Suppose
that $\deg F = \deg G = 100$.  Then it is easy to find complete intersections
$I_X$ whose generators have degree $\leq 100$; a simple example is $I_X =
(x_0^2, x_1^2)$.  However, any such complete intersection will have degree at
least 4 along the line $x_0 = x_1 = 0$, so it cannot provide a geometric link
for $C_1$: it is impossible to write $X = C_1 \cup C_2$ as schemes, no matter
what $C_2$ is.  However, once we look in degrees $\geq 101$, geometric links
are possible (since the fourth generator then enters the picture).
\end{example}

As this example illustrates, geometric links are too restrictive.  We have to
allow common components somehow.  However, an algebraic observation that we
made
above (Remark \ref{observations} (\ref{ideal quotient})) gives us the
solution.  That is, we will build our definition and theory around ideal
quotients.
Note first that if $X$ is merely \ACM, problems can arise, as mentioned in
Remark~\ref{observations} (\ref{ACM probs}).

\begin{example} \label{bad acm}
Let $I_X = (x_0,x_1)^2 \subset K[x_0,x_1,x_2,x_3]$, let $C_1$ be the double
line of Example~\ref{dble line} and let $C_2$ be the line defined by $I_{C_2} =
(x_0,x_1)$.  Then
\[
[I_X : I_{C_1}] = I_{C_2}, \ \ \ \hbox{but } \ [I_X : I_{C_2} ] = I_{C_2} \neq
I_{C_1}.
\]
\end{example}

As we will see, this sort of problem does not occur when our links are by {\em
arithmetically Gorenstein} schemes (e.g.\ complete intersections).  We make the
following definition.

\begin{definition} \label{def of alg link}
Let $C_1,C_2 \subset \proj{n}$ be subschemes with $X$ arithmetically
Gorenstein.  Assume that $I_X \subset I_{C_1} \cap I_{C_2}$ and that $[I_X :
I_{C_1}] = I_{C_2}$ and $[I_X : I_{C_2}] = I_{C_1}$.
Then  $C_1$ and $C_2$ are said to be {\em (directly) algebraically G-linked by
$X$}, and we say that $C_2$ is {\em residual} to $C_1$ in $X$.  We write $C_1
\stackrel{X}{\sim} C_2$.  If
$X$ is a complete intersection, we say that $C_1$ and $C_2$ are {\em (directly)
algebraically CI-linked}.  In either case, if $C_1 = C_2$ then we say that the
subscheme is {\em self-linked} by $X$.
\end{definition}

\begin{remark} \label{rem-alg-link}
An amazing fact, which we will prove later, is that when $X$ is arithmetically
Gorenstein (e.g.\ a complete intersection), then such a problem as
illustrated in Example \ref{bad acm} and Remark \ref{observations} (\ref{acm
not pres}) and (\ref{ACM probs}) does not arise.  That is, if
$I_X \subset I_{C_1}$ is arithmetically Gorenstein, and if we define
$I_{C_2} := [I_X : I_{C_1}]$ then it automatically follows that \linebreak
$[I_X
: I_{C_2}] = I_{C_1}$ whenever $C_1$ is equidimensional (i.e.\ $I_{C_1}$ is
unmixed).  It also follows that $\deg C_1 + \deg C_2 = \deg X$.

One might wonder what happens if $C_1$ is not equidimensional.  Then it turns
out that
\[
I_X : [I_X : I_{C_1}] = \hbox{top dimensional part of $C_1$},
\]
in other words this double ideal quotient is equal to the intersection of the
primary components of $I_{C_1}$ of minimal height (see \cite{migbook} Remark
5.2.5).
\end{remark}

\begin{example}
Let $I_X = (x_0x_1,x_0+x_1) = (x_0^2,x_0+x_1) = (x_1^2,x_0+x_1)$.  Let $I_{C_1}
= (x_0,x_1)$.  Then $I_{C_2} := [I_X : I_{C_1}] = I_{C_1}$.
\begin{figure}[ht]
  \begin{picture}(500,140)(-40,0)
  \put (70,20){\line (0,1){40}}
\put (70,60){\line (1,0){60}}
\put (130,20){\line (0,1){40}}
\put (70,20){\line (1,0){15}}
\put (85,5){\line (0,1){40}}
\put (100,20){\line (0,1){40}}
\put (85,5){\line (1,1){15}}
\put (85,45){\line (1,1){30}}
\put (100,20){\line (1,0){30}}
\put (115,60){\line (0,1){15}}
\put (140,30){$\cap$}
\put (155,35){\line (0,1){40}}
\put (155,75){\line (1,-1){30}}
\put (155,35){\line (1,-1){30}}
\put (185,5){\line (0,1){40}}
\put (195,30){$=$}
\put (215,30){$?$}
\put (135,90){\oval(70,60)[t]}
\put (99.5,90){\vector(0,-1){25}}
\put (170,90){\line (0,-1){10}}
  \end{picture}
\caption{Algebraic Link}\label{self-link fig}
\end{figure}
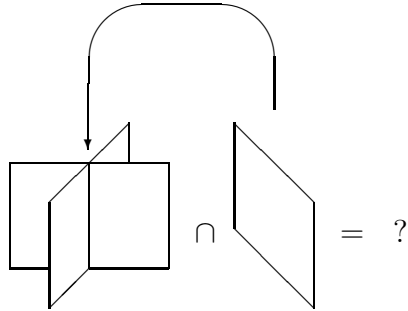
That is, $C_1$ is self-linked by $X$ (see Figure \ref{self-link fig}).
The question of when a scheme can be self-linked is a difficult one that has
been addressed by several papers, e.g.\  \cite{BE}, \cite{CC}, \cite{FKL},
\cite{KU}, \cite{dble}, \cite{rao3}.  Most schemes are not self-linked.  See
also Question \ref{stci question} of Section \ref{open prob sect}, and Example
\ref{hartshorne ex}.
\end{example}

Part of Definition \ref{def of alg link} is that the notion of direct linkage
is symmetric.  The observation above is that for most schemes it is not
reflexive (i.e.\ most schemes are not self-linked).  It is not hard to see that
it is rarely transitive.  Hence it is not, by itself, an equivalence relation.
Liaison is the equivalence relation {\em generated} by direct links, i.e.\ the
transitive closure of the direct links.

\begin{definition} Let $C \subset \proj{n}$ be an equidimensional subscheme.
The {\em Gorenstein liaison class of $C$} (or the {\em G-liaison class of $C$})
is the set of subschemes which can be obtained from $C$ in a finite number of
direct links.  That is, $C'$ is in the G-liaison class of $C$ if there exist
subschemes $C_1,\dots,C_r$ and arithmetically Gorenstein schemes
$X_1,\dots,X_r,X_{r+1}$ such that
\[
C \stackrel{X_1}{\sim} C_1 \stackrel{X_2}{\sim}  \dots \stackrel{X_r}{\sim}
  C_r \stackrel{X_{r+1}}{\sim} C'.
\]
If $r+1$ is even then we say that $C$ and $C'$ are {\em evenly G-linked}, and
the set of all subschemes that are evenly linked to $C$ is the {\em even
G-liaison class} of $C$.  If all the links are by complete intersections then
we talk about the {\em CI-liaison class of $C$} and the {\em even CI-liaison
class of $C$} respectively.  {\em Liaison} is the study of these equivalence
relations.
\end{definition}

\begin{remark}
Classically liaison was restricted to CI-links.  The most complete results have
been found in codimension two, especially for curves in $\proj{3}$ (\cite{BBM},
\cite{MP}, \cite{rao1}, \cite{rao2}, \cite{N-gorliaison}, \cite{nollet}).
However, Schenzel \cite{schenzel} and later Nagel \cite{N-gorliaison} showed
that the set-up and basic results for complete intersections continue to hold
for G-liaison as well, in any codimension.

As we noted earlier, in codimension two every arithmetically Gorenstein scheme
is a complete intersection.  Hence the complete picture which is known in
codimension two belongs just as much to Gorenstein liaison theory as it does to
complete intersection liaison theory!

The recent monograph \cite{KMMNP} began the study of the important differences
that arise, and led to the recent focus on G-liaison in the literature.  We
will describe much of this work.  In particular, we will see how several
results in G-liaison theory neatly generalize standard results in codimension
two theory, while the corresponding statements for CI-liaison are false!
\end{remark}

Here are some natural questions about this equivalence class, which we will
discuss and answer (to the extent possible, or known) in these lectures.  In
the last section we will discuss several open questions.  We will see that the
known results very often hold for {\em even} liaison classes, so some of our
questions focus on this case.

\begin{enumerate}
\item \label{nec cond} Find necessary conditions for $C_1$ and $C_2$ to be in
the same (even) liaison class (i.e.\ find {\em (even) liaison invariants}).  We
will see that the dimension is invariant, the property of being \ACM\ is
invariant, as is the property of being locally Cohen-Macaulay, and that more
generally, for an even liaison class the graded modules
$H^i_* ({\mathcal I}_C)$ are essentially invariant (modulo shifts), for $1
\leq i \leq \dim C$.  The situation is somewhat simpler when we assume that
the schemes are locally Cohen-Macaulay.  There is also a condition in terms of
stable equivalence classes of certain reflexive sheaves.

\item Find sufficient conditions for $C_1$ and $C_2$ to be in the same (even)
liaison class.  We will see that for instance being linearly equivalent is a
sufficient condition for even liaison, and that for codimension two the problem
is solved.  In particular, for codimension two there is a condition which is
both necessary and sufficient for two schemes to be in the same even liaison
class.  An important question is to find a condition which is both necessary
and sufficient in higher codimension, either for CI-liaison or for G-liaison.
Some partial results in this direction will be discussed.

\item Is there a structure common to all even liaison classes?  Again, this is
known in codimension two.  It is clear that the structure, as it is commonly
stated in codimension two, does not hold for even G-liaison.  But perhaps some
weaker structure does hold.

\item Are there good applications of liaison?  In codimension two we will
mention a number of applications that have been given in the literature, but
there are fewer known in higher codimension.

\item \label{diff and sim} What are the differences and similarities between
G-liaison and CI-liaison?  What are the advantages and disadvantages of either
one?  See Remark \ref{some nice results}  and Section \ref{compare section}.

\item Do geometric links generate the same equivalence relation as algebraic
links?  For CI-liaison the answer is ``no'' if we allow schemes that are not
local complete intersections.  Is the answer ``yes'' if we restrict to local
complete intersections?  And is the answer ``yes'' in any case for G-liaison?

\item \label{geom ACM liaison}  We have seen that there are fewer nice
properties when we try to allow links by \ACM\ schemes.  It is possible to
define an equivalence relation using ``geometric ACM links."  What does this
equivalence relation look like?  See Remark \ref{Walter result}.
\end{enumerate}

\begin{remark} \label{Walter result}
We now describe the answer to Question \ref{geom ACM liaison} above.  Clearly
if we are going to study {\em geometric} ACM links, we have to restrict to
schemes that are locally Cohen-Macaulay in addition to being equidimensional.
Then we quote
the following three results:
\begin{itemize}
\item (\cite{walter}) Any locally Cohen-Macaulay equidimensional
subscheme $C \subset \proj{n}$ is ACM-linked in finitely many steps to some
\ACM\ scheme.

\item (\cite{KMMNP} Remark 2.11) Any \ACM\ scheme is CM-linked to a complete
intersection of the same dimension.

\item (Classical; see \cite{schwartau}) Any two complete intersections of the
same dimension are CI-linked in finitely many steps.  (See Open Question
\ref{gor g-linked?} on page \pageref{gor g-linked?} for an interesting related
question for G-liaison.)
\end{itemize}
The first of these is the deepest result.  Together they show that there is
only one ACM-liaison class, so there is not much to study here.  Walter
\cite{walter} does give a bound on the number of steps needed to pass from an
arbitrary locally Cohen-Macaulay scheme to an \ACM\ scheme, in terms of the
dimension.  In particular, for curves it can be done in one step!
\end{remark}

So the most general kind of linkage for subschemes of projective space seems to
be Gorenstein liaison.  Recent contributions to this theory have been made by
Casanellas, Hartshorne, Kleppe, Lesperance, Migliore, Mir\'o-Roig, Nagel,
Notari, Peterson, Spreafico, and others.  We will describe this work in the
coming sections.

\begin{remark} \label{some nice results}
To end this section, as a partial answer to Question \ref{diff and sim}, we
would like to mention two results about G-liaison from \cite{KMMNP} that are
easy to state, cleanly generalize the codimension two case, and are {\em false}
for CI-liaison.

\begin{itemize}
\item Let $S \subset \proj{n}$ be \ACM\ satisfying property $G_1$ (so that
linear equivalence is well-defined; see \cite{gen-div}).  Let $C_1 ,C_2
\subset S$ be divisors such that $C_2 \in | C_1 + tH|$, where $H$ is the
class of a hyperplane section and $t \in {\mathbb Z}$.  Then $C_1$ and $C_2$
are G-linked in two steps.

\item Let $V \subset \proj{n}$ be a subscheme of codimension $c$ such that
$I_V$ is the ideal of maximal minors of a $t \times (t+c-1)$ homogeneous
matrix.  Then $V$ can be G-linked to a complete intersection in finitely many
steps.
\end{itemize}

\end{remark}


\section{Preliminary results} \label{prelim result section}

The purpose of this section is to recall some concepts and results
we will use later on. Among them we include a comparison of local
and sheaf cohomology, geometric and algebraic hyperplane
sections, local duality and $k$-syzygies. Furthermore, we discuss
the structure of deficiency modules and introduce the notion of
(cohomological) minimal shift.
\smallskip

Throughout we will use the following notation. $A$ will always
denote a (standard)
  graded $K$-algebra, i.e.\ $A = \oplus_{i \geq 0} [A]_i$ is generated (as
algebra)
  by its elements of degree $1$, $[A]_0 = K$ is a field and $[A]_i$ is the
vector
  space of elements of degree $i$ in $A$. Thus, there is a homogeneous ideal $I
  \subset R = K[x_0,\ldots,x_n]$ such that $A \cong R/I$. The
irrelevant maximal ideal  of $A$ is $\fm := \fm_A :=  \oplus_{i >
0} [A]_i$.

If $M$ is a graded module over the ring $A$ it is always assumed
that $M$ is $\mathbb{Z}$-graded and $A$ is a graded $K$-algebra as
above. All $A$-modules will be finitely generated unless stated
otherwise. Furthermore, it is always understood that homomorphisms
between graded $R$-modules are morphisms in the category of graded
$R$-modules, i.e.\ are graded of degree zero.
\bigskip

\noindent {\it Local cohomology}
\medskip

There will be various instances where it is preferable to use
local cohomology instead of the (possibly more familiar) sheaf
cohomology. Thus we recall the definition of local cohomology and
describe the comparison between both cohomologies briefly.

We start with the following

\begin{definition}
Let $M$ be an arbitrary $A$ module. Then we set $$ \HH^0(M) := \{m
\in M \s \fm_A^k \cdot  m = 0 \; \mbox{for some} \ k \in \NN \}.
$$
\end{definition}

This construction provides the functor $\HH^0(\_)$ from the
category of $A$-modules into itself. It has the following
properties.

\begin{lemma} \mbox{}
\begin{itemize}
\item[(a)] The functor $\HH^0(\_)$ is left-exact.
\item[(b)] $\HH^0(M)$ is an Artinian module.
\item[(c)] If $M$ is graded then $\HH^0(M)$ is graded as well.
\end{itemize}
\end{lemma}

\begin{example}
Let $I \subset R$ be an ideal with saturation $I^{sat} \subset R$
then $$ \HH^0(R/I) = I^{sat}/I. $$ This is left as an exercise to
the reader.
\end{example}

Since the functor $\HH^0(\_)$ is left-exact one can define its
right-derived functors using injective resolutions.

\begin{definition} \label{def-local-coh}
The $i$-th right derived functor of $\HH^0(\_)$ is called the {\it
$i$-th local cohomology functor} and denoted by $\HH^i(\_)$.
\end{definition}

Thus, to each short exact sequence of $A$-modules $$ 0 \to M' \to
M \to M'' \to 0 $$ we have a the induced long exact cohomology
sequence $$ 0 \to \HH^0(M') \to \HH^0(M) \to \HH^0(M'') \to
\HH^1(M') \to \ldots $$ We note some further properties.

\begin{lemma} \label{lem-prop-coho} \mbox{}
\begin{itemize}
\item[(a)] All $\HH^i(M)$ are Artinian $A$-modules (but often not finitely
generated).
\item[(b)] If $M$ is graded then all $\HH^i(M)$ are graded as well.
\item[(c)] The Krull dimension and the depth of $M$ are cohomologically
characterized by
\begin{eqnarray*}
  \dim M & = & \max \{i \s \HH^i(M) \neq 0 \} \\
\depth M & =  & \min \{i \s \HH^i(M) \neq 0 \}
\end{eqnarray*}
\end{itemize}
\end{lemma}

Slightly more than stated in part (b) is true: The cohomology
sequence associated to a short exact sequence of graded modules is
an exact sequence of graded modules as well.

Part (a) implies that a local cohomology module is Noetherian if
and only if it has finite length.
Part (c) immediately provides the following.

\begin{corollary} \label{cohen-macaulay}
The module $M$ is \CM\ if and only if $\HH^i(M) = 0$ for all $i
\neq \dim M$.
\end{corollary}

As mentioned in the last section, a subscheme $X \subset \PP^n$ is called {\em
\aCM\ } if its homogeneous coordinate ring $R/I_X$ is \CM, i.e.\ a \CM-module
over itself.
\smallskip

Now we want to relate local cohomology to sheaf cohomology.

The projective spectrum $ X = \Proj A$ of a graded $K$-algebra $A$
is a projective scheme of dimension $(\dim A - 1)$. Let $\cF$ be a
sheaf of modules over $X$. Its cohomology modules are denoted by
$$ H^i_*(X,\cF) = \bigoplus_{j \in \mathbb{Z}} H^i(X,\cF(j)). $$ If
there is no ambiguity on the scheme $X$ we simply write
$H^i_*(\cF)$.

There are two functors relating graded $A$-modules and sheaves of
modules over $X$. One is the ``sheafification'' functor which
associates to each graded $A$-module $M$ the sheaf $\tilde{M}$.
This functor is exact.

In the opposite direction there is the ``twisted global sections''
functor which associates to each sheaf $\cF$ of modules over $X$
the graded $A$-module $H^0_*(X,\cF)$. This functor is only left
exact. If $\cF$ is quasi-coherent then the sheaf
$\widetilde{H^0_*(X,\cF)}$ is canonically isomorphic to $\cF$.
However, if $M$ is a graded $A$-module then the module
$H^0_*(X,\tilde{M})$ is not isomorphic to $M$ in general. In fact,
even if $M$ is finitely generated, $H^0_*(X,\tilde{M})$ needs not
to be finitely generated. Thus the functors  $\tilde{\_\_}$  and
$H^0_*(X,\_\_)$ do not establish an equivalence of categories
between graded $A$-modules and quasi-coherent sheaves of modules
over $X$. However, there is the following comparison result (cf.\
\cite{SV2}).

\begin{proposition} \label{comp} Let $M$ be a graded $A$-module. Then there
is an
exact sequence $$ 0 \to \HH^0(M) \to M \to H^0_*(X,\tilde{M}) \to
\HH^1(M) \to 0 $$ and for all $i \geq 1$ there are isomorphisms $$
H^i_*(X,\tilde{M}) \cong \HH^{i+1}(M). $$
\end{proposition}

The result is derived from the exact sequence $$ 0 \to \HH^0(M)
\to M \to H^0(M) \to \HH^1(M) \to 0 $$ where $ \displaystyle H^0(M) =
\lim_{\stackrel{\longrightarrow}{{\scriptstyle n}}}
\Hom_R(\fm^n,M)$. Note that $H^0(M) \cong H^0_*(X,\tilde{M})$.
\smallskip

\begin{corollary} \label{cor-defic-mod}
Let $X \subset \PP^n = \Proj R$ be a closed subscheme of dimension $d \leq
n-1$. Then there are graded isomorphisms $$ H^i_*(\cI_X) \cong
\HH^i(R/I_X) \quad \mbox{for all} \; i = 1,\ldots,d+1. $$
\end{corollary}

\begin{proof}
Since $\HH^i(R) = 0$ if $i \leq n$ the cohomology sequence of $$ 0
\to I_X \to R \to R/I_X \to 0 $$ implies $\HH^i(R/I_X) \cong
\HH^{i+1}(I_X)$ for all $i < n$. Thus, the last proposition yields
the claim.
\end{proof}

\begin{remark} \label{rem-regularity}
Let $M$ be a graded $R$-module. Then its {\it Castelnuovo-Mumford
   regularity} is the number
$$
\reg M := \min \{m \in \ZZ \s [\HH^i(M)]_{j-i} = 0 \; \; \mbox{for all} \; j >
m \}.
$$
For a subscheme $X \subset \PP^n$ we put $\reg \cI_X = \reg I_X$. The
preceding corollary shows that this last definition agrees with Mumford's in
\cite{mumford-book}.
\end{remark}

It is convenient and common to use the following names.

\begin{definition}
Let $X \subset \PP^n$ be a closed subscheme of dimension $d$. Then
the graded $R$-modules $H^i_*(\cI_X),\ i=1,\ldots,d,$ are called
the {\it deficiency modules} of $X$. If $X$ is $1$-dimensional
then $H^1_*(\cI_X)$ is also called the {\em \HR\ } of $X$.
\end{definition}

The deficiency modules reflect properties of the scheme.  For example, as
mentioned in the first section, it follows from what we have now said
(Corollary \ref{cohen-macaulay} and Corollary \ref{cor-defic-mod}) that
$X$ is \acm\ if and only if $H^i_*({\mathcal I}_X) = 0$ for $1 \leq i \leq
\dim X$. Note that a scheme
$X
\subset
\PP^n$ is said to be {\it equidimensional} if its homogeneous ideal $I_X
\subset R$ is unmixed, i.e.\ if all its components have the same dimension.
In particular, an equidimensional scheme has no embedded components.

\begin{lemma}\label{lem-unmixed-crit} For a subscheme $X \subset \PP^n$ we have
\begin{itemize}
\item[(a)] $X$ is equidimensional and locally \CM\ if and only if all its
deficiency
modules have finite length.
\item[(b)] $X$ is equidimensional if and only if \; $\dim R/ \Ann
H^i_*(\cI_X) \leq
i-1$ for all $i = 1,\ldots,\dim X$.
\end{itemize}
\end{lemma}

By a curve we always mean an equidimensional scheme of dimension
$1$. In particular, a curve is \lCM\ since by definition it does
not have embedded components.
Thus, we have.

\begin{corollary}
A $1$-dimensional scheme $X \subset \PP^n$ is a curve if and only
if its \HR\ $H^1_*(\cI_X)$ has finite length.
\end{corollary}
\bigskip

\noindent {\it Hyperplane sections}
\medskip

Let $H \subset \PP^n$ be the hyperplane defined by the linear form
$l \in R$. The {\em geometric hyperplane section} (or simply the {\em
hyperplane
section}) of a scheme $X \subset \PP^n$ is the subscheme  $X \cap
H$.  We usually consider $X \cap H$ as a subscheme of
$H \cong \PP^{n-1}$, i.e.\ its homogeneous ideal $I_{X \cap H}$ is
an ideal of $\bar{R} = R/l R$. The {\em algebraic hyperplane section} of
$X$ is given by the ideal $\overline{I_X} := (I_X + l R)/l R
\subset \bar{R}$. $\overline{I_X}$ is not necessarily a saturated
ideal. In fact, the saturation of $\overline{I_X}$ is just $I_{X
\cap H}$.  The difference between the
hyperplane section and the algebraic hyperplane section is measured by
cohomology.

\begin{lemma} \label{hyperpl sect lemma}
$$ \HH^0(R/I_X + l R) \cong I_{X \cap H}/\overline{I_X} $$
\end{lemma}

If the ground field $K$ contains sufficiently many elements we can
always find a hyperplane which is general enough with respect to a
given scheme $X$. In particular, we get $\dim X \cap H = \dim X -
1$ if $X$ has positive dimension. In order to relate properties of
$X$ to the ones of its hyperplane section we note some useful
facts. We use the following notation.

For a graded $A$-module $M$ we denote by $h_M$ and $p_M$ its
Hilbert function and Hilbert polynomial, respectively, where $h_M
(j) = \rank [M]_j$. The Hilbert function and Hilbert polynomial of
a subscheme $X \subset \PP^n$ are the corresponding functions of
its homogeneous coordinate ring $R/I_X$. For a numerical function
$h: \ZZ \to \ZZ$ we define its first difference by $\Delta h (j) =
h(j) - h(j-1)$ and the higher differences by $\Delta^i h = \Delta
(\Delta^{i-1} h)$ and $\Delta^0 h = h$.

\begin{remark} \label{rem-hyp-sect}
Suppose $K$ is an infinite field and let $H \subset \PP^n$ be a
hyperplane.

(i) If $\dim X > 0$ and $H$ is general enough then we have $$ I_{X
\cap H}= \overline{I_X} \quad \mbox{if and only if} \quad H^1_*
(\cI_X) = 0. $$

(ii) If $X  \subset \PP^n$ is locally or arithmetically
Cohen-Macaulay of positive dimension then $X \cap H$ has the same
property for a general hyperplane $H$. The converse is false in
general.

(iii) Suppose $X  \subset \PP^n$ is \aCM\ of dimension $d$. Let \linebreak
$l_1, \ldots ,l_{d+1} \in R$ be linear forms such that $\bar{A} =
R/(I_X + (l_1,\ldots,l_{d+1}))$ has dimension zero. Then $\bar{A}$
is called an {\it Artinian reduction} of $R/I_X$. For its Hilbert
function we have $h_{\bar{A}} = \Delta^{d+1} h_{R/I_X}$.
\end{remark}
\bigskip

\noindent {\it Minimal free resolutions}
\medskip

Let $R = K[x_0,\ldots,x_n]$ be the polynomial ring. By our
standard conventions a homomorphism $\ffi: M \to N$ of graded
$R$-modules is graded of degree zero, i.e.\ $\ffi([M]_j) \subset
[N]_j$ for all integers $j$. Thus, we have to use degree shifts
when we consider the homomorphism $R(-i) \to R$ given by
multiplication by $x_0^i$. Observe that $R(-i)$ is not a graded
$K$-algebra unless $i = 0$.

\begin{definition} \label{k-syz} Let $M$ be a graded $R$-module. Then
$N \neq 0 $ is said to be a {\itshape $k$-syzygy of $M$} (as
$R$-module)
  if there is an exact sequence of graded $R$-modules
$$ 0 \to N \to F_k \stackrel{\ffi_k}{\longrightarrow} F_{k-1} \to
\ldots \to F_1 \stackrel{\ffi_1}{\longrightarrow} M \to 0 $$ where
the modules $F_i, i = 1, \ldots , k,$ are free $R$-modules.
A module is called a {\itshape $k$-syzygy} if it is a $k$-syzygy
of some module.
\end{definition}

Note that a $(k+1)$-syzygy is also a $k$-syzygy (not for the same module).
Moreover,
every
$k$-syzygy $N$ is a maximal $R$-module, i.e.\ $\dim N = \dim R$.

Chopping long exact sequences into short ones we easily obtain

\begin{lemma} \label{lem-coho-of-syzy}
If  $N$ is a $k$-syzygy of the $R$-module $M$ then $$ \HH^i(N)
\cong \HH^{i-k}(M) \quad \mbox{for all} \; i < \dim R. $$
\end{lemma}

If follows that the depth of a $k$-syzygy is at least $k$.

The next concept ensures uniqueness properties.

\begin{definition}
Let $\ffi: F \to M$ be a homomorphism of $R$-modules where $F$ is
free. Then $\ffi$ is said to be a {\itshape minimal}  homomorphism
if $\ffi \otimes id_{R/\fm} : F/\fm F \to M/\fm M$ is the zero map
in case $M$ is free and an isomorphism in case $\ffi$ is
surjective.

In the situation of the definition above, $N$ is said to be a
{\itshape
   minimal $k$-syzygy} of $M$ if the morphisms $\ffi_i, i = 1, \ldots , k$,
are minimal.  If  $N$ happens to be free then the exact sequence
is called a {\it minimal free resolution} of $M$.
\end{definition}

Nakayama's lemma implies easily that minimal $k$-syzygies of $M$
are unique up to isomorphism and that a minimal free resolution is
unique up to isomorphism of complexes.

Note that every finitely generated projective $R$-module is free.

\begin{remark} Let
$$ 0 \to F_s \stackrel{\ffi_s}{\longrightarrow} F_{s-1} \to \ldots
\to F_1 \stackrel{\ffi_1}{\longrightarrow} F_0 \to M \to 0 $$ be a
free resolution of $M$. Then it is minimal if and only if (after
choosing bases for $F_0,\ldots,F_s$) the matrices representing
$\ffi_1,\ldots,\ffi_s$ have entries in the maximal ideal $\fm =
(x_0,\ldots,x_n)$ only.
\end{remark}
\bigskip

\noindent {\it Duality results}
\medskip

Later on we will often use some duality results. Here we state
them only for the
  polynomial ring $R = K[x_0, \dots , x_n]$. However, they are
  true, suitably adapted, over any graded Gorenstein $K$-algebra.

Let $M$ be a graded $R$-module. Then we will consider two types of
dual modules, the $R$-dual $M^* := \Hom_R (M, R)$ and the $K$-dual
$M^{\vee} := \oplus_{j \in \ZZ} \Hom_K ([M]_{-j}, K)$.

Now we can state a version of Serre duality (cf.\ \cite{S},
\cite{SV2}).

\begin{proposition} \label{prop-serredual} Let $M$ be a graded  $R$-module.
  Then for all $i \in \mathbb{Z}$, we have   natural isomorphisms of
   graded $R$-modules
$$ H^i_{{\mathfrak m}} (M)^{\vee} \cong \Ext^{n+1-i}_R
(M,R)(-n-1). $$
\end{proposition}

The $K$-dual of the top cohomology module plays a particular role.

\begin{definition} \label{def-can-module}
The module $K_M := \Ext^{n+1-\dim M}_R (M,R)(-n-1)$ is called the
{\it canonical module} of $M$. The canonical module $K_X$ of a
subscheme $X \subset \PP^n$ is defined as $K_{R/I_X}$.
\end{definition}

\begin{remark} \label{rem-can-mod}
(i) For a subscheme $X \subset \PP^n$ the sheaf $\omega_X :=
\widetilde{K_X}$ is the dualizing sheaf of $X$.

(ii) If $X \subset \PP^n$ is \aCM\ with minimal free resolution $$
0 \to F_c \stackrel{\ffi_c}{\longrightarrow} F_{c-1} \to \ldots
\to F_1 \to I_X \to 0 $$ then dualizing with respect to $R$
provides the complex $$ 0 \to R \to F^*_1 \to \ldots \to F^*_{c-1}
\stackrel{\ffi^*_c}{\longrightarrow} F^*_c \to \coker \ffi^*_c \to
0 $$ which is a minimal free resolution of $\coker \ffi^*_c \cong
K_X (n+1)$.
\end{remark}

If the scheme $X$ is equidimensional and \lCM, one can relate the
cohomology modules of $X$ and its canonical module. More
generally, we have
  (\cite{S},
Korollar 3.1.3).

\begin{proposition} \label{prop-duality} Let $M$ be a graded $R$-module
such that $\HH^i(M)$ has finite length if $i \neq
   d = \dim M$. Then there are canonical isomorphisms for $i = 2, \ldots ,
   d-1$
$$ \HH^{d+1-i}(K_M) \cong \HH^i(M)^{\vee}. $$
\end{proposition}

Observe that the first cohomology $\HH^1(M)$ is not involved in
the statement above.
\bigskip

\noindent {\it Restrictions for deficiency modules}
\medskip

Roughly speaking, it will turn out that there are no restrictions
on the module structure of deficiency modules, but there are
restrictions on the degrees where non-vanishing pieces can occur.
\smallskip

In the following result we will assume $c \leq n-1$ because
subschemes of $\PP^n$ with codimension $n$ are \aCM.

\begin{proposition} \label{prop-ex-lCM}
Suppose the ground field $K$ is infinite. Let $c$ be
an integer with $2 \leq c \leq n-1$ and let $M_1,\ldots,M_{n-c}$ be
graded $R$-modules of finite length. Then there is an integral
locally Cohen-Macaulay subscheme $X \subset \PP^n$ of codimension
$c$ such that $$ H^i_*(\cI_X) \cong M_i (-t) \quad \mbox{for all}
\; i = 1,\ldots,n-c $$ for some integer $t$.
\end{proposition}

\begin{proof}
Choose a smooth complete intersection $V \subset \PP^n$ such that
$$ I_V = (f_1,\ldots,f_{c-2}) \subset \bigcap_{i=1}^{n-c} \Ann M_i
$$ where $I_V = 0$ if $c = 2$.

Let $N_i$ denote a $(i+1)$-syzygy of $M_i$ as $R/I_V$-module and
let $r$ be the rank of $N = \oplus_{i=1}^{n-c} N_i$. For $s \ll 0$
the cokernel of a general map $\ffi: R^{r-1} \to N$ is
torsion-free of rank one, i.e.\ isomorphic to $I(t)$ for some
integer $t$ where $I \subset A = R/I_V$ is an ideal such that
$\dim A/I = \dim A -2$. Moreover, $I$ is a prime ideal by
Bertini's theorem. The preimage of $I$ under the canonical
epimorphism $R \to A$ is the defining ideal of a subscheme $X
\subset \PP^n$ having the required properties. For details we
refer to \cite{MNP2}.
\end{proof}

\begin{remark} \label{rem-prescr-coho}
(i) The previous result can be generalized as follows. Let
$M_1,\ldots,M_{n-c}$ be graded (not necessarily finitely
generated) $R$-modules such that $M_i^{\vee}$ is finitely
generated of dimension $\leq i-1$ for all $i = 1,\ldots,n-c$. Then
there is an equidimensional subscheme $X \subset \PP^n$ of
codimension $c$ such that $$ H^i_*(\cI_X) \cong M_i (-t) \quad
\mbox{for all} \; i = 1,\ldots,n-c $$ for some integer $t$. Details will
appear in \cite{N-Eil-McLane}.  Note
that the condition on the modules $M_1,\ldots,M_{n-c}$ is
necessary according to Lemma \ref{lem-unmixed-crit}.

(ii) A more general version of Proposition \ref{prop-ex-lCM} for
subschemes of codimension two is shown in \cite{EG}.
\end{remark}

Now we want to consider the question of which numbers $t$ can occur
in Proposition \ref{prop-ex-lCM}. The next result implies that
with $t$ also $t+1$ occurs. The name of the statement will be
explained later on.

\begin{lemma}[Basic double link]  \label{lem-basic-d-link}
Let $0 \neq J \subset I \subset R$ be homogeneous ideals such that
$\codim I = \codim J + 1$ and $R/J$ is Cohen-Macaulay. Let $f \in
R$ be a homogeneous element of degree $d$ such that $J : f = J$.
Then the ideal $\tilde{I} := J + f I$ satisfies $\codim \tilde{I}
= \codim I$ and $$ \HH^i(R/\tilde{I}) \cong \HH^i(R/I) (-d) \quad
\mbox{for all}\; i < \dim R/I. $$ In particular, $I$ is unmixed if
and only if $\tilde{I}$ is unmixed.
\end{lemma}

\begin{proof}
Consider the sequence
\begin{equation} \label{bdl exact sequence}
  0 \to J(-d)
\stackrel{\ffi}{\longrightarrow} J \oplus I(-d)
\stackrel{\psi}{\longrightarrow} \tilde{I} \to 0
\end{equation}
  where $\ffi$
and $\psi$ are defined by $\ffi(j) = (f j, j)$ and $\psi(j, i) = j
- f i$. It is easy to check that this sequence is exact. Its
cohomology sequence implies the claim on the dimension and
cohomology of $R/\tilde{I}$. The last claim follows by Lemma
\ref{lem-unmixed-crit}.
\end{proof}

\begin{proposition} \label{prop-coho-min-shift}
Suppose that $K$ is infinite. Let $M_{\bullet} =
(M_1,\ldots,M_{n-c})$ $(2 \leq c < n)$ be a vector of graded (not
necessarily finitely generated) $R$-modules such that $M_i^{\vee}$
is finitely generated of dimension $\leq i-1$ for all $i =
1,\ldots,n-c$ and not all of these modules are trivial. Then there
is an integer $t_0$ such that there is an equidimensional
subscheme $X \subset \PP^n$ of codimension $c$ with $$
H^i_*(\cI_X) \cong M_i (-t) \quad \mbox{for all} \; i =
1,\ldots,n-c $$ for some integer $t$ if and only if $t \geq t_0$.
\end{proposition}

\begin{proof}
If the ground field $K$ is infinite we can choose the element $f$
in Lemma \ref{lem-basic-d-link} as a linear form. Thus, in spite
of this lemma and Remark \ref{rem-prescr-coho} it suffices to show
that $$ H^i_*(\cI_X) \cong M_i (-t) \quad \mbox{for all} \; i =
1,\ldots,n-c $$ is impossible for a subscheme $X \subset \PP^n$ of
codimension $c$ if $t \ll 0$. But this follows if $\dim X = 1$
because we have for every curve $C \subset \PP^n$
\begin{equation} \label{incr-in-neg-deg}
  h^1(\cI_C (j-1)) \leq \max \{0, h^1(\cI_C(j)) - 1\} \quad \mbox{if} \; j
\leq 0
\end{equation}
  by \cite{BN}, Lemma 3.4 or \cite{Mi-invar}. By taking
general hyperplane sections of $X$, the general case is easily
reduced to the case of curves.  See also Proposition 1.4 of \cite{BM4}.
\end{proof}

The last result allows us to make the following definition.

\begin{definition} \label{def-coho-min-shift}
The integer $t_0$, which by Proposition \ref{prop-coho-min-shift}
is uniquely determined, is called the {\it (cohomological)
minimal shift} of
$M_{\bullet}$.
\end{definition}

\begin{example} \label{ex-min-shift}
Let $M_{\bullet} = (K)$. Then the estimate (\ref{incr-in-neg-deg}) for the
first cohomology of a curve in the last proof shows that the minimal
shift $t_0$ of $M_{\bullet}$ must be non-negative. Since we have
for a pair $C$ of skew lines $H^1_*(\cI_C) \cong K$ we obtain $t_0
= 0$ as minimal shift of $(K)$.
\end{example}


\section{Gorenstein Ideals}

Before we can begin the discussion of Gorenstein liaison, we will need some
basic facts about Gorenstein ideals and Gorenstein algebras.  In this
section we
will give the  definitions, properties, constructions, examples and
applications
which will be used or discussed in the coming sections.  Most of the material
discussed here is treated in more detail in \cite{migbook}.

We saw in Remark \ref{rem-hyp-sect} that if $X$ is \ACM\ of dimension $d$
with coordinate ring $A = R/I_X$ then we have the {\em Artinian reduction}
$\bar A$ of $X$ (or of $R/I_X$).  Its Hilbert function was given as
$h_{\bar{A}} = \Delta^{d+1} h_{R/I_X}$.  Since $\bar A$ is finite dimensional
as a $K$-vector space, we have that $h_{\bar A}$ is a finite sequence of
integers
\[
1 \ \ c \ \ h_2 \ \ h_3 \ \ \dots \ \ h_s \ \ 0 \ \ \dots
\]
This sequence is called the {\em $h$-vector} of $X$, or of $A$.  In
particular, $c$ is the {\em embedding codimension} of $X$.  In other words,
$c$ is the codimension of $X$ inside the smallest linear space containing
it.  Of course, the Hilbert function of $X$ can be recovered from the
$h$-vector by ``integrating.''

Now suppose that $X$ is \acm\ and non-degenerate in $\proj{n}$, of
codimension $c$, and that $R/I_X$ has minimal free resolution
\[
0 \rightarrow F_c \rightarrow F_{c-1} \rightarrow \dots \rightarrow F_1
\rightarrow R \rightarrow R/I_X \rightarrow 0.
\]
Suppose that $F_c = \bigoplus_{i=1}^r R(-a_i)$ and let $a = \max_i \{ a_i
\}$.  As mentioned in Section 2, $r = \rank F_c$ is called the {\em
Cohen-Macaulay type} of $X$ (or of $A$).  Furthermore, we have the relation
\begin{equation} \label{relation betw a, c and t}
a -c = s = \reg {\mathcal I}_X -1
\end{equation}
where $s$ is the last degree in which the $h$-vector is non-zero and $\reg
{\mathcal I}_X$ is the Castelnuovo-Mumford regularity of ${\mathcal I}_X$
(cf.\ Remark \ref{rem-regularity}). We now formally make the definition
referred to in Section 2:

\begin{definition}
The subscheme $X \subset \proj{n}$ is {\em arithmetically Gorenstein} if it is
\acm\ of Cohen-Macaulay type 1.  We often say that $I_X$ is Gorenstein or $I_X$
is arithmetically Gorenstein.
\end{definition}

\begin{example}
A line in $\proj{3}$ is arithmetically Gorenstein since its minimal free
resolution is
\[
0 \rightarrow R(-2) \rightarrow R(-1)^2 \rightarrow I_X \rightarrow 0,
\]
and $R(-2)$ has rank 1.  More generally, any complete intersection in
$\proj{n}$
is arithmetically Gorenstein thanks to the Koszul resolution.  The last free
module in the resolution of the complete intersection of forms of degree
$d_1,\dots,d_c$ is $R(-d_1 - \dots - d_c)$.
\end{example}

\begin{remark}
In Remark \ref{rem-hyp-sect} (ii) it was noted that if $X$ is arithmetically
Cohen-Macaulay of dimension $\geq 1$ then the general hyperplane section $X
\cap
H$ is also \acm.  (In fact this is true for any proper hyperplane section.)  It
was remarked that the converse is false in general.  However, there are
situations in which the converse does hold.

First, if $X$ is assumed to be equidimensional (i.e.\ $I_X$ is unmixed) and
locally Cohen-Macaulay of dimension $\geq 2$ then it is not hard to show
that the
converse holds.  Indeed, let $L$ be a general linear form defining the
hyperplane $H$ and consider the exact sequence
\[
H^i ({\mathcal I}_X (s-1)) \stackrel{\times L}{\longrightarrow}
H^i({\mathcal I}_X (s)) \rightarrow H^i ({\mathcal I}_{X \cap H|H} (s))
\rightarrow H^{i+1}({\mathcal I}_X (s-1)) \stackrel{\times L}{\longrightarrow}
H^{i+1}({\mathcal I}_X(s))
\]
If $X\cap H$ is \acm\ and $1 \leq i \leq \dim X\cap H = \dim X -1$ then the
multiplication map on the left is surjective for all $s$ and the one on the
right is injective for all $s$.  Both of these are impossible unless $X$ is
itself \acm, because $H^i ({\mathcal I}_X)$ has finite length for $1 \leq i
\leq
\dim X$ (cf.\ Lemma \ref{lem-unmixed-crit}).

Obviously if $X$ is the union of an \acm\ scheme and a point (possibly
embedded)
then it is not \acm\ but its general hyperplane section is \acm.  Also,
clearly if $X$ is a curve which is not \acm\ then its general hyperplane
section
is
\acm\ since it is a finite set of points, but again, $X$ is not \acm.  A
fascinating question is whether there are conditions on $X \cap H$ which {\em
force} a curve $X$ to be \acm.  The best results in this direction come when $X
\cap H$ is arithmetically Gorenstein.  Several authors have contributed to this
question, but we mention in particular \cite{HU} and \cite{strano}.
\end{remark}

There are several other  conditions which are equivalent to being \acm\ with
Cohen-Macaulay type 1, and which could be used in the definition of
arithmetically Gorenstein subschemes of $\proj{n}$.

\begin{proposition} \label{equiv to gor}
Let $X \subset \proj{n}$ be \acm.  The following are equivalent:
\begin{itemize}
\item[(i)] $X$ has Cohen-Macaulay type 1 (i.e.\ is arithmetically Gorenstein);
\item[(ii)] $R/I_X \cong K_X (\ell)$ for some $\ell \in {\mathbb Z}$, where
$K_X$ is the canonical module of $X$ (cf.\ Definition \ref{def-can-module});
\item[(iii)] The minimal free resolution of $R/I_X$ is self-dual up to twisting
by $n+1$.
\end{itemize}
\end{proposition}

\begin{proof}
Note that $\ell$ is whatever twist moves the module so that it starts in degree
0.  The main facts used in the proof are that
\[
\begin{array}{rcl}
K_X & = & \Ext_R^c (R/I_X,R)(-n-1) \\
\hbox{and} \ \ \ I_X & = & \Ann_R(K_X)
\end{array}
\]
Details of the proof can be found in \cite{migbook}.
\end{proof}

\begin{corollary}
Let $X$ be arithmetically  Gorenstein.  Then ${\mathcal O}_X \cong \omega_X
(\ell)$ for some $\ell \in {\mathcal Z}$.
\end{corollary}

\begin{corollary}
Let $X$ be arithmetically Gorenstein.  Then the $h$-vector of $X$ is symmetric.
\end{corollary}

\begin{proof}
This follows from the fact that the Gorenstein property is preserved in passing
to the Artinian reduction, and the Hilbert function of the canonical module of
the Artinian reduction is given by reading the $h$-vector backwards (cf.
\cite{migbook}).
\end{proof}

The integer $\ell$ in Proposition \ref{equiv to gor} is related to the integers
in the equation (\ref{relation betw a, c and t}).  In fact, we have

\begin{corollary} \label{cor-shift}
Let $X$ be arithmetically Gorenstein with minimal free resolution
\[
0 \rightarrow R(-a) \rightarrow F_{c-1} \rightarrow \dots \rightarrow F_1
\rightarrow R \rightarrow R/I_X \rightarrow 0
\]
and assume that ${\mathcal O}_X \cong \omega_X (\ell)$.
Then $\ell = n+1-a$.
\end{corollary}

If $A$ is Gorenstein then the integer $s$, the last degree in which the
$h$-vector is non-zero, is called the {\em socle degree} of $\bar A$, the
Artinian reduction of $A = R/I_X$.

There is a very useful criterion for zeroschemes to be arithmetically
Gorenstein.  To explain it, we will need a new notion.  For now we will assume
that our zeroschemes are reduced, although the necessity for this was
removed by
Kreuzer \cite{kreuzer}.

\begin{definition} \label{def of CB and UPP}
Let $Z \subset \proj{n}$ be a finite reduced set of points.  Assume that $s+1 =
\reg ({\mathcal I}_X)$, i.e.\ $s$ is the last degree in which the $h$-vector of
$Z$ is non-zero.  Then $Z$ has the {\em Cayley-Bacharach property} (CB) if, for
every subset $Y \subset Z$ consisting of $|Z| -1$ points, we have $h_{R/I_Y}
(s-1) = h_{R/I_Z}(s-1).$  $Z$ has the {\em Uniform Position property} (UPP) if
any two subsets $Y, Y'$ of (the same) arbitrary cardinality have the same
Hilbert
function, which necessarily is
\[
h_{R/I_Y} (t) = \min \{ h_{R/I_Z} (t), |Y| \} \ \ \ \hbox{ for all } t.
\]
\end{definition}

\begin{example}
The Cayley-Bacharach property is a weaker version of the Uniform Position
Property.  For example, in $\proj{2}$ consider the following examples.

\bigskip

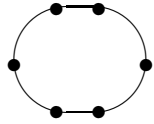
\newsavebox{\UPP}
\savebox{\UPP}(400,65)[tl]
{
\begin{picture}(400,125)
\put(140,40){\oval(50,40)}

\put (112,35){$\bullet$}
\put (162,35){$\bullet$}
\put (128,56){$\bullet$}
\put (128,17){$\bullet$}
\put (144,56){$\bullet$}
\put (144,17){$\bullet$}

\put (190,47){$h$-vector $1 \ 2 \ 2 \ 1$ (complete intersection on a conic)}
\put (190,27){This has UPP.}
\end{picture}
}

\newsavebox{\sUPP}
\savebox{\sUPP}(400,65)[tl]
{
\begin{picture}(400,125)
\put (110,30){\line (1,0){62}}
\put (110,50){\line (1,0){62}}
\put (120,20){\line (0,1){40}}
\put (140,20){\line (0,1){40}}
\put (160,20){\line (0,1){40}}

\put (117,27){$\bullet$}
\put (117,47){$\bullet$}
\put (137,27){$\bullet$}
\put (137,47){$\bullet$}
\put (157,27){$\bullet$}
\put (157,47){$\bullet$}

\put (190,47){$h$-vector $1 \ 2 \ 2 \ 1$ (complete intersection)}
\put (190,27){This has CB but not UPP.}
\end{picture}
}

\newsavebox{\CB}
\savebox{\CB}(400,65)[tl]
{
\begin{picture}(400,125)
\put (110,30){\line (1,0){62}}
\put (110,50){\line (1,0){62}}

\put (117,27){$\bullet$}
\put (117,47){$\bullet$}
\put (137,27){$\bullet$}
\put (137,47){$\bullet$}
\put (157,27){$\bullet$}
\put (157,47){$\bullet$}

\put (190,47){$h$-vector $1 \ 2 \ 2 \ 1$ (complete intersection)}
\put (190,27){This has CB but not UPP.}
\end{picture}
}

\newsavebox{\notCB}
\savebox{\notCB}(400,65)[tl]
{
\begin{picture}(400,125)
\put (90,30){\line (1,0){82}}
\put (90,50){\line (1,0){82}}

\put (97,27){$\bullet$}
\put (117,27){$\bullet$}
\put (117,47){$\bullet$}
\put (137,27){$\bullet$}
\put (137,47){$\bullet$}
\put (157,27){$\bullet$}

\put (190,47){$h$-vector $1 \ 2 \ 2 \ 1$}
\put (190,27){This has neither CB nor UPP.}
\end{picture}
}

\begin{picture}(400,65)
\put(-50,50){\usebox{\UPP}}
\end{picture}

\begin{picture}(400,65)
\put(-50,50){\usebox{\CB}}
\end{picture}

\begin{picture}(400,65)
\put(-50,50){\usebox{\notCB}}
\end{picture}

\end{example}

\begin{theorem}[\cite{DGO}] \label{DGO theorem}
A reduced set of points $Z$ is arithmetically Gorenstein if and only if its
$h$-vector is symmetric and it has the Cayley-Bacharach property.
\end{theorem}

\begin{example}
A set of $n+2$ points in $\proj{n}$ in linear general position is
arithmetically
Gorenstein.  In particular, a set of 5 points in $\proj{3}$ is arithmetically
Gorenstein, so we see that 4 points in linear general position are G-linked to
one point.  This was the first illustration of the fact that G-liaison behaves
quite differently from CI-liaison, since it follows from work of Ulrich and
others that 4 points in linear general position are not CI-linked to a single
point in any number of steps.
\end{example}

\begin{remark}
Theorem \ref{DGO theorem} was used by Bocci and Dalzotto \cite{BD} to
produce (and
verify) nice concrete examples of arithmetically Gorenstein sets of points in
$\proj{3}$, and this work is described in this
volume.  Generalizations of this construction have been given by Bocci,
Dalzotto, Notari and Spreafico \cite{BDNS}.
\end{remark}

A very useful construction of arithmetically Gorenstein schemes is the
following.

\begin{theorem}\label{sum of linked ideals thm}
{\bf (Sums of geometrically linked ideals)}
Let $V_1, V_2 \subset \proj{n}$ be \acm\ subschemes of codimension $c$ with no
common component.  Assume that $V_1 \cup V_2 = X$ is arithmetically Gorenstein,
i.e.\ $I_{V_1} \cap I_{V_2} = I_X$ with $R/I_X$ Gorenstein.  Then $I_{V_1} +
I_{V_2}$ is Gorenstein of codimension $c+1$ (i.e.\ $V_1 \cap V_2$ is
arithmetically Gorenstein).
\end{theorem}

\begin{proof}
  From the exact sequence (\ref{int dirsum sum}) we can build up the diagram
\[
\begin{array}{cccccccccccccccccccc}
&&      0     &&& 0 \\
&& \downarrow &&& \downarrow \\
&& R(-a) && A_c & \oplus & B_c \\
&& \downarrow &&& \downarrow \\
&& F_{c-1} && A_{c-1} & \oplus & B_{c-1} \\
&& \downarrow &&& \downarrow \\
&& \vdots &&& \vdots \\
&& \downarrow &&& \downarrow \\
&& F_1 && A_1 & \oplus & B_1 \\
&& \downarrow &&& \downarrow \\
0 & \rightarrow & I_X & \longrightarrow & I_{V_1} & \oplus & I_{V_2} &
\longrightarrow & I_{V_1} + I_{V_2} & \rightarrow & 0 \\
&& \downarrow &&& \downarrow \\
&&0 &&& 0
\end{array}
\]
The mapping cone then gives the long exact sequence
\[
\begin{array}{c}
\begin{array}{rclccccccccccccccccccccccc}
0 \rightarrow R(-a) \rightarrow F_{c-1} \oplus A_c \oplus B_c \rightarrow
F_{c-2} \oplus A_{c-1} \oplus B_{c-1} \rightarrow \dots \hbox{\hskip 2in}
\end{array}
\\ \\
\begin{array}{rclccccccccccccccccccccccc}
\hbox{\hskip 1in} \dots \rightarrow F_1 \oplus A_2 \oplus B_2 \rightarrow A_1
\oplus B_1 & \hskip -.4cm
\longrightarrow & \hskip -.4cm R  \rightarrow R/(I_{V_1} + I_{V_2})
\rightarrow 0
\\
& \hskip -1.3cm \searrow & \hskip -.9cm \nearrow \\
& \hskip 0cm I_{V_1} + I_{V_2} \\
& \hskip -1.3cm \nearrow & \hskip -.9cm \searrow \\
& \hskip -.4cm  0 \hfill & \hskip -.3cm 0
\end{array}
\end{array}
\]
Of course there may be some splitting.  However, $V_1 \cap V_2$ has codimension
$\geq c+1$ since $V_1$ and $V_2$ have no common component.  This resolution has
homological dimension at most $c+1$.  Therefore it has homological dimension
exactly $c+1$ and $V_1 \cap V_2$ is \acm\ of codimension $c+1$ with
Cohen-Macaulay type 1, i.e.\ is arithmetically Gorenstein.
\end{proof}

This construction has been used to good effect in constructing arithmetically
Gorenstein schemes with nice properties.  To illustrate, let us consider some
natural questions.

\begin{question} \label{possible artin gor hf}
What are the possible Hilbert functions (resp.\ minimal free resolutions) of
Artinian Gorenstein ideals?
\end{question}

\begin{question} \label{possible reduced gor hf}
What are the possible Hilbert functions (resp.\ minimal free resolutions)
of the
ideals of {\em reduced} arithmetically Gorenstein subschemes of $\proj{n}$?
\end{question}

The general question of which Artinian ideals, or which properties of Artinian
ideals, can be lifted to reduced sets of points is a very interesting one.  We
will discuss some of the known answers to Questions \ref{possible artin gor hf}
and \ref{possible reduced gor hf} according to the codimension.

\bigskip

\noindent \underline{Case I}: Codimension Two

What are the possible arithmetically Gorenstein subschemes $X$?  We know the
beginning and the end of the resolution:
\[
0 \rightarrow R(-a) \rightarrow ( \ ?? \ ) \rightarrow R \rightarrow R/I_X
\rightarrow 0.
\]
By considering the rank, the middle term in this resolution has to have
rank 2.
Therefore, we have established the well known fact (mentioned before) that {\em
every arithmetically Gorenstein subscheme of $\proj{n}$ of codimension two is a
complete intersection.}  This answers the question about the minimal free
resolution, so the Hilbert functions are known as well.  In fact, the
$h$-vectors must be symmetric of the form
\[
1 \ \ 2 \ \ 3 \ \ \dots s-1 \ \ s \ \ s \ \ \dots \ \ s \ \ s-1 \ \ \dots \ \ 3
\ \ 2 \ \ 1.
\]

\bigskip

\noindent \underline{Case II}: Codimension Three

Everything that is known in this case follows from the famous structure theorem
of Buchsbaum and Eisenbud \cite{BuEi}.  For a Gorenstein ideal $I$ we have a
minimal free resolution
\[
0 \rightarrow R(-a) \rightarrow F_2 \stackrel{A}{\longrightarrow} F_1
\rightarrow R \rightarrow R/I \rightarrow 0.
\]
One can choose bases so that $A$ is skew-symmetric.  In particular, the number
of generators must be odd!  Diesel used this result to completely describe the
possible graded Betti numbers for Artinian Gorenstein ideals.  De Negri and
Valla (and others) described the possible Hilbert functions.  In
particular, not
only must it be symmetric, but the ``first half'' must be a so-called {\em
differentiable O-sequence}.  This means that the first difference of the
``first
half'' of the Hilbert function must grow in a way that is permissible for
standard $K$-algebras.  For example, the sequence
\[
1 \ \ 3 \ \ 6 \ \ 7 \ \ 9 \ \ 7 \ \ 6 \ \ 3 \ \ 1
\]
is not a possible Hilbert function for an Artinian Gorenstein algebra (even
though it itself satisfies Macaulay's growth condition) since the first
difference of the ``first half'' is
\[
1 \ \ 2 \ \ 3 \ \ 1 \ \ 2
\]
and the growth from degree 3 to degree 4 in the first difference exceeds
Macaulay's growth condition (cf.\ \cite{mr.macaulay}).  This describes the
answers to Question \ref{possible artin gor hf}.

For Question \ref{possible reduced gor hf}, Geramita and Migliore \cite{GM5}
showed that any set of graded Betti numbers which occurs at the Artinian level
in fact occurs for a reduced set of points (or for a stick figure curve, or
more
generally a ``generalized stick figure'' configuration of linear varieties).
The idea was to use Theorem \ref{sum of linked ideals thm} and add the
ideals of
geometrically linked stick figure curves in $\proj{3}$ (or suitable surfaces in
$\proj{4}$, etc.) in suitable constructed complete intersections.  Ragusa and
Zappal\`a \cite{RZ1}, \cite{RZ2} have used the ``sum of geometrically linked
ideals'' construction to obtain other nice results on the Hilbert functions and
resolutions of codimension three Gorenstein ideals.

\bigskip

\noindent \underline{Case III}: Codimension $\geq 4$

To date no one has determined what Hilbert functions can occur, so certainly we
do not know what minimal free resolutions can occur.  In codimension $\geq
5$ it
is known that the Hilbert function of an Artinian Gorenstein algebra does not
even have to be unimodal \cite{ber-iar}, \cite{boij}, \cite{BL}.  This is
open in
codimension 4.  However, the situation that one would expect to be the
``general'' one is better understood:

\begin{definition} \label{wlp definition}
An Artinian algebra $R/I$ has the {\em Weak Lefschetz property} if, for a
general
linear form $L$, the multiplication map
\[
\times L : (R/I)_i \rightarrow (R/I)_{i+1}
\]
has maximal rank, for all $i$.
\end{definition}

When  the socle degree is fixed,  a result of Watanabe
\cite{watanabe} says that the ``general'' Artinian Gorenstein algebra has
the Weak Lefschetz property.

When the whole Hilbert function is fixed, a similar result is not possible
in general because the parameter space for the corresponding Gorenstein
algebras can have several components if the codimension is at least
four. However, since having the Weak Lefschetz property is an open
condition by semicontinuity, the general Artinian Gorenstein algebra of a
component has the Weak Lefschetz property if and only if the component
contains one algebra with this property.

In any case, Harima \cite{harima} classified the possible
Hilbert functions for Artinian Gorenstein algebras with the Weak Lefschetz
property, in any codimension.  In particular, he showed that these Hilbert
functions are precisely the Stanley-Iarrobino (SI) sequences, \label{def of
SI-sequence} namely they are symmetric, unimodal and the ``first half'' is a
differentiable O-sequence.

\label{MN3 ref} For Question \ref{possible reduced gor hf}, Migliore and Nagel
\cite{MN3} have shown that any SI-sequence is the $h$-vector of some
arithmetically Gorenstein reduced set of points, or more generally a reduced
union of linear varieties. The method of proof again used sums of geometrically
linked ideals, but the new twist here was that the ideals were G-linked and not
CI-linked.  Furthermore, they gave sharp bounds on the graded Betti numbers of
Gorenstein ideals of any codimension, {\em among ideals with the Weak Lefschetz
property}.  Partial results along these lines had been obtained by Geramita,
Harima and Shin
\cite{GHS}.  In codimension 4, Iarrobino and Srinivasan (in progress) have some
results on the possible resolutions.  There remains a great deal to do in this
area.

\begin{remark}
Theorem \ref{sum of linked ideals thm} shows how to use geometrically linked,
codimension $c$, \acm\ subschemes of $\proj{n}$ to construct a codimension
$c+1$ \aG\
subscheme.  Later, in Corollary \ref{cor-aci-gor}, we will see how to use very
special linked \acm\ codimension $c$ subschemes (not necessarily geometrically
linked) to construct an \aG\ subscheme which is also of codimension $c$.
In fact,
every Gorenstein ideal arises in this way (Remark \ref{rem-constr-gor})!
\end{remark}

One problem with the construction of Theorem \ref{sum of linked ideals thm} is
that it is very desirable, from the point of view of liaison, to be able to
start with a scheme $V$ and find a ``good'' (which often means ``small'')
Gorenstein scheme $X$ containing it.  This is not so easy to do with sums of
geometrically linked ideals.   Another very useful construction for Gorenstein
ideals potentially will solve this problem (based on experimental
evidence).  To
describe it we will need a little preparation.

Consider a homogeneous map
\[
\bigoplus_{i=1}^{t+r} R(-a_i) \stackrel{\phi}{\longrightarrow}
\bigoplus_{j=1}^t
R(-b_j).
\]
The map $\phi$ is represented by a homogeneous $t \times (t+r)$ matrix.  We
assume furthermore that the ideal of maximal minors of $\phi$ defines a scheme
of the ``expected'' codimension $r+1$.  Let $B_\phi$ be the kernel of $\phi$.
Then $B_\phi$ is a {\em Buchsbaum-Rim module}.  Let $M_\phi$ be the cokernel of
$\phi$.  We have an exact sequence
\[
\begin{array}{ccccccccccccccccccc}
0 & \rightarrow & B_\phi & \rightarrow & \bigoplus_{i=1}^{t+r} R(-a_i)
& \stackrel{\phi}{\longrightarrow} & \bigoplus_{j=1}^t R(-b_j) & \rightarrow
& M_\phi & \rightarrow & 0. \\
&&&& || && || \\
&&&& F && G
\end{array}
\]
Sheafifying this gives
\[
0 \rightarrow \tilde B_\phi \rightarrow \bigoplus_{i=1}^{t+r} {\mathcal
O}_{\proj{n}} (-a_i)
\stackrel{\phi}{\longrightarrow} \bigoplus_{j=1}^t {\mathcal O}_{\proj{n}}
(-b_j)
\rightarrow \tilde M_\phi \rightarrow 0.
\]
If $r=n$ then $\tilde M = 0
$ and $\tilde B_\phi$ is locally free.  In any case,
$\tilde B_\phi$ is the {\em Buchsbaum-Rim sheaf associated to $\phi$}.

\begin{theorem}[\cite{MNP1}] \label{sections of BR sheaves}
{\bf (Sections of Buchsbaum-Rim sheaves)}
Assume that $r$ is odd.  Let $s$ be a {\em regular} section of $\tilde
B_\phi$.
Let $I$ be the ideal corresponding to the vanishing of~$s$.  Then the top
dimensional part of $I$ is arithmetically Gorenstein of codimension $r$.
Denoting by $J$ this top dimensional part, the minimal free resolution of $R/J$
can be written in terms of $F$ and $G$.
\end{theorem}

This can be used to find an arithmetically Gorenstein scheme (of the same
dimension) containing a given one by  means of the following corollary.

\begin{corollary}
If $\codim V = r$, a regular section of $H^0_* (\proj{n}, \tilde B_\phi \otimes
{\mathcal I}_V)$ has top dimensional part which is an arithmetically Gorenstein
scheme, $X$, containing $V$.
\end{corollary}

Theorem \ref{sections of BR sheaves} is just a small sample (the application
relevant to liaison) of the possible results on sections of Buchsbaum-Rim
sheaves, and we refer the interested reader to \cite{MNP1} for more general
results.

Our final construction requires a little preparation.

\begin{definition} \label{def of Gr}
A subscheme $S \subset \proj{n}$ satisfies condition $G_r$ if every
localization
of $R/I_S$ of dimension $\leq r$ is a Gorenstein ring.  $G_r$ is sometimes
referred to as ``Gorenstein in codimension $\leq r$, i.e.\ the ``bad
locus'' has
codimension $\geq r+1$.
\end{definition}

\begin{definition} Let $S \subset \proj{n}$ be an \acm\ subscheme and let
$F$ be
a homogeneous polynomial of degree $d$ not vanishing on any component of $S$
(i.e.\ $I_S : F = I_S$).  Then $H_F$ is the divisor on $S$ cut out by $F$.  We
call $H_F$ the {\em hypersurface section of $S$ cut out by $F$.}  As a
subscheme
of $\proj{n}$, $H_F$ is defined by the ideal $I_S + (F)$.  Note that this ideal
is saturated, since $S$ is \acm.  (The idea is the same as that in Lemma
\ref{hyperpl sect lemma} and Remark \ref{rem-hyp-sect}.)
\end{definition}

Hartshorne \cite{gen-div} has developed the theory of divisors, and in
particular linear equivalence, on schemes having at least $G_1$.  Using the
notion of linear equivalence, the following theorem gives our construction.  In
subsequent sections we will give some important applications for liaison.

\begin{theorem}[\cite{KMMNP}] \label{tw anticanon div}
{\bf (Twisted anticanonical divisors)} Let $S \subset \proj{n}$ be
an \acm\ subscheme satisfying $G_1$ and let $K$ be a canonical
divisor of $S$.   Then every effective divisor in the linear
system $|dH-K|$, viewed as a subscheme of $\proj{n}$, is
arithmetically Gorenstein.
\end{theorem}

\begin{proof} (Sketch)

Let $X \in |dH-K|$ be an effective divisor. Choose a sufficiently
large integer $\ell$ such that there is a regular section of
$\omega_S (\ell)$ defining a twisted canonical divisor $Y$.  Let
$F \in I_Y$ be a homogeneous polynomial of degree $d + \ell$ such
that $F$ does not vanish on any component of $S$ and let $H_F$ be
the corresponding hypersurface section.

Then $X$ is linearly equivalent to the effective divisor $H_F - Y$
and we have isomorphisms
\[
\left ( {\mathcal I}_X/{\mathcal I}_S \right ) (d) \cong {\mathcal
I}_{X|S} (d) \cong {\mathcal O}_S ((d+\ell) H-X) \cong {\mathcal
O}_S (Y) \cong \omega_S (\ell).
\]
Because $S$ is \acm, this gives
\[
0 \rightarrow I_S \rightarrow I_X \rightarrow H^0_* (\omega_S)(\ell -d)
\rightarrow 0.
\]
Then considering a minimal free resolution of $I_S$ and the corresponding one
for $K_S = H^0_*(\omega_S)$ (cf.\ Remark \ref{rem-can-mod}) we have a diagram
(ignoring twists)
\[
\begin{array}{ccccccccccccccccccccc}
&&&&&& 0 \\
&&&&&& \downarrow \\
&& 0 &&&& R \\
&& \downarrow &&&& \downarrow \\
&& F_c &&&& F_1^* \\
&& \downarrow &&&& \downarrow \\
&& \vdots &&&& \vdots \\
&& \downarrow &&&& \downarrow \\
&& F_1 &&&& F_c^* \\
&& \downarrow &&&& \downarrow \\
0 & \rightarrow & I_S & \rightarrow & I_X & \rightarrow & K_S & \rightarrow & 0
\\
&& \downarrow &&&& \downarrow \\
&& 0 &&&& 0
\end{array}
\]
Then the Horseshoe Lemma (\cite{weibel} 2.2.8, p.\ 37) shows that $I_X$ has a
free resolution in which the last free module has rank one.  Since $\codim X =
c$, this last free module cannot split off, so $X$ is arithmetically Gorenstein
as claimed.
\end{proof}

\begin{example}
Let $S$ be a twisted cubic curve in $\proj{3}$.  Then a canonical divisor
$K$ has degree $-2$, so the linear
system
$|-K+dH|$ (for $d \geq 0$) consists of all effective divisors of degree $\equiv
2$ (mod 3).  Any such scheme is arithmetically Gorenstein.
\end{example}


\section{First relations between linked schemes} \label{sec-relations}

In this section we begin to investigate the relations between
linked ideals. In particular, we will compare the Hilbert
functions of directly linked ideals and cover some of the results
announced in Section \ref{sec-overv}.

All the ideals will be homogeneous ideals of the polynomial ring
$R = K[x_0,\ldots,x_n]$. Analogously to linked schemes we define
linked ideals. This includes linkage of Artinian ideals
corresponding to empty schemes.

\begin{definition} \label{def-link-ideals}
(i) Two unmixed ideals $I, J \subset R$ of the same codimension
are said to be {\it geometrically CI-linked} (resp.\ {\it geometrically
G-linked}) by the ideal $\fc$ if $I$ and $J$ do not have
associated prime ideals in common and $\fc = I \cap J$ is a
complete intersection (resp.\ a Gorenstein ideal).

(ii) Two ideals $I, J \subset R$ are said to be {\it (directly) CI-linked}
(resp.\ {\it (directly) G-linked}) by the ideal $\fc$ if  $\fc$ is a
complete intersection (resp.\ a Gorenstein ideal) and
$$
\fc : I = J \quad \mbox{and} \quad \fc : J = I.
$$
In this case we write $I \stackrel{\fc}{\sim} J$.
\end{definition}

If a statement is true for CI-linked ideals and G-linked ideals we
will just speak of linked ideals.

\begin{remark}
(i) If we want to stress the difference between (i) and (ii) we
say in case (ii) that the ideals are algebraically linked.

(ii) If two ideals are geometrically linked then they are also
algebraically linked.

(iii) Since Gorenstein ideals of codimension two are complete
intersections CI-linkage is the same as G-linkage for ideals of
codimension two.
\end{remark}

If the subschemes $V$ and $W$ are geometrically linked by $X$ then
  $\deg X = \deg V + \deg W$. We will see that this equality is
  also true if $V$ and $W$ are only algebraically linked. For this
  discussion we will use the following.

\begin{notation}
$I, \fc \subset R$ denote homogeneous ideals where $\fc$ is a
Gorenstein ideal of codimension $c$. Excluding only trivial cases
we assume $c \geq 2$.
\end{notation}

\begin{lemma} \label{lem-colon-unm}
If $I \nsubseteqq \fc$ then $\fc : I$ is an unmixed ideal of
codimension $c$.
\end{lemma}

\begin{proof}
Let $\fc = \fq_1 \cap \ldots \cap \fq_s$ be a shortest primary
decomposition of $\fc$. Then the claim follows because
$$
\fc : I = \bigcap_{i=1}^s (\fq_i : I)
$$
and
$$
\fq_i : I = \left \{ \begin{array}{cl}
R & \mif I \subset \fq_i \\
\Rad (\fq_i)-\mbox{primary} & \mbox{otherwise}.
\end{array} \right.
$$
\end{proof}

The next observation deals with the difference between geometric
and algebraic linkage.

\begin{corollary} \label{cor-alg-vers-geo}
Suppose the ideals $I$ and $J$ are directly linked by $\fc$. Then we have:
\begin{itemize}
\item[(a)] $\Rad (I \cap J) = \Rad \fc$.
\item[(b)] $I$ and $J$ are unmixed of codimension $c$.
\item[(c)] If $I$ and $J$ do not have associated prime ideals in
common  then $I$ and $J$ are geometrically linked by $\fc$
\end{itemize}
\end{corollary}

\begin{proof}
(a) By definition we have
$$
I \cdot J \subset \fc \subset I \cap J.
$$
Since $\Rad (I \cdot J) = \Rad (I \cap J)$ the claim follows.

(b) is an immediate consequence of the preceding lemma.

(c) We have to show that $I \cap J = \fc$.

The inclusion $\fc \subset I \cap J$ is clear. For showing the other
inclusion assume on the contrary that there is a homogeneous polynomial $f
\neq 0$ in $(I \cap J) \setminus \fc$. Let $\fc = \fq_1 \cap \ldots \cap
\fq_s$ be a shortest primary  decomposition of $\fc$. We may assume that $f
\notin \fq_1$ and $\Rad \fq_1 \in \Ass_R (R/I)$. The assumption on the
associated prime ideals of $I$ and $J$ guarantees that there is a
homogeneous polynomial $g \in J$ such that $g \notin \fp$ for all $\fp \in
\Ass_R (R/I)$. Since $I = \fc : J$ we get $f g \in \fc \subset
\fq_1$. Thus, $g \notin \Rad \fq_1$ implies $f \in \fq_1$, a
contradiction.
\end{proof}

Claim (a) of the statement above says for schemes $V, W$ linked by
$X$ that we have (as sets) $V_{red} \cup W_{red} = X_{red}$.

In order to identify certain degree shifts we need the following number. It
is well-defined because the Hilbert function equals the Hilbert polynomial
in all sufficiently large degrees.

\begin{definition}  \label{def-ref-index}
The {\it regularity index} of a finitely generated graded $R$-module $M$ is
the number
$$
r(M) := \min \{ t \in \ZZ \s h_M(j) = p_M (j) \; \; \mbox{for all} \; j
\geq t \}.
$$
\end{definition}

\begin{example} \label{ex-reg-ind}
(i) $ r(K[x_0,\ldots,x_n]) = -n$.

(ii) If $A$ is an Artinian graded $K$-algebra with $s = \max \{ j \in \ZZ
\s [A]_j \neq 0 \}$ then $r(A) = s+1$.

(iii) Let $\fc \in R$ be a Gorenstein ideal with minimal free resolution
$$
0 \to R(-t) \to F_{c-1} \to \ldots \to F_1 \to \fc \to 0.
$$
Then it is not to difficult to see that $r(R/\fc) = t - n$ (cf.\ also
Corollary \ref{cor-shift}).
\end{example}

The index of regularity should not be confused with the Castelnuovo-Mumford
regularity. There is the following comparison result.

\begin{lemma} \label{lem-reg-v-index}
Let $M$ be a graded $R$-module. Then we have
$$
\reg M - \dim M + 1 \leq r(M) \leq \reg M - \depth M  + 1.
$$
In particular, $r(M) = reg M - \dim M +1$ if $M$ is \CM.
\end{lemma}

This lemma generalizes Example \ref{ex-reg-ind}(ii) and is a consequence of the
following version of the Riemann-Roch theorem \cite{serre}  which we will
use again soon.

\begin{lemma} \label{lem-pol-v-fkt}
Let $M$ be a graded $R$-module. Then we have for all $j \in \ZZ$
$$
h_M (j) - p_M (j) = \sum_{i \geq 0} (-1)^i \rankk [\HH^i(M)]_j.
$$
\end{lemma}

We are now ready for a crucial observation.

\begin{lemma}[Standard exact sequences] \label{lem-st-ex-seq}
Suppose that $\fc \subset I$ and both ideals have the same
codimension $c$. Put $J := \fc : I$. Then there are exact
sequences (of graded $R$-modules)
$$
0 \to \fc \to J \to K_{R/I} (1-r(R/\fc)) \to 0
$$
and
$$
0 \to K_{R/I} (1-r(R/\fc)) \to R/\fc \to R/J \to 0.
$$
\end{lemma}

\begin{proof}
We have to show that $J/\fc \cong K_{R/I} (1-r(R/\fc))$.

There are the following isomorphisms
$$
J/\fc \cong (\fc : I)/\fc \cong \Hom_R (R/I, R/\fc) \cong
\Hom_{R/\fc} (R/I, R/\fc)
$$
and
$$
K_{R/I} (1-r(R/\fc)) \cong \Ext^c_R (R/I, R)(- r(R/\fc) - n).
$$
Thus, our claim follows from the isomorphism
$$
\Hom_R (R/I, R/\fc) \cong \Ext^c_R (R/I, R)(- r(R/\fc) - n),
$$
which is easy to see if $\fc$ is a complete intersection. In the
general case it follows from a more abstract characterization of
the canonical module.
\end{proof}

Before drawing first consequences we recall that the Hilbert
polynomial of the graded module $M$ can be written in the form
$$
p_M (j) = h_0(M) \binom{j}{d-1} +  h_1(M) \binom{j}{d-2} + \ldots
+ h_{d-1}(M)
$$
where $d = \dim M$ and $h_0(M),\ldots,h_{d-1}(M)$ are integers.
Moreover, if $d > 0$ then $\deg M = h_0(M)$ is positive and called
the degree of $M$. However, if $M = R/I$ for an ideal $I$ then by
abuse of notation we define $\deg I := \deg R/I$. For a subscheme
$X \subset \PP^n$ we  have then $\deg X = \deg I_X$.

\begin{corollary} \label{cor-dir-linked-id}
Let $I$ be an ideal of codimension $c$ which contains the
Gorenstein ideal $\fc$. Put $J := \fc : I$ and $s := r(R/\fc) -1$.
Then we have
\begin{itemize}
\item[(a)] $$
\deg J = \deg \fc - \deg I,
$$
and if $c < n$ and $I$ is unmixed then
$$
h_1(R/J) = \frac{1}{2} (s- n+c+1) [\deg I - \deg J] + h_1(R/I).
$$
\item[(b)] If $I$ is unmixed and $R/I$ is \lCM\ then also $R/J$ is
\lCM\ and
$$
\HH^i(R/J) \cong \HH^{n+1-c-i}(R/I)^{\vee} (-s) \quad (i =1,\ldots,n-c)
$$
and
$$
p_{R/J} (j) = p_{R/\fc}(j) + (-1)^{n-c} p_{R/I}(s-j).
$$
\item[(c)] If $R/I$ is \CM\ then also $R/J$ has this property and
$$
h_{R/J} (j) = h_{R/\fc}(j) + (-1)^{n-c} h_{R/I}(s-j).
$$
\end{itemize}
\end{corollary}

\begin{proof}
  Consider the version of Riemann-Roch (Lemma \ref{lem-pol-v-fkt})
$$
h_{R/I} (j) - p_{R/I} (j) = \sum_{i \geq 0} (-1)^i \rankk [\HH^i({R/I})]_j.
$$
Since the degree of the Hilbert polynomial
of $\HH^i(R/I)$ is at most $i-1$, we obtain for all  $j \ll 0$
$$
\begin{array}{rcl}
-p_{R/I} (j) & = & (-1)^{n+1-c} \rankk [\HH^{n+1-c}({R/I})]_j +
O(j^{n-1-c}) \\[1ex]
& = & (-1)^{n+1-c} \rankk [K_{R/I}]_{-j} + O(j^{n-1-c}).
\end{array}
$$
Combined with the standard exact sequence this provides
$$
p_{R/J} (j) = p_{R/\fc} (j) + (-1)^{n-c} p_{R/I} (s-j) + O(j^{n-1-c}).
$$
Comparing coefficients we get if $c \leq n$
$$
\deg J = \deg \fc - \deg I.
$$
Assume now that the ideal $I$ is unmixed. Lemma \ref{lem-unmixed-crit} implies
that then the degree of the Hilbert polynomial of $\HH^i(R/I)$ is at most
$i-2$.  Thus,
we obtain as above
$$
p_{R/J}S (j) = p_{R/\fc} (j) + (-1)^{n-c} p_{R/I} (s-j) + O(j^{n-2-c}).
$$
Hence, we get if $c < n$
$$
h_1(R/J) = (s-n+1+c) \deg I + h_1(R/I) + h_1(R/\fc).
$$
But by duality we have
$$
h_1(R/\fc) = (s-n+1+c) \deg \fc.
$$
Combining the last two equalities proves the second statement in (a).

The isomorphisms
$$
\HH^i(R/J) \cong \HH^{n+1-i}(R/I)^{\vee} (-s) \quad (i =1.\ldots,n-c)
$$
follow essentially from the long exact cohomology sequence induced by the
standard
sequence taking into
account Proposition \ref{prop-duality}.

The remaining claims in (a) - (c) are proved similarly as above. For
details
we refer to \cite{N-gorliaison}.
\end{proof}

\begin{remark} \label{rem-coho-liai}
If $I$ is unmixed but not \lCM\ then the formula in claim (b) relating the
local cohomologies of $I$ and $J$ is not true in general. For example, it
is never true if $I$ defines a non-\lCM\ surface in $\PP^4$.
\end{remark}

The last statement applies in particular to directly linked ideals. Thus,
we obtain.

\begin{corollary} \label{cor-CM-under-li}
Let $V$ and $W$ be directly linked. Then we have:
\begin{itemize}
\item[(a)] $V$ is \aCM\ if and only if $W$ is \aCM.
\item[(b)] $V$ is \lCM\ if and only if $W$ has this property.
\end{itemize}
\end{corollary}

Before turning to examples we want to rewrite Claim (a) in Corollary
\ref{cor-dir-linked-id} for curves in a more familiar form.

\begin{remark} \label{rem-genus-u-lia}
Let $C_1, C_2 \subset \PP^n$ be curves linked by an \aG\ subscheme $X$. Let
$g_1$ and $g_2$ denote the arithmetic genus of $C_1$ and $C_2$,
respectively.  Since by definition $g_i = 1 - h_1(R/I_{X_i})$,
Corollary \ref{cor-dir-linked-id} provides the formula
$$
g_1 - g_2 = \frac{1}{2} (r(R/I_X) - 1) \cdot [\deg C_1 - \deg
C_2].
$$
In particular, if $X$ is a  complete intersection
cut out by hypersurfaces of degree $d_1,\ldots,d_{n-1}$  we obtain
$$
g_1 - g_2 = \frac{1}{2} (d_1 + \ldots + d_{n-1} - n - 1) \cdot [\deg C_1 - \deg
C_2]
$$
because the index of regularity of $X$ is
$$
r(R/I_X) = d_1 + \ldots + d_{n-1} - n.
$$
\end{remark}

\begin{example} \label{ex-linkage}
(i) Let $C \subset \PP^3$ be the twisted cubic parameterized by $(s^3, s^2t,
st^2, t^3)$. It is easy to see that $C$ is contained in the complete
intersection $X$ defined by
$$
\fc := (x_0 x_3 - x_1 x_2,  x_0 x_2 - x_1^2) \subset I_C.
$$
Then Corollary \ref{cor-dir-linked-id} shows that $\fc : I_C$ has degree
$1$. In fact, we  easily get $\fc : I_C = (x_0, x_1)$ defining the line
$L$. Thus, $C$ and $L$ are geometrically linked by $X$ and $C$ is \aCM.

(ii)  Let $C \subset \PP^3$ be the rational quartic parameterized by
$(s^4, s^3t, st^3, t^4)$. $C$ is contained in the complete
intersection $X$ defined by
$$
\fc := (x_0 x_3 - x_1 x_2,  x_0 x_2^2 - x_1^2 x_3) \subset I_C.
$$
Hence  $\fc : I_C$ has degree $2$. Indeed, it is easy to see that
$$
\fc : I_C = (x_0, x_1) \cap (x_2, x_3).
$$
This implies that $C$ is geometrically linked to a pair of skew
lines. Therefore $C$ is not \aCM\ (thanks to Example \ref{ex-min-shift}).

(iii) We want to illustrate  Corollary \ref{cor-dir-linked-id} (c).
Let $I := (x^2, xy, y^4) \subset R := K[x, y]$ and let $\fc := (x^3, y^4)
\subset I$. We want to compute the Hilbert function of $R/J$. Consider the
following table

\medskip

\begin{center}
\begin{tabular}{c|cccccccc}
$j$ & 0 & 1 & 2 & 3 & 4 & 5 & 6  \\
\cline{1-8}
$h_{R/I} (j)$ &  1 & 2 & 1 & 1 & 0 & 0 & 0 \\[2pt]
$h_{R/\fc} (j)$ &  1 & 2 & 3 & 3 & 2 & 1 & 0 \\ [2pt]
$h_{R/J} (5-j)$ &  0 & 0 & 2 & 2 & 2 &  1& 0
\end{tabular}
\end{center}

\medskip
The second row shows that $r(R/\fc) = 6$. Thus by Corollary
\ref{cor-dir-linked-id} (c),  the last row is the second row minus the first
row. We get that $R/J$ has the Hilbert function $1, 2, 2, 2, 0, \ldots$.
\end{example}

We justify now Remark \ref{rem-alg-link}.

\begin{corollary} \label{cor-alg-link}
Suppose $I$ is an unmixed ideal of codimension $c$.  If the Gorenstein
ideal $\fc$ is properly contained in $I$ then the ideals $\fc : I$ and $I$
are directly linked by $\fc$.
\end{corollary}

\begin{proof}
We have to show the equality
$$
\fc : (\fc : I) = I.
$$
It is clear that $I \subset \fc : (\fc : I)$. Corollary
\ref{cor-dir-linked-id} (a) provides
$$
\deg [\fc : (\fc : I)] = \deg \fc - \deg(\fc : I) = \deg I.
$$
Since $\fc : (\fc : I)$ and $I$ are unmixed ideals of the same codimension
they must be equal.
\end{proof}

Finally, we want to show how CI-linkage can be used to produce Gorenstein
ideals. To this end we introduce.

\begin{definition} \label{def-aci}
An ideal $I \subset R$ is called an {\it \aci} if $R/I$ is \CM\ and $I$ can
be generated by $\codim I + 1$ elements.
\end{definition}

\begin{example}
The ideal $(x_0, x_1)^2$ is an \aci.

The twisted cubic $C \subset \PP^3$ is also an \aci.
\end{example}

\begin{corollary} \label{cor-aci-gor}
Let $I \subset R$ be an \aci\ and let $\fc \subsetneq I$ be a complete
intersection such that $\codim I = \codim \fc$ and $I = \fc + f R$ for some
$f \in R$. Then $J := \fc : I$ is a Gorenstein ideal.
\end{corollary}

\begin{proof}
Consider the standard exact sequence
$$
0 \to \fc \to  \fc + f R \to K_{R/J}(1 - r(R/\fc)) \to 0.
$$
It shows that $K_{R/J}$ has just one minimal generator (as $R$-module).

Let
$$
0 \to F_c \stackrel{\ffi_c}{\longrightarrow} \ldots \to F_1 \to J \to 0
$$
be a minimal free resolution. Then the beginning of the minimal free resolution
of $K_{R/J}$ has the form
$$
\ldots \to F_{c-1}^* \stackrel{\ffi_c^*}{\longrightarrow} F_c^* \to
K_{R/J}(n+1) \to 0.
$$
It follows that $F_c$ must have rank $1$, i.e.\ $J$ is a Gorenstein ideal.
\end{proof}

\begin{remark} \label{rem-constr-gor}
Every Gorenstein ideal arises as in the previous corollary.

We only sketch the argument. Given a Gorenstein $J$ of codimension $c$ we
choose a complete intersection $\fc$ of codimension $c$ which is properly
contained in $J$. Then $I := \fc : J$ is an \aci\ and $J = \fc : I$.
\end{remark}


\section{Some Basic Results and Constructions}

We begin this section by proving one of the results mentioned in Remark
\ref{some nice results}.

\begin{theorem} \label{linear equiv theorem}
Let $S \subset \proj{n}$ be \ACM\ satisfying property $G_1$ (so that
linear equivalence is well-defined; see Definition \ref{def of Gr}  and
the discussion preceding Theorem \ref{tw anticanon div}).  Let $C_1 ,C_2
\subset S$ be divisors such that $C_2 \in | C_1 + tH|$, where $H$ is the
class of a hyperplane section and $t \in {\mathbb Z}$.  Then $C_1$ and $C_2$
are G-linked in two steps.
\end{theorem}

\begin{proof}
Let $Y$ be an effective twisted canonical divisor.  Choose an integer $a
\in {\mathbb Z}$ such that $[I_Y]_a$ contains a form $A$ not vanishing on
any component of $S$.  Hence $H_A - Y$ is effective on $S$.

Now, recall that
\[
\begin{array}{rcl}
C_2 \in |C_1| & \Leftrightarrow & C_2 - C_1 \hbox{ is the divisor of a
rational function on $S$} \\
& \Leftrightarrow & \hbox{there exist $F,G$ of the same degree such that
} \left ( \frac{F}{G} \right ) = C_2 - C_1 \\
&& \hbox{(where $\left ( \frac{F}{G} \right )$ is the divisor of the
rational function $ \frac{F}{G}$)} \\
& \Leftrightarrow & \hbox{there exists a divisor $D$ such that } \\
&&
\hbox{\hskip .35in}
\begin{array}{rcl}
H_F & = & C_2 + D \\
H_G & = & C_1 + D
\end{array}
\\
&& \hbox{(in particular, $F \in I_{C_2}$ and $G \in I_{C_1}$).}
\end{array}
\]
Similarly,
\[
\begin{array}{rcl}
C_2 \in |C_1 + tH| & \Leftrightarrow & \hbox{there exist $F,G$ with $\deg
F = \deg G + t$ and a divisor $D$ such that}\\
&&
\hbox{\hskip .35in}
\begin{array}{rcl}
H_F & = & C_2 + D \\
H_G & = & C_1 + D
\end{array} \\
&& \hbox{(in particular, $F \in I_{C_2}$ and $G \in I_{C_1}$).}
\end{array}
\]
Note that the effective divisor $H_{AF} - Y$ is arithmetically Gorenstein,
by Theorem \ref{tw anticanon div}.  Then one checks if $S$ is smooth that
\[
(H_{AF} - Y) - C_2 = (H_A - Y) + (H_F - C_2) = (H_A - Y) + D
\]
and
\[
(H_{AG} - Y) - C_1 = (H_A - Y) + (H_G - C_1) = (H_A - Y) + D
\]
Therefore $C_2$ is directly linked to $(H_A - Y)+D$ by the Gorenstein
ideal $H_{AF} - Y$ and $C_1$ is directly linked to the same $(H_A - Y)+D$
by the Gorenstein ideal $H_{AG} - Y$.  This concludes the proof in the
special case. For the general case we refer to \cite{KMMNP}, Proposition 5.12.
\end{proof}

Theorem \ref{linear equiv theorem} was the first result that really
showed that Gorenstein liaison is a theory about divisors on \acm\
schemes, just as Hartshorne \cite{gen-div} had shown that CI-liaison is a
theory about divisors on complete intersections.   It is fair to say that
most of the results about Gorenstein liaison discovered in
the last few years use this result either directly or at least indirectly.

\begin{remark} \label{rem-Hartsh} 
As pointed out to us by R.\ Hartshorne, there is an interesting point
lurking in
the background here.  Following \cite{MP} and \cite{gen-div}, we say that a
subscheme $V_2 \subset \proj{n}$ is obtained from a subscheme $V_1 \subset
\proj{n}$ by an {\em elementary CI-biliaison} if there is a complete
intersection $S$ in $\proj{n}$ such that $V_2 \sim V_1+hH$ on $S$ for some
integer $h \geq 0$, where $\sim$ denotes linear equivalence.  It is not hard
to show, and has long been known, that $V_1$ and $V_2$ are CI-linked in two
steps.  It is a theorem (\cite{gen-div}) that the equivalence relation
generated by elementary CI-biliaisons is the same as the equivalence relation
of  even CI-liaison (see Definition \ref{def-even-li-cl}).

Now, Theorem \ref{linear equiv theorem} naturally suggests the idea of saying
that $V_2$ is obtained from  $V_1$ by an {\em elementary G-biliaison} if there
is an \acm\ scheme $S$ with property $G_1$ such that $V_2 \sim V_1+hH$ on $S$
for some integer $h \geq 0$.  It is an open problem to determine if the
equivalence relation generated by elementary G-biliaisons is the same as the
equivalence relation of even G-liaison.  It is conceivable that schemes $V_1$
and $V_2$ could be evenly G-linked, but no sequence of elementary G-biliaisons
beginning with $V_1$ can arrive at $V_2$.

As mentioned above, many of the results about Gorenstein liaison in fact use
elementary G-biliaisons, so the results are actually slightly stronger in this
sense.
\end{remark}

Theorem \ref{linear equiv theorem} clearly needs the $G_1$ assumption
since linear equivalence was used.  In general, without $G_1$, it is not
always possible to talk about divisors of the form $H_A - Y$.  However,
we will see now that there is a notion of ``adding'' a hypersurface
section even if the $G_1$ assumption is relaxed.  This construction was
given in Lemma \ref{lem-basic-d-link} and was called ``Basic Double
Link" there.  Now we will see why this name was chosen.  In Lemma
\ref{lem-basic-d-link} almost no assumption was made on the ideal $J$.
Here we present it in more geometric language, and we have to assume at
least $G_0$ in order to get a liaison result.

\begin{proposition}[Basic Double G-Linkage] \label{bdl is linked}
Let $S \subset \proj{n}$ be an \acm\ subscheme satisfying $G_0$.  Let 
$C \subset S$
be an equidimensional subscheme of codimension 1 and let $A \in R$ be
homogeneous with $I_S : A = I_S$ (i.e.\ $A$ does not vanish on any
component of $S$).  Then $I_C$ and $I_S + A \cdot I_C$ are G-linked in
two steps.
\end{proposition}

\begin{remark}
As we will see, the ideal $I_S + A \cdot I_C$ represents the divisor $C +
H_A$ on $S$.  If $\deg A = d$ and $S$ satisfies $G_1$ then the scheme $Y$
defined by the ideal $I_Y + I_S + A \cdot I_C$ is in the linear system
$|C + dH|$.  But in our level of generality, linear systems may not make
sense.  See \cite{kyoto} for a more detailed discussion of these divisors.
\end{remark}

\noindent {\em Proof of  Proposition \ref{bdl is linked}} (sketch)

The unmixedness statement in Lemma \ref{lem-basic-d-link} shows that in
particular, $I_Y = I_S + A \cdot I_C$ is saturated.  Furthermore, from
the exact sequence (\ref{bdl exact sequence}) we can make a Hilbert
function calculation:
\[
h_{R/I_Y} (t) = h_{R/I_S}(t) - h_{R/I_S}(t-d) + h_{R/I_C}(t-d).
\]
It follows that
\[
\begin{array}{rcl}
\deg Y & = & \deg C + d \cdot \deg S \\
& = & \deg C + \deg H_A.
\end{array}
\]
The idea of the proof is to mimic the proof of Theorem \ref{linear equiv
theorem} in an algebraic way.  We proceed in four steps:

\bigskip

\noindent \underline{\em Step I}: Let $c = \codim S$.  The $G_0$
hypothesis is enough to guarantee that there exists a {\em Gorenstein}
ideal $J \subset R$ with $I_S \subset J$, $\codim J = c+1$ and $J/I_S$ is
Cohen-Macaulay of Cohen-Macaulay type 1 (cf. \cite{KMMNP}).  Since
$\codim C > \codim S$, there exists $B \in I_C$ of some degree, $e$, such
that $I_S : B = I_S$ (i.e.\ $B$ does not vanish on any component of $S$).

\bigskip

\noindent \underline{\em Step II}: One checks that $I_S + B \cdot J$ is
Gorenstein and $I_S + B \cdot J \subset I_C$.  Hence $I_S + B \cdot J$
links $I_C$ to some ideal $\fa$ which is unmixed.

\bigskip

\noindent \underline{\em Step III}: $I_S + AB \cdot J$ is Gorenstein and
is contained in $I_Y = I_S + A \cdot I_C$.  Hence $I_S + AB \cdot J$
links $I_Y$ to some ideal $\fb$ which is unmixed.

\bigskip

\noindent \underline{\em Step IV}: One can check that $\fa \subset \fb$
and compute that $\deg \fa = \deg \fb$.  Since both are unmixed of the
same dimension, it follows that $\fa = \fb$.  Hence $I_C$ is G-linked to
$I_Y$ in two steps.
\hfill \qed

\begin{remark} \label{basic double CI-link}
A special case of Proposition \ref{bdl is linked} is worth mentioning.
Suppose that $S$ is a complete intersection, $I_S = (F_1,\dots,F_c)$,
and $I_Y = A \cdot I_C + (F_1,\dots,F_c)$.  Then all the links in
Proposition \ref{bdl is linked} are complete intersections.  This
construction is called {\em Basic Double CI-Linkage}; cf.\ \cite{LR},
\cite{BM4}, \cite{GM4}.  As an even more special case, suppose that $c = 1$
and hence $\codim C = 2$.  Let $F \in I_C$ and assume that $A,F$ have no
common factor.  Then $I_Y = A \cdot I_C + (F)$.  This construction is central
to the Lazarsfeld-Rao property, which we will discuss below.  This property
is only known in codimension two.

A different way of viewing Basic Double Linkage, as a special case of Liaison
Addition, will be discussed next.
\end{remark}

Liaison Addition was part of the Ph.D.\ thesis of Phil Schwartau
\cite{schwartau}.  The problem which he originally considered was the
following.  Consider curves $C_1, C_2 \subset \proj{3}$.  Suppose that $H_*^1
({\mathcal I}_{C_1}) = M_1$ and $H_*^1 ({\mathcal I}_{C_2}) = M_2$.  Find an
explicit construction of a curve $C$ for which $H^1_*({\mathcal I}_C) = M_1
\oplus M_2$.

The first observation to make is that this is impossible in general!!
We give a simple example.

\begin{example}
Let $C_1$ and $C_2$ be two disjoint sets of two skew lines.  We have noticed
(Example \ref{ex-min-shift}) that
$H^1_*({\mathcal I}_{C_1}) = H^1_*({\mathcal I}_{C_2}) = K$, a graded module
of dimension 1 occurring in degree~0.  So the question is whether there
exists a curve $C$ with $H^1({\mathcal I}_C) = K^2$, a 2-dimensional module
occurring in degree 0.  Suppose that such a curve exists.  Let $L$ be a
general linear form defining a hyperplane $H$.  We have the long exact
sequence
\[
\begin{array}{ccccccccccccccccccc}
0 & \! \rightarrow & \! H^0({\mathcal I}_C) & \! \rightarrow & \!
H^0({\mathcal I}_C(1)) & \! \rightarrow & \! H^0({\mathcal I}_{C \cap H | H}
(1)) & \!
\rightarrow & \! H^1({\mathcal I}_C) & \! \rightarrow & \! H^1({\mathcal
I}_C(1)) & \! \rightarrow \dots \\
&& || && || &&&& || && || \\
&& 0 && 0 &&&& 2 && 0
\end{array}
\]
Note that $h^0({\mathcal I}_C (1)) = 0$ since otherwise $C$ would be a plane
curve, hence a hypersurface and hence \acm.  But this means that $C$ is a
curve whose general hyperplane section lies on a pencil of lines in
$\proj{2}$.  This forces $\deg C = 1$, hence $C$ is a line and is thus \acm.
Contradiction.
\end{example}

However, an important idea that we have seen in Section \ref{prelim result
section} is that the {\em shift} of the modules is of central importance.
Hence the refined problem that Schwartau considered is whether there is a
construction of a curve $C$ for which $H_*^1 ({\mathcal I}_C) = M_1 \oplus
M_2$ {\em up to shift}.  As we will see, he was able to answer this question
and even a stronger one (allowing the modules to individually have different
shifts), and his work was for codimension two in general.  The version that
we will give is a more general one, however, from \cite{GM4}.  The statement,
but not the proof, were inspired by \cite{schwartau}, which proved the case
$r=2$.

\begin{theorem}
Let $V_1,\dots,V_r$ be closed subschemes of $\proj{n}$, with $2
\leq r \leq n$.  Assume that $\codim V_i \geq r$ for all $i$.
Choose homogeneous polynomials
\[
\hbox{\hskip 1in}
F_i \in \bigcap_{
\begin{array}{c}
  1 \leq j \leq r \\ j \neq i
\end{array}}
I_{V_j} \hbox{\hskip 1in for $1 \leq i \leq r$}
\]
such that $(F_1,\dots,F_r)$ form a regular sequence, hence defining a
complete intersection, $V \subset \proj{n}$.  Let $d_i = \deg F_i$.  Define
the ideal $I = F_1 \cdot I_{V_1} + \dots + F_r \cdot I_{V_r}$.  Let $Z$ be
the closed subscheme of $\proj{n}$ defined by $I$ (which a priori is not
saturated).  Then

\begin{itemize}
\item[(a)] As sets, $Z = V_1 \cup \dots \cup V_r \cup V$;

\item[(b)] For all $1 \leq j \leq n-r = \dim V$ we have
\[
H^j_*({\mathcal I}_Z) \cong H^j_*({\mathcal I}_{V_1})(-d_1) \oplus \dots
\oplus H^j_*({\mathcal I}_{V_r})(-d_r);
\]

\item[(c)] $I$ is saturated.
\end{itemize}
\end{theorem}

We will not give the proof of this theorem, but refer the reader to
\cite{GM4} or to \cite{migbook}.

\begin{remark}\label{basic double ci linkage}
Note that we allow $V_1,\dots,V_r$ to be of different codimensions, we allow
them to fail to be locally Cohen-Macaulay or equidimensional, and we even
allow them to be empty.  In fact, this latter possibility gives another
approach to Basic Double CI-Linkage (cf.\ Remark \ref{basic double
CI-link}).  Indeed, if we let $V_2 = \dots = V_r = \emptyset$, with $I_{V_2}
= \dots = I_{V_r} = R$, and set $V_1 = C$ and $V = S$ as in Remark \ref{basic
double CI-link}, then the ideal
\[
\begin{array}{rcl}
I & = & F_1 \cdot I_{V_1} + \dots + F_r \cdot I_{V_r} \\
& = & F_1 \cdot I_{V_1} + (F_2,\dots,F_r)
\end{array}
\]
is precisely the ideal of the Basic Double CI-Linkage.
\end{remark}

An application of Liaison Addition is to the construction of arithmetically
Buchsbaum curves, or more generally arithmetically Buchsbaum subschemes of
projective space.  We will give the basic idea here and come back to it in
Section \ref{application section}.

\begin{definition}
A curve $C \subset \proj{n}$ is {\em \aBM} if $H^1_*({\mathcal I}_C)$ is
annihilated by the ``maximal ideal'' $\fm = (x_0,\dots,x_n)$.  A subscheme $V
\subset \proj{n}$ of dimension $d \geq 2$ is \aBM\ if $H^i_*({\mathcal I}_V)$
is annihilated by $\fm$ for $1 \leq i \leq d$ and furthermore the general
hyperplane section $V \cap H$ is an \aBM\ scheme in $H = \proj{n-1}$.
\end{definition}

Buchsbaum curves in $\proj{3}$, especially, are fascinating objects which
have been studied extensively.  A rather large list of references can be
found in \cite{migbook}.  Liaison Addition can be used to construct examples
of \aBM\ curves with modules whose components have any prescribed
dimensions, up to shift.  Indeed, Schwartau's original work \cite{schwartau}
already produced examples of modules of any dimension, concentrated in one
degree.  The more general result was obtained in \cite{BM2}.  The idea is to
use sets of two skew lines as a ``building block'' to build bigger curves.  We
will give the basic idea with an example, omitting the proof of the general
result.

\begin{example}\label{construction of aBM curves}
Recall (Example \ref{ex-min-shift}) that the deficiency module of a set of two
skew lines is one-dimensional as a $K$-vector space, and the non-zero
component occurs in degree 0.  Furthermore, this is the minimal shift for
that module.  Let $C_1$ and $C_2$ be two such curves.  How $C_1$ and $C_2$
meet (i.e.\ whether they are disjoint from each other, meet in finitely
many points
or contain common components) is not important.  Let
$F_1 \in I_{C_2}$ and $F_2 \in I_{C_1}$ such that $F_1$ and $F_2$ have no
common factor.  Note that $\deg F_i \geq 2$.  Then the curve $C$ obtained by
$I_C = F_1 \cdot I_{C_1} + F_2 \cdot I_{C_2}$ is \aBM\ and its deficiency
module is the direct sum of twists of the deficiency modules of $C_1$ and
$C_2$.  In particular, this module is 2-dimensional as a
$K$-vector space.  The components occur in degrees $\deg F_1$ and $\deg
F_2$.  In particular, the components can be arbitrarily far apart, and
regardless of how far apart they are, the leftmost component occurs in degree
$\geq 2$.  Furthermore, an example can be obtained for which this component is
in degree exactly 2 (for instance by choosing $F_1$ with $\deg F_1 = 2$ and
then $F_2$ of appropriate degree).  Note that $2 = 2 \cdot 2 - 2$ (see
Proposition \ref{buchs leftmost shift}).

An \aBM\ curve whose module is 3-dimensional as a $K$-vector space can then
be constructed by taking the Liaison Addition of $C$ with another set of two
skew lines, and it is clear that this process can be extended to produce any
module (up to shift) which is annihilated by $\fm$.  Furthermore, if we
produce in this way a curve whose deficiency module has dimension $N$ as a
$K$-vector space, a little thought shows that this can be done in
such a way that the leftmost non-zero component occurs in degree $\geq 2N-2$,
and that a sharp example can be constructed.  We refer to \cite{BM2} for
details.
\end{example}

This approach was also used in \cite{GM4} to construct \aBM\ curves in
$\proj{4}$, and in \cite{BM5} to construct certain \aBM\ surfaces in
$\proj{4}$ with nice properties.

One would like to have an idea of ``how many'' of the Buchsbaum curves can be
constructed using this approach together with Basic Double Linkage (which
preserves the module but shifts it to the right, adding a complete
intersection to the curve).  The first step was obtained in \cite{GM1}:

\begin{proposition} \label{buchs leftmost shift}
Let $C \subset \proj{3}$ be an \aBM\ curve.  Let
\[
N = \dim H^1_*({\mathcal
I}_C) = \sum_{t \in {\mathbb Z}} h^i({\mathcal I}_C(t)).
\]
   Then the first non-zero component of $H^1_*({\mathcal I}_C)$ occurs in
degree $\geq 2N-2$.
\end{proposition}

\begin{proof}
The proof is an easy application of a result of Amasaki (\cite{amasaki},
\cite{GM2}) which says that $C$ lies on no surface of degree $< 2N$.  We
refer to \cite{GM1} for the details.
\end{proof}

It follows that the construction of Example \ref{construction of aBM curves}
provides curves which are in the minimal shift for their module.  We will see
below that the Lazarsfeld-Rao property then gives an incredible amount of
information about all \aBM\ curves, once we know even one curve in the
minimal shift.  It will turn out that this construction together with Basic
Double Linkage, gives {\em all} \aBM\ curves in $\proj{3}$ up to deformation.
Furthermore, this construction will even give us information
about \aBM\ stick figures.


\section{Necessary conditions for being linked} \label{sec-necess}

Since (direct) linkage is symmetric, the transitive closure of this
relation generates an equivalence relation, called liaison.
However, it will be useful to study slightly different
equivalence classes.

As in the previous sections we will restrict ourselves to
subschemes of $\PP^n$ and ideals of $R = K[x_0,\ldots,x_n]$,
although the results are more generally true for subschemes of an
\aG\ scheme.

\begin{definition} \label{def-even-li-cl}
Let $I \subset R$ denote an unmixed ideal. Then the {\it even
G-liaison class} of $I$ is the set
$$
\cL_I = \{ J \subset R \s I \stackrel{\fc_1}{\sim} I_1
\stackrel{\fc_2}{\sim} \ldots \stackrel{\fc_{2k}}{\sim} J \}
$$
where $\fc_1,\ldots,\fc_{2k}$ are Gorenstein ideals.
If we require that all ideals $\fc_1,\ldots,\fc_{2k}$ are complete
intersections that we get the {\it even CI-liaison class} of $I$.

The even G-liaison class and the even CI-liaison class of an
equidimensional subscheme $V \subset \PP^n$ are defined
analogously.
\end{definition}

\begin{remark} \label{rem-even-cl}
It is clear from the definition that every G-liaison class
consists of at most two even G-liaison classes and that every
CI-liaison class consists of at most two even CI-liaison classes.

A liaison class can agree with an even liaison class, for example,
if it contains a self-linked element.  On the other hand, it was shown by Rao
\cite{rao3} that there are liaison classes that coincide with even liaison
classes
but contain no self-linked elements, the simplest being the liaison class
of two
skew lines in $\proj{3}$.
\end{remark}

The next result has been shown in various levels of generality by
Chiarli, Schenzel, Rao, Migliore.

\begin{lemma} \label{lem-coho-lCM-cl}
If $V \subset \PP^n$ is an equidimensional \lCM\ subscheme and $W
\in \cL_V$ then there is an integer $t$ such that
$$
H^i_*(\cI_W) \cong H^i_*(\cI_V) (t) \quad \mbox{for all} \;\; i =
1,\ldots,\dim V.
$$
\end{lemma}

\begin{proof}
This follows immediately from the comparison of the cohomology of
directly linked \lCM\ subschemes (Corollary
\ref{cor-dir-linked-id}) because we have for every finitely
generated graded $R$-module $M$ that $(M^{\vee})^{\vee} \cong M$.
\end{proof}

Our next goal is to show that there is a stronger result which is
true even if $V$ is not \lCM. For this we have to consider certain
types of exact sequences. The names have been coined by
Martin-Deschamps and Perrin \cite{MP}.

\begin{definition} \label{def-E-N-type}
Let $I \subset R$ be a homogeneous ideal of codimension $c \geq 2$.

(i)   An {\it $E$-type resolution}
   of $I$ is  an exact sequence of finitely generated graded $R$-modules
$$
0 \to E \to F_{c-1} \to \ldots \to F_1 \to I \to 0
$$
where the modules $F_1, \ldots, F_{c-1}$ are free.

(ii) An {\it $N$-type resolution} of $I$ is an exact sequence of finitely
generated graded $R$-modules
$$
0 \to G_c \to \ldots \to G_2 \to N \to I \to 0
$$
where $G_2, \ldots, G_c$ are free $R$-modules and $\HH^i(N) = 0$ for all
$i$ with $n+2-c
\leq i \leq n$.
\end{definition}

\begin{remark} \label{rem-E-N-type}
(i) The existence of an $E$-type resolution is clear because it is
just the beginning of a free resolution of $I$. For the existence
of an $N$-type resolution we refer to \cite{N-gorliaison}.
However, we will see that for an unmixed ideal the existence
follows by liaison.

(ii) It is easy to see that
$$
\HH^i(N) \cong \HH^{i-1}(R/I) \quad \mbox{for all}\;\; i \leq n+1-c
$$
and
$$
\HH^i(E) \cong \HH^{i-c}(R/I) \fora i \leq n.
$$
This shows that the modules $E$ and $N$ ``store'' the deficiency
modules of $R/I$.

(iii) The sheafifications $\tilde{E}$ and $\tilde{N}$ are vector
bundles if and only if $I$ defines an equidimensional \lCM\
subscheme.
\end{remark}

By definition, the module $E$ is a $c$-syzygy. It is not to
difficult to check that $N$ must be a torsion-free module.
However, if the ideal $I$ is unmixed then these modules have
better properties.

\begin{lemma} \label{lem-unmixed-by-E-N}
Using the notation of Definition \ref{def-even-li-cl} the
following conditions are equivalent:
\begin{itemize}
\item[(a)] The ideal $I$ is unmixed.
\item[(b)] The module $N$ is reflexive, i.e.\ the bilinear map $N
\times N^* \to R, (m, \ffi) \mapsto \ffi(m)$, induces an
isomorphism $N \to N^{**}$.
\item[(c)] The module $E$ is a $(c+1)$-syzygy.
\end{itemize}
\end{lemma}

\begin{proof}
This follows from the cohomological characterization of these
concepts.
\end{proof}

The next result establishes the crucial fact that $E$-type and
$N$-type resolutions are interchanged under direct linkage.

\begin{proposition} \label{E-and-N-resolution-under-linkage} Let $I, J
   \subset R$ be homogeneous ideals of codimension $c$ which are directly
   linked by
   $\fc$. Suppose $I$ has resolutions of $E$- and
   $N$-type as in Definition \ref{def-even-li-cl}. Let
$$
0 \to D_c \to \ldots \to D_1 \to \fc \to 0
$$
be a minimal free resolution of $\fc$.  Put $s = r(R/\fc) + n$.  Then
$J$ has an $N$-type resolution
$$
0 \to F_1^*(-s) \to D_{c-1}  \oplus F_2^*(-s) \to \ldots \to D_2 \oplus
F_{c-1}^*(-s)
\to D_1 \oplus E^*(-s) \to J \to 0
$$
and an $E$-type resolution
$$
0 \to N^*(-s) \to D_{c-1} \oplus G_2^*(-s) \to \ldots \to D_1 \oplus
G_c^*(-s) \to J \to 0.
$$
\end{proposition}

\begin{proof}
  We want to produce an $N$-type resolution of $J$. We
proceed in several steps.

(I) Resolving $E$ we get an exact sequence
$$
\begin{array}{rcl}
\ldots \to F_{c+1} \stackrel{\ffi_{c+1}}{\longrightarrow} F_c &
\stackrel{\ffi_{c}}{\longrightarrow} & F_{c-1} \to \ldots \to F_1
\stackrel{\ffi_{1}}{\longrightarrow} I \to 0. \\
& \searrow \quad \nearrow & \\
& E & \\
& \nearrow \quad \searrow & \\
0 & & 0
\end{array}
$$
Dualizing with respect to $R$ gives a complex
$$
0 \to R \to F_1^* \to \ldots \to F_{c-1}^*
\stackrel{\ffi_{c}^*}{\longrightarrow}
F_c^* \stackrel{\ffi_{c+1}^*}{\longrightarrow} F_{c+1}^*
$$
and an exact sequence
$$
0 \to E^* \to F_c^* \stackrel{\ffi_{c+1}^*}{\longrightarrow} F_{c+1}^*.
$$
If follows $\ker \ffi_{c+1}^* \cong E^*$, thus $\ER^c(R/I,R) \cong E^*/\im
\ffi_c^*$.  Moreover, we get by duality that
$$
\ker \ffi_{i+1}^*/ \im \ffi_i^* \cong \ER^i(R/I,R) \cong
\HH^{n+1-i}(R/I)^{\vee}(1-r(R)) = 0 \quad \mif i < c
$$
because $\dim R/I = n+1-c$. Therefore we can splice the two complexes above
together and the resulting diagram
$$
\begin{array}{rcl}
0 \to R \to F_1^* \to \ldots \to F_{c-1}^* & \longrightarrow & E^* \to
\ER^c(R/I,R) \to 0 \\
& \searrow \quad \nearrow & \\
& \im \ffi_c^* & \\
& \nearrow \quad \searrow & \\
0 & & 0
\end{array}
$$
is exact.

(II) The self-duality of the minimal free resolution of $R/\fc$ and
Corollary \ref{cor-shift} provide $D_c = R(-s)$ and $D_{c-i}^*
\cong D_i (s)$ for all $i = 1,\ldots, c-1$.

(III) Put $r := r(R/\fc) -1$. The standard exact sequence
provides
the following diagram with exact rows and column:
$$
\begin{array}{rcl}
& 0 &\\
& \downarrow & \\
& K_{R/J}(-r) & \\
& \downarrow & \\
0 \to R(-s) \to D_{c-1} \to \ldots \to D_1 \to R \to & R/\fc & \to 0 \\
& \downarrow & \\
0 \to E \to F_{c-1} \to \ldots \to F_1 \to R \to & R/I & \to  0 \\
& \downarrow & \\
& 0. &
\end{array}
$$
Since the modules $D_1, \ldots , D_{c-1}$ are free the epimorphism $R/\fc
\to R/I$ lifts to a morphism of complexes. Thus, using steps (I) and  (II)
we get by
dualizing with respect to $R$ the commutative exact diagram:
$$
\begin{array}{ccccccccrcl}
& & & & & & & & & {\ER^{c-1}(K_{R/J},R)(-r)} &\\
& & & & & & & & & \downarrow & \\
0 \to &  R & \to & F_1^* & \to \ldots \to & F_{c-1}^* & \to & E^* & \to &
\ER^c(R/I,R) & \to 0 \\
& \| & & \downarrow & & \downarrow & & \downarrow & & \downarrow\
{\scriptstyle \alpha} &
\\
0 \to & R & \to & D_{c-1}(s) & \to \ldots \to & D_1(s) & \to & R(s) & \to &
\ER^c(R/\fc,R) & \to 0.
\end{array}
$$
Since $K_{R/J}$ has dimension $n+1-c$ we obtain  by duality $
\ER^{c-1}(K_{R/J},R) = 0$. Moreover, we have already seen that
$\ER^c(R/I,R) \cong K_{R/I}(n+1)$ and $\ER^c(R/\fc,R) \cong
R/\fc(s)$.  It follows
that $\alpha$ is injective and, by comparison with the standard exact
sequence,  $\coker \alpha \cong R/J(s)$. Thus,
the mapping cone procedure provides an exact sequence which begins
with $R/J(s)$ and ends with $R$. However, it can be shown that
this last module can be canceled. The result is  the $N$-type
resolution of $J$ as claimed because $E^*$ meets the cohomological
requirements (cf.\ \cite{Auslander-Bridger}, Theorem 4.25 and
\cite{EG_buch}, Theorem 3.8).

The claimed $E$-type resolution of $J$ can be obtained by similar
arguments.
\end{proof}

Using Remark \ref{rem-E-N-type} we obtain as first consequence the
generalization of Lemma \ref{lem-coho-lCM-cl}.

\begin{corollary} \label{cor-coho-even}
Let $V, W \subset \PP^n$ be equidimensional subschemes. If $W \in
\cL_V$ then
there is an integer $t$ such that
$$
H^i_*(\cI_W) \cong H^i_*(\cI_V) (t) \fora i = 1,\ldots,\dim V.
$$
\end{corollary}

This result gives necessary conditions for $V$ and $W$ being in the same even
liaison class. In order to state a stronger consequence of Proposition
\ref{E-and-N-resolution-under-linkage} we need.

\begin{definition} \label{def-min-E-N-res}
The $E$-type and the $N$-type resolution, respectively, of $I$ are said to {\it
minimal}
if it is not possible to cancel free direct summands. They are uniquely
determined
by $I$ up to isomorphism of complexes.

Let $\ffi(I)$ denote the last non-vanishing module in a minimal $E$-type
resolution
of $I$, and let $\psi(I)$ denote the second non-vanishing module in a minimal
$N$-type resolution of $I$.
\end{definition}

We consider $\ffi$ and $\psi$ as maps from the set of ideals into the set
of maximal
$R$-modules.
If the $E$-type and $N$-type resolutions in Definition \ref{def-E-N-type} are
minimal then we have $\ffi(I) = E$ and $\psi(I) = N$. Basically we will
ignore possible free direct summands of $\ffi(I)$ and $\psi(I)$. This is
formalized as
follows.

\begin{definition} \label{def-stable-equiv}
Two  graded maximal $R$-modules $M$ and $N$ are said to be {\itshape stably
equivalent} if there are finitely generated, free $R$-modules $F, G$ and an
integer $t$ such that
$$
M \oplus F \cong N(t) \oplus G.
$$
The {\it stable equivalence class} of $M$ is the set
$$
[M] := \{ N \s N \ \mbox{is stably equivalent to} \ M \}.
$$
\end{definition}

Using Proposition \ref{E-and-N-resolution-under-linkage} repeatedly  we get the
following relation between even liaison and certain stable equivalence classes.

\begin{theorem}[Rao's correspondence] \label{thm-rao-corres}
The map $\ffi$ induces a well-defined  map $\Phi: \cM_c \to \cM_E,
\cL_I \mapsto [\ffi(I)],$
from the set $\cM_c$ of even
   liaison classes of unmixed ideals in $R$ of codimension $c$ into the set
$\cM_E$
   of stable
   equivalence classes of finitely generated $(c+1)$-syzygies.

The map $\psi$ induces a well-defined map $\Psi: \cM_c \to \cM_N,
\cL_I \mapsto [\psi(I)],$
from the set $\cM_c$ of even
   liaison classes   into the set $\cM_N$ of stable
   equivalence classes of finitely generated reflexive modules $N$ which
   satisfy $\HH^i(N) = 0$ for all $i$ with $n-c+2 \leq i \leq n$.
\end{theorem}

\begin{remark} \label{rem-rao-corr}
(i) Rao's correspondence provides the following  diagram with two
commuting squares $$
\begin{array}{rcl}
\cM_c & \stackrel{\Phi}{\longrightarrow} & \cM_E \\
\downarrow  {\scriptstyle \alpha} & & \downarrow {\scriptstyle \beta} \\
\cM_c & \stackrel{\Psi}{\longrightarrow} & \cM_N \\
\downarrow {\scriptstyle \alpha} & & \downarrow {\scriptstyle \beta} \\
\cM_c & \stackrel{\Phi}{\longrightarrow} & \cM_E
\end{array}
$$ where $\alpha$ is induced by direct linkage and $\beta$ is
induced by dualization with respect to $R$.

(ii) Combining Rao's correspondence with Horrocks' classification
of stable equivalence classes of vector bundles on $\PP^n$  in
terms of cohomology groups and extensions \cite{Horrocks-corr}
gives a stronger result than Corollary \ref{cor-coho-even} in the
case of \lCM\ subschemes. Horrocks' result gives for example: if
$C$ is a curve then the stable equivalence class of $\psi(I_C)$ is
determined by $H^1_*(\cI_C) \cong \HH^2(\psi(I_C))$. For a \lCM\
surface $S \subset \PP^n$ the stable equivalence class of
$\psi(I_S)$ is determined by the triple $(\HH^2(\psi(I_S)),
\HH^3(\psi(I_S)), \eta)$ where $\eta \in
\Ext^2_R(\HH^3(\psi(I_S)), \HH^2(\psi(I_S)))$. For the modules
associated to schemes of dimension $\geq 3$ Horrocks'
classification becomes less elegant.

In particular, the modules $H^1_*(\cI_S)$ and $H^2_*(\cI_S)$ are not
enough to determine the even liaison class.  This is illustrated, for
instance, in  Example \ref{ex-rao-corr}.
\end{remark}

The next result gives more information on Rao's correspondence.

\begin{proposition} \label{prop-Rao-surj}
For every $c \geq 2$ the maps $\Phi$ and $\Psi$ occurring in Rao's
correspondence
are surjective.
\end{proposition}

\begin{proof}
Fixing $c$  it suffices to show the claim for one of the maps
according to the previous remark.

Let $c=2$. Let $M \in \cM_N$ be a module of rank $r$. Then for $s \gg 0$ a
sufficiently general map $R^{r-1}(-s) \to M$ provides an exact sequence
$$
  0 \to R^{r-1}(-s) \to M \to I(t) \to 0
$$
where $t$ is an integer and $I$ is an unmixed ideal of codimension two.
This result
is sometimes referred to as Theorem of Bourbaki.

For $c \geq 3$ the claim is shown in \cite{amasaki-pams}.
\end{proof}

Rao's correspondence gives the strongest known necessary conditions for two
subschemes
belonging to the same even G-liaison class. For even CI-liaison classes of
ideals of
codimension $c \geq 3$ there are
additional necessary conditions (cf.\ \cite{HU-Annals},  \cite{KMMNP}).
\smallskip

The next example illustrates the fact that Rao's correspondence provides
stronger
necessary conditions than Corollary \ref{cor-coho-even}.

\begin{example} \label{ex-rao-corr}
The Koszul complex resolves the ideal $\fm = (x_0,\ldots,x_4)$ over $R :=
K[x_0,\ldots,x_4]$
$$
\begin{array}{rcccl}
0 \to R(-5) \to R^5(-4) & \longrightarrow & R^{10}(-3) \to R^{10}(-2) &
\longrightarrow & R^5(-1) \to \fm \to 0. \\
& \searrow \quad \nearrow & &  \searrow \quad \nearrow & \\
& \Omega^3 &  &  \Omega^1 & \\
& \nearrow \quad \searrow &  & \nearrow \quad \searrow &\\
0 & & 0 \hfill 0 & & 0
\end{array}
$$
The modules $\Omega^3$ and $\Omega^1$ are defined as the indicated syzygy
modules.

There is a surface $S \subset \PP^4$ admitting an exact sequence
$$
0 \to (\Omega^3(-1))^2 \stackrel{\alpha}{\longrightarrow} R(-4) \oplus
(\Omega^1(-3))^2 \to I_S \to 0.
$$
If the map $\alpha$ is general enough then $S$ is a smooth rational surface
of degree
$10$ (\cite{DES}, Example B1.15). Moreover, its deficiency modules are
$$
\begin{array}{rcl}
H^1_*(\cI_S) & \cong & K^2(-3) \\
H^2_*(\cI_S) & \cong & K^2(-1).
\end{array}
$$
Denote by $\Omega^2$ the second syzygy module of $\fm$ in the Koszul
complex above.
There is a surface $T \subset \PP^4$ such that there is an exact sequence
$$
0 \to R^{15} (-4) \to (\Omega^2)^2 \oplus (\Omega^1(-2))^2 \to I_T \to 0.
$$
The deficiency modules of $T$ are
$$
\begin{array}{rcl}
H^1_*(\cI_T) & \cong & K^2(-2) \\
H^2_*(\cI_T) & \cong & K^2.
\end{array}
$$
Thus, we obtain
$$
\begin{array}{rcl}
H^i_*(\cI_S) & \cong & H^i_*(\cI_T) (-1) \quad \mbox{for} \; \; i = 1,2.
\end{array}
$$
This leaves open the possibility that $S$ and $T$ are evenly linked. But in
fact,
$S$ and $T$ belong to different even liaison classes because $\ffi(I_S)$ and
$\ffi(I_T)$ are not stably equivalent (\cite{N-characterization}, Example
7.4). Taking
concrete examples for $S$ and $T$ this can be checked by looking at their
hyperplane
sections. Let $H \subset \PP^4$ be a general hyperplane. It is not
difficult to see
that $S \in \cL_T$ would imply $S \cap H \in \cL_{T \cap H}$, but the latter is
impossible because
$$
\fm \cdot H^1_*(\cI_{T \cap H}) = 0, \quad \mbox{but} \quad \fm \cdot
H^1_*(\cI_{S \cap H}) \neq 0.
$$
In other words, the surface $T$ is \aBM, but the surface $S$ is not \aBM. Using
Rao's correspondence one can show that the property of being \aBM\ is
preserved under
direct linkage (cf.\ \cite{schenzel}. \cite{N-characterization}) which
again implies
$S \notin \cL_T$.
\end{example}


\section{Sufficient conditions for being linked} \label{sec-suff}

In Section \ref{sec-necess} we have seen that Rao's correspondence relates even
liaison classes to certain stable equivalence classes. Moreover, this
correspondence
is surjective. Thus, an ideal result would be an affirmative answer to the
following.

\begin{mquestion} \label{mque-inj}
Are the maps $\Phi$ and $\Psi$ in Rao's correspondence injective for all $c
\geq 2$?
\end{mquestion}

In this section we will discuss this question.
It is worthwhile to point out special cases of the Main question.  We begin
with
a definition.

\begin{definition}
A subscheme $V \subset \proj{n}$ is {\em licci} if it is in the CI-liaison
class
of a complete intersection.  $V$ is {\em glicci} if it is in the G-liaison
class
of a complete intersection.
\end{definition}

\begin{remark} \label{rem-spec-cases}
(i) Since for an \aCM\ subscheme $V$ the modules $\Phi(I_V)$ and
$\Psi(I_V)$ are free, i.e.\ stably equivalent to the zero module,
in this case the Main question takes the form:

{\bf Question 1:} Is it true that a subscheme $V$ is \aCM\
if and only if it is {glicci}, i.e.\ in the G-liaison class of
a complete intersection?

(ii) Let $C, D \subset \PP^n$ be two curves. Then it is not to
difficult to see, and it is a special case of  Horrocks' results
\cite{Horrocks-corr}, that the following conditions are equivalent:
\begin{itemize}
\item[(a)] $\Phi(I_C)$ and $\Phi(I_D)$ are stably equivalent.
\item[(b)] $\Psi(I_C)$ and $\Psi(I_D)$ are stably equivalent.
\item[(c)] $H^1_*(\cI_C) \cong H^1_*(\cI_D) (t)$ for some $t \in
\ZZ$.
\end{itemize}
Hence, for curves the main question specializes to:

{\bf Question 2:} Is it true that two curves $C, D \subset \PP^n$
belong to the same even G-liaison class if and only if
$H^1_*(\cI_C) \cong H^1_*(\cI_D) (t)$ for some $t \in \ZZ$?

(iii) Strictly speaking we should ask the Main question and the
two questions above for G-liaison and CI-liaison, separately. We
state only  the Main question for  CI-liaison.

{\bf Question 3:} Let $V, W \subset \PP^n$ be two equidimensional
subschemes of the same codimension. Is it true that $V$ and $W$
belong to the same even CI-liaison class if and only if
$\Phi(I_V)$ and $\Phi(I_W)$ are stably equivalent?
\end{remark}

For subschemes of codimension two the answer to all these
questions is `yes' because of the following result which is
essentially due to Rao \cite{rao1} (cf.\ also \cite{nollet},
\cite{N-gorliaison}).

\begin{theorem} \label{thm-rao-inj}
Let $I,I' \subset R$ be unmixed homogeneous ideals of
   codimension $2$ with $N$-type resolutions
$$
0 \to \bigoplus_{i=1}^s R(-a_i) \stackrel{\delta}{\longrightarrow}
N \to I \to 0
$$
and
$$
0 \to \bigoplus_{i=1}^s R(-b_i) \stackrel{\veps}{\longrightarrow}
N(h) \to I' \to 0. $$ Then $I$
and $I'$ belong to the same even liaison class.
\end{theorem}

\begin{proof} We only outline the proof but give enough details to
see where the problems are in extending the argument for ideals of
higher codimension.

{\it Case 1:} Suppose $N$ is a free $R$-module. Then $I$ and $I'$
are standard determinantal ideals and the claim follows from the
more general Theorem \ref{thm-gaeta}.

{\it Case 2:} Suppose that $N$ is not free. Then, possibly after
linking $I$ and $I'$ in an even number of steps to new ideals, we
may assume that $N$ does not have a free direct summand.

Write $\im \delta = (m_1,\ldots,m_s)$ and $\im \veps =
(n_1,\ldots,n_s)$ where $m_i R = \delta(R(-a_i))$ and $n_i R =
\veps(R(-b_i))$. Suppose that $m_i = n_i$
for $i < t \leq s$.

We want to show that we can find ideals $I_1 \in \cL_I$ and
$I_1' \in \cL_{I'}$  having $N$-type
resolutions where $m_i = n_i$ for $i \leq t$. Then, repeating this
process at most $s$ times  our statement follows.

Choose an integer $p \gg 0$ and elements $u,v \in [N]_p$  whose images
$f_1,f_2$ in $I$
and $g_1,g_2$ in $I'$ generate complete intersections $\fc$ and $\fc'$,
respectively.
Put $J = \fc : I$ and $J' = \fc' : I'$. According to
Proposition \ref{E-and-N-resolution-under-linkage} these ideals have
$E$-type resolutions as follows:
$$
0 \to N^*(-2p) \to R^2(-p) \oplus \bigoplus_{i=1}^s R(a_i-2p) \to J \to
0,
$$
$$
0 \to N^*(2h-2p) \to R^2(h-p) \oplus \bigoplus_{i=1}^s R(b_i+2h-2p) \to J' \to
0.
$$
Let $f \in J$ be the generator of the image of $R(a_t-2p)$ in $J$ and let
$g \in J'$ be the generator of the image of $R(b_t+2h-2p)$; $f$ and $g$
are not zero because $N$ does not have a free direct summand.  Since
$\{f_1,f_2\}$ and $\{g_1,g_2\}$ are regular sequences it is possible to
find $\lambda,\mu \in K$ such that $\fd = (f,f')$ and $\fd' = (g,g')$ are
complete intersections where $f' =
\lambda f_1 + \mu f_2$ and $g' = \lambda g_1 + \mu g_2$. Put $I_1 = \fd :
J$ and $I_1' = \fd' : J'$.
Since $N$ does not have a free direct summand the  $E$-type resolutions
   of $J$
and $J'$ above must be minimal. It follows that $f,f'$ are minimal
generators of
$J$ and that $g,g'$ are minimal generators of $J'$. Therefore we can split
off $R(-p) \oplus R(a_t-2p)$ respectively $R(h-p) \oplus R(b_t+2h-2p)$ in
the $N$-type resolution of $I_1$ respectively $I_1'$ given by Proposition
\ref{E-and-N-resolution-under-linkage}. The resulting resolutions are:
$$
0 \to \bigoplus_{i \neq t} R(a_t-a_i-p) \oplus R(a_t-2p)
\stackrel{\alpha}{\longrightarrow} N(a_t-p) \to
I_1 \to 0,
$$
$$
0 \to \bigoplus_{i \neq t} R(b_t-b_i+h-p) \oplus R(b_t+2h-2p)
\stackrel{\alpha'}{\longrightarrow} N(b_t+h-p) \to
I_1' \to 0
$$
where the image of $\alpha$ is generated by
$m_1,\ldots,m_{t-1},m_{t+1},\ldots,m_s,\lambda u + \mu v$ and the image of
$\alpha'$ is generated by $n_1,\ldots,n_{t-1},n_{t+1},\ldots,n_s,\lambda u +
\mu v$. This means that we have reached our goal by replacing $m_t$ and $n_t$,
respectively,  by $\lambda u + \mu v$.
\end{proof}

Let us look at this proof in case  the codimension of $I$ and $I'$
is at least three. We can still split off terms in the $N$-type resolutions of
$I_1$and $I'_1$ at the end of these resolutions (as in the proof
above). But since the resolutions are longer this is not enough to
guarantee the splitting at the beginning of the resolution which
would be needed to complete the argument.
\smallskip

As pointed out in Remark \ref{rem-spec-cases}, the last result has an
implication for space curves.

\begin{corollary} \label{cor-li-sp-curves}
Let $C, D \subset \PP^3$ be two curves. Then $D \in \cL_C$
  if and only if
$H^1_*(\cI_C) \cong H^1_*(\cI_D) (t)$ for some $t \in \ZZ$.
\end{corollary}

We still have to prove Case 1 in the previous proposition, which is
a result of Gaeta. It can be generalized to arbitrary codimension
\cite{KMMNP}. For this we recall the following.

\begin{definition} \label{def-stand-det}
If $A$ is a homogeneous matrix, we denote by $I(A)$ the ideal of
maximal minors of $A$.  If $\varphi :F \rightarrow G$ is a
homomorphism of free graded $R$-modules then we define $I(\varphi)
= I(A)$ for any homogeneous matrix $A$ representing $\varphi$
after a choice of basis for $F$ and $G$. A codimension $c+1$ ideal
$I \subset R$
will be called a {\em standard determinantal ideal} if $I_X =
I(A)$ for some homogeneous $t \times (t+c)$ matrix, $A$. In a
similar way we define standard determinantal subscheme of $\PP^n$.
\end{definition}

It is well-known that every standard determinantal subscheme is
\aCM.

Now we can state one of the main results of \cite{KMMNP}.  The case of
codimension two was due to Gaeta, and this generalization thus bears his name.

\begin{theorem}[Generalized Gaeta theorem] \label{thm-gaeta}
Every standard determinantal ideal is glicci.
\end{theorem}

\begin{proof}
The proof is essentially an algorithm describing how the required
links can be achieved. We outline the steps of this algorithm but
refer  to \cite{KMMNP}, Theorem 3.6 for the complete proof. The
interested reader is invited to run the algorithm with a concrete
example.
\smallskip

Let $I \subset R$ be a standard determinantal
   ideal of codimension $c+1$. Thus, there is a homogeneous $t \times (t+c)$
   matrix $A$  with entries in $R$ such that $I = I(A)$. If $t =1$ then $I$
   is a complete intersection and there is
   nothing to prove. Let $t > 1$. Then our assertion follows by induction
   on $t$ if
   we have shown that $I$ is evenly G-linked to a standard determinantal
   scheme $I'$  generated by the maximal minors of
   a  $(t-1) \times (t+c-1)$ matrix $A'$. Actually, we will see
   that
   $A'$ can be chosen as the matrix which we get after deleting an appropriate
   row and column of the matrix $A$ and that then $I$ and $I'$ are directly
   G-linked in two steps. In order to do that we proceed in several
   steps.
   \smallskip

{\it Step} I: Let $B$ be the matrix consisting of the first
$t+c-1$ columns of $A$.  Then the ideal $\fa := I(B)$ has
codimension $c$, i.e.\ it is a standard determinantal ideal.

   Possibly after elementary row operations we may assume that the
maximal minors of the matrix $A'$ consisting of the first $t-1$
rows of $B$ generate an ideal of (maximal) codimension $c+1$. Denote this
standard determinantal ideal by $I' := I(A')$.
\smallskip

{\it Step} II: Possibly after elementary column operations we may assume
that the maximal minors of the matrix $A_1$  consisting of the first
$t-1$ columns of $A$ generate an ideal of (maximal) codimension two.
Put $J = I(A_1)$. Let $d$ be the determinant of the matrix which consists
of the first $(t-1)$ and the last column of $A$. Then one can show
that
\begin{itemize}
\item[(i)] $\fa : d = \fa$.
\item[(ii)] $I = (\fa +d R) : J$.
\item[(iii)] $\fa + d J^{c-1}$ is a Gorenstein ideal of
codimension $c+1$.
\item[(iv)] $\deg \fa + d J^{c-1} = \deg d \cdot \deg \fa  + \deg
(\fa + J^{c-1})$.
\end{itemize}
\smallskip

{\it Step} III: Consider for  $i = 0,\ldots,c$ the ideals $I_S +
J^{i}$. These are Cohen-Macaulay ideals of degree
$$
\deg (\fa + J^{i}) = i \cdot [\deg d \cdot \deg \fa  - \deg I].
$$
The proof involves in particular a deformation argument.
\smallskip

{\it Step} IV: Comparing degrees it is now not to difficult to
check that
$$
(\fa + d J^{c-1}) : I = \fa + J^c.
$$
\smallskip

{\it Step} V: Let  $d'$ be the determinant of the matrix which consists
of the first $(t-1)$  columns of $A'$. Then, similarly as above,
$\fa + d' J^{c-1}$ is a Gorenstein ideal of codimension $c+1$ and
$$
(\fa + d' J^{c-1}) : I' = \fa + J^c.
$$
\smallskip

{\it Step} VI: Step V says that the ideal $I'$ is directly
G-linked to $\fa + J^c$ while Step IV gives that $I$ is directly
G-linked to $\fa + J^c$. Hence the proof is complete.
\end{proof}

\begin{example} \label{ex-rat-n-curve}
Let $C \subset \PP^n$ denote a rational normal curve. It is
well-known that after a change of coordinates we may assume that
the homogeneous ideal of $C$ is generated by the maximal minors of
the matrix
$$
\left ( \begin{array}{lcr}
x_0 & \ldots & x_{n-1} \\
x_1 & \ldots & x_n
\end{array} \right ).
$$
Hence, $C$ is standard determinantal and therefore glicci by
Gaeta's theorem.

On the other hand the curve $C$ has a linear free resolution,
i.e.\ its minimal free resolution has the shape
$$
0 \to R^{\beta_{n-1}} (-n) \to \ldots \to R^{\beta_1} (-2) \to I_C
\to 0.
$$
Hence, \cite{HU-Annals}, Corollary 5.13 implies for $n \geq 4$
that the curve $C$ is not
licci, i.e.\ not in the CI-liaison class of a complete
intersection.
\end{example}

Let us look back to the questions posed at the beginning of this
section. The previous example shows that the answer to Question 3
is `no', i.e.\ Rao's correspondence is not injective for
CI-liaison in codimension $\geq 3$. Gaeta's theorem indicates that
the situation might be different for G-liaison. In fact, there is
more evidence that Question 1 could have an affirmative answer. To
this end we consider certain monomial ideals.

\begin{definition} \label{def-stable}
A monomial ideal $J \subset R$ is said to be  {\em stable} if
\[
m = x_0^{a_0}\cdots x_n^{a_n} \in J \mbox{ and }  a_i >0  \hbox{ imply }
\frac{x_j}{x_i} \cdot m \in J
\]
for all $1 \leq j < i \leq n$.
\end{definition}

\begin{theorem}\label{thm-stable}
Suppose that the ground field $K$ is infinite. Then every
Cohen-Macaulay Borel-fixed monomial ideal is glicci.
\end{theorem}

\begin{proof}
The main tools are basic double links and liftings of monomial
ideals (cf.\ \cite{MN2}). We only outline the main steps of the proof and
refer for details to \cite{MN4}. Moreover, we assume for
simplicity that $K$ has characteristic zero.
\smallskip

{\it Step I}: Let $J \subset R$ be a Cohen-Macaulay stable monomial ideal of
codimension $c+1$. Denote by
$$
\alpha := \min \{ t \in \ZZ \s [I]_t \neq 0 \}
$$
its initial degree. Then there are uniquely determined Artinian
stable ideals $I_0,\ldots,I_{\alpha} \subset T :=
K[x_1,\ldots,x_c]$ such that
$$
I_0 \subset I_1 \subset \ldots \subset I_{\alpha} = T
$$
and
\[
\begin{array}{rcl}
J & = & I_0 R + x_0 I_1 R + x_0^2 I_2 R+ \dots + x_0^\alpha I_\alpha
R\\[1ex]

      & = & I_0 R + x_0 I'
\end{array}
\]
where $I' = I_1 R + x_0 I_2 R + \dots + x_0^{\alpha -1} I_\alpha R$.

It follows that
\begin{itemize}
\item[(i)] $I_0 R$ is a Cohen-Macaulay ideal of codimension $c$.
\item[(ii)] $I_0 R \subset I'$ because $I_0 R \subset I_1 R
\subset I'$.
\item[(iii)] $I'$ is a Cohen-Macaulay ideal of codimension $c+1$.
\end{itemize}
\smallskip

{\it Step II}: Now we want to lift monomial ideals in $T$ to
reduced ideals in $S := T[x_0]$.

Consider the lifting map $\lambda: \{\mbox{monomials in } T\} \to
\{\mbox{monomials in } S\}$ given by
$$
\prod_{j=1}^c x_j^{a_j} \mapsto \prod_{j=1}^c \left (
\prod_{i=0}^{a_j -1} (x_j + i x_0) \right ).
$$
(Here the assumption on the characteristic is used. In
general, one just has to choose sufficiently general linear forms in
order to replace monomials by products of linear forms as above.)
For example, we get $\lambda(x_1^3 x_2^2) = x_1 (x_1 + x_0) (x_1 +
2 x_0) x_2 (x_2 + x_0)$.

The properties of the lifting map ensure that $\lambda(I_0)$ is a
reduced ideal defining a set of points in $\PP^c$. Therefore this
set has the property $G_1$.
\smallskip

{\it Step III}: Using the stability of $J$ one can show that
$\lambda(I_0)R \subset I'$ and
$$
J = \lambda(I_0) R + x_0 I'.
$$
Thus, $J$ is a basic double link of $I'$ and Proposition \ref{bdl is
linked} shows that $I' \in \cL_J$. But the initial degree of $I'$
is $\alpha - 1$. Repeating this argument successively we see that
$I_0 R + x_0 R \in \cL_J$. Hence it is sufficient to show that
$I_0 R + x_0 R$ is glicci. But this follows because $I_0 R$ is glicci
by induction on the codimension. The claim is clearly true for
ideals of codimension one.
\end{proof}

Theorem \ref{thm-stable} is of a more general nature than it is
apparent from its formulation.

\begin{remark} \label{rem-stable}
Let $V \subset \PP^n$ be an \aCM\ subscheme.
  It is
well-known its generic initial ideal $gin (I_V)$ is a stable ideal
and defines an \aCM\ subscheme  which is a deformation of the original
scheme $V$.  Indeed, the fact that $gin (I_V)$ is stable is due to Galligo
\cite{galligo};  that it gives a flat deformation is due to Bayer
\cite{bayer}; that it is again Cohen-Macaulay follows from a result of Bayer
and Stillman (cf.\ \cite{eisenbud}, Theorem 15.13).  Thus our result says that
every \aCM\ subscheme admits a flat deformation which is glicci.  In other
words, we have found an affirmative answer to Question 1 ``up to
flat deformation.''
\end{remark}

In view of Remark \ref{rem-stable}, we consider Theorem
\ref{thm-stable} as the strongest evidence that Question~1 might
have an affirmative answer. However, there is also other evidence.

\begin{remark} \label{rem-acm-glicci}
The results about linear systems (Theorem \ref{linear equiv
theorem}) can be used to show that many \aCM\ subschemes are glicci.
This becomes particularly effective for divisors on \aCM\
subschemes with known  Picard group. Some typical results of this
approach are
\begin{itemize}
\item[(i)] All \aCM\ curves on a general smooth rational \aCM\
surface in $\PP^4$ are glicci (\cite{KMMNP}, Corollary 8.9).
\item[(ii)] Let $S \subset \PP^4$ be a general \aCM\ surface such
that all the entries of its Hilbert-Burch matrix have positive
degree. Then all \aCM\ curves on $S$ are glicci (\cite{CM2}).
\item[(iii)] Effective \aCM\ divisors on a smooth rational normal
scroll are glicci (\cite{Cas-thesis}).
\item[(iv)] Every general set of points in $\PP^3$ on a
nonsingular quadric surface is glicci (\cite{Hart-co-3}). More
generally, every general set of points on a smooth rational
surface scroll is glicci (\cite{Cas-thesis}, Theorem 3.4.2).
\end{itemize}
\end{remark}

One of the few sufficient conditions for linkage in higher codimension was
mentioned in Remark \ref{Walter result}, and now we sketch the proof.

\begin{proposition}
Any two complete intersections of the same codimension are CI-linked.
\end{proposition}

\begin{proof} (Sketch of proof from \cite{schwartau})

The proof rests on the following observation:  If $I_{X_1} =
(F_1,\dots,F_{c-1},F)$ and $I_{X_2} = (F_1,\dots,F_{c-1},G)$ are two complete
intersections of codimension $c$ then they are directly linked by the complete
intersection $I_X = (F_1,\dots,F_{c-1},FG)$.  Then the proof follows by
changing one entry at a time.
\end{proof}

  From this it would follow that if one could show that every \aG\ scheme is
glicci, then all \aG\ schemes are in the same G-liaison class.  However,
this is
not known. It is true if the codimension is at most three. Then an \aG\
subscheme
is even licci (\cite{Wata-liai}).

Moreover, Hartshorne \cite{Hart-co-3} has proposed interesting
examples. He suspects that a set of $20$ general points in $\PP^3$
as well as the general curve in the irreducible component of the
Hilbert scheme of curves in $\PP^4$ of degree $20$ and genus $26$
containing standard determinantal curves
is not glicci.

There are also some results for non-\aCM\ subschemes indicating
that even the Main question might have an affirmative answer:

\begin{itemize}

\item  Hartshorne \cite{Hart-co-3} and Lesperance \cite{lesperance1}
independently showed that any two sets of two skew lines in $\proj{4}$ are
G-linked.  (See also Conjecture \ref{2 skew lines in p4}.)   Hartshorne also
obtained partial results on other curves with Rao module $k$.

\item Lesperance \cite{lesperance1} showed that curves in $\proj{4}$ consisting
of unions of two plane curves are (at least ``usually'') linked if and only if
they have the same Rao module.

\item Lesperance \cite{lesperance2} showed that if $C$ and $C'$ are degenerate
\aBM\ curves in $\proj{4}$ (not necessarily in the same hyperplane) then
$C$ and
$C'$ are evenly G-linked if and only if they have isomorphic Rao modules up to
shift.

\item Casanellas and Mir\'o-Roig \cite{CM1}, \cite{CM2} showed the same for
many
subschemes of small degree (not necessarily curves), especially unions of
linear
varieties; their idea was to view them as divisors on a suitable rational
normal
scroll.

\item Nagel, Notari and Spreafico \cite{NNS2} proved for double lines in
$\proj{n}$ and for some other non-reduced curves on lines, that they are evenly
linked if and only if they have isomorphic Rao modules up to shift.
\end{itemize}

The proof of the last result differs from the others by not using the
result about the G-liaison classes of divisors on \aCM\ subschemes with the
property $G_1$. Indeed, the non-reduced curves that are considered are not even
divisors on a generically Gorenstein surface.

Note also that the Hartshorne-Rao modules of the curves considered in the
first results mentioned above are rather simple while the curves studied in
\cite{NNS1}, \cite{NNS2} can have a  rather complicated Hartshorne-Rao
module.


\section{The Structure of an Even Liaison Class} \label{str-even-liaison}

We have seen a rather complete description of when two subschemes are
linked in codimension two.  The main result is Theorem \ref{thm-rao-inj},
and it is one of the main results of liaison theory.  We have discussed
to some extent the possibility of extending this result to higher
codimension (e.g.\ Theorem \ref{thm-gaeta}), and we will continue to
discuss it below.  As we saw, it is more natural to consider {\em even}
liaison.

Another natural question is whether the even liaison classes possess a
common structure of any sort.  We will see that in codimension two
there is a nice answer.  Again, one can try to extend it to higher
codimension, and we will also discuss the evidence for and against this
idea.  The following remark sets up the background.

\begin{remark} \label{setup for LR}
Let $V \subset \proj{n}$ be an equidimensional closed subscheme of
codimension $c$.  Let $M_i = H^i_*({\mathcal I}_V)$ for
$1 \leq i \leq \dim V = n-c.$  Let ${\mathcal L}_V$ be the even liaison
class of $V$.  Note that

\begin{itemize}

\item The vector of graded modules $M_\bullet = (M_1,\dots,M_{n-c})$ is an
invariant of ${\mathcal L}_V$, up to
shift (Lemma \ref{lem-coho-lCM-cl}).

\item There is a minimal shift of this vector that can occur among subschemes
of $\proj{n}$ (Proposition \ref{prop-coho-min-shift}; see Definition
\ref{def-coho-min-shift} for the definition of minimal shift).

\item Hence there is a minimal shift of this vector among
elements of ${\mathcal L}_V$ (which is not necessarily the same as the
minimal shift among all subschemes in $\proj{n}$ with vector ${\mathcal
L}_V$, except for curves in $\proj{3}$).

\item Although leftward shifts of $M_\bullet$ may not exist, any rightward
shift of $M_\bullet$ does exist thanks to Basic Double Linkage (Lemma
\ref{lem-basic-d-link}, Remark \ref{basic double ci linkage}).
\end{itemize}
\end{remark}

\begin{definition} \label{def of L0}
If an element $W \in {\mathcal L}_V$ has cohomology which achieves the minimal
shift, among elements of ${\mathcal L}_V$, guaranteed by Proposition
\ref{prop-coho-min-shift}, we say that $W$ is a {\em minimal element of its
even liaison class}, or that {\em $W$ is in the minimal shift of ${\mathcal
L}_V$}.  We write $W \in {\mathcal L}_V^0$.
\end{definition}

\begin{example}
Let $Z_1$ be the disjoint union in $\proj{3}$ of a line, $\lambda$, and a
conic, $Y$.  We have the exact sequence
\[
\begin{array}{cccccccccccccccccccc}
0 & \rightarrow & I_{Z_1} & \rightarrow & I_\lambda \oplus I_Y &
\longrightarrow & \ \ R & \rightarrow & H^1_*({\mathcal I}_{Z_1}) &
\rightarrow & 0 \\
&&&&\hfill \searrow && \nearrow \hfill \\
&&&&&I_\lambda + I_Y \\
&&&& \hfill \nearrow && \searrow \hfill \\
&&&& 0 && \hfill \ \ \ \ 0
\end{array}
\]
Since $I_\lambda + I_Y$ contains three independent linear forms, we
conclude that $H^1_*({\mathcal I}_{Z_1}) \cong K[x]/(x^2)$ for some
linear form $x$.  Therefore the module is one-dimensional in each of
degrees 0 and 1, and zero everywhere else.  However, notice that the
module structure is not trivial: multiplication from the degree 0
component to the degree 1 component by the linear form $x$ is not zero.

Note further that this curve is in the minimal shift of its even liaison class,
thanks to the bound  (\ref{incr-in-neg-deg}) which says that in negative
degree the dimensions have to be strictly increasing.

Now consider a Buchsbaum curve $Z_2$, obtained via Liaison Addition as
in Example \ref{construction of aBM curves}, with deficiency module which is
1-dimensional in each of two consecutive degrees.  The smallest such curve
that can be so constructed is obtained by choosing $C_1$ and $C_2$ in Example
\ref{construction of aBM curves} to each be a pair of skew lines,  and
$\deg F_1 = 2$, $\deg F_2 = 3$.  Then the first non-zero component of
$H^1_*({\mathcal I}_{Z_1})$ occurs in degree 2 which, thanks to Proposition
\ref{buchs leftmost shift}, is the minimal shift.

Note that the {\em structure} of these two modules, $H^1_*({\mathcal
I}_{Z_1})$ and $H^1_*({\mathcal I}_{Z_2})$, is different (the latter is
annihilated by all linear forms), even though dimensionally they are the
same.  Hence they are not in the same (even) liaison class.

For an example of surfaces where even the modules are isomorphic but the
liaison classes are different, see Example \ref{ex-rao-corr}.
\end{example}

\begin{remark} \label{rem-shifts}
Let $V \subset \PP^n$ be temporarily an equidimensional scheme of
codimension $c \geq 2$. Then it is clear how to adapt the above definition
of $\cL^0_V$. Note, that we have already defined the (cohomological)
minimal shift in Definition \ref{def-coho-min-shift}. Strictly speaking we
should distinguish even a third notion of minimal shift suggested by Rao's
correspondence. This provides the
following list:
\begin{itemize}
\item[(i)] The (cohomological) minimal shift
of $V$ is the integer
$$
c(V) := \min \left \{ t \in \ZZ \s \begin{array}{c}\mbox{There is a
       subscheme}\  W \subset \PP^n \ \mbox{ of codimension}\ c\ \mbox{with} \\
  H^i_*(\cI_W) \cong H^i_*(\cI_V)(-t)
\end{array}  \right \}.
$$
\item[(ii)] The minimal Rao shift of $V$ is the integer
$$
r(V) := \min \left \{ t \in \ZZ \s \begin{array}{c}\mbox{There is a
       subscheme}\  W \subset \PP^n \ \mbox{ of codimension}\ c\ \mbox{with} \\
  \ffi(I_V) \oplus F \cong \ffi(I_W)(-t) \oplus G\ \mbox{for free
   $R$-modules}\ F, G
\end{array}  \right \}.
$$
\item[(iii)] The minimal shift of the even G-liaison class $\cL_V$ is the
   integer
$$
l(V) := \min \left \{ t \in \ZZ \s \begin{array}{c}\mbox{There is a
       subscheme}\  W \in \cL_V\ \mbox{with} \\
  \ffi(I_V) \oplus F \cong \ffi(I_W)(-t) \oplus G\ \mbox{for free
   $R$-modules}\ F, G
\end{array}  \right \}.
$$
\end{itemize}
According to Remark \ref{rem-E-N-type} and Rao's correspondence we have the
following inequalities
$$
c(V) \leq r(V) \leq l(V).
$$
Moreover, if $V$ is a curve then $c(V) = r(V)$, but if the dimension of $V$
is at least $2$ we can have $c(V) < r(V)$.

If the codimension of $V$ is two we get $r(V) = l(V)$ due to Theorem
\ref{thm-rao-inj}. It would be interesting to know if this equality is also
true in codimension $c \geq 3$. This would follow from an affirmative
answer to the Main question \ref{mque-inj}, but it is conceivable that the
Main question has a negative answer and $r(V) = l(V)$ is still always true.

In \cite{N-gorliaison}, Proposition 5.1  a lower bound for $r(V)$ is given
which cannot be improved in general. It would be interesting to have a
priori estimates for $c(V)$ and $l(V)$ as well.
\end{remark}

We now describe a structure of an even liaison class, generally called the
{\em Lazarsfeld-Rao property}.  As remarked above, this property is only known
to hold in codimension two, so we now make this assumption.  Later we will
discuss the possibility of extending it.

Let $\mathcal L$ be an even liaison class of codimension two subschemes of
$\proj{n}$.  For simplicity we will assume that the elements of $\mathcal L$
are locally Cohen-Macaulay, and of course they must be equidimensional.  (The
locally Cohen-Macaulay assumption was removed by Nagel \cite{N-gorliaison} and
by Nollet \cite{nollet}.)

As we have seen (e.g.\ Theorem \ref{thm-gaeta}), the \acm\ codimension two
subschemes {\em form} an even liaison class.  (In this case any two schemes
are both evenly and oddly linked.)  We thus assume that the elements of
$\mathcal L$ are {\em not} \acm, so $M_\bullet$ is not zero (i.e.\ at least one
of the modules, not necessarily all, is non-zero).  Then it follows from Remark
\ref{setup for LR} that we can partition $\mathcal L$ according to the shift
of $M_\bullet$:
\[
{\mathcal L} = {\mathcal L}^0 \cup {\mathcal L}^1 \cup {\mathcal L}^2 \cup
\dots \cup {\mathcal L}^h \cup \dots.
\]
Here, ${\mathcal L}^0$ was defined in Definition \ref{def of L0} and consists
of the minimal elements.  Then ${\mathcal L}^h$ consists of those elements of
$\mathcal L$ whose deficiency modules are shifted $h$ degrees to the right of
the minimal shift.

In Remark \ref{basic double CI-link} we saw the notion of Basic Double
CI-Linkage and in particular we gave the version for codimension two:
Let $V_1$ be a codimension two subscheme of $\proj{n}$ and choose $F_2 \in
I_{V_1}$ of degree $d_2$ and $F_1 \in R$ of degree $d_1$ such that $(F_1,F_2)$
forms a regular sequence (i.e.\ a complete intersection).  Then $F_1 \cdot
I_{V_1} + (F_2)$ is the saturated ideal of a scheme $Z$ which is CI-linked to
$V_1$ in two steps.  Furthermore,
\[
\phantom{\hbox{\hskip .8in for }
i=1,\dots,n-2.}
H_*^i({\mathcal I}_Z) \cong H^i_*({\mathcal I}_{V_1})(-d_1) \hbox{\hskip .8in
for } i=1,\dots,n-2.
\]
As sets, $Z = V_1 \cup V$ where $V$ is the complete intersection defined by
$(F_1,F_2)$.  Note that if $V_1 \in {\mathcal L}^h$ then $Z \in {\mathcal
L}^{h+d_1}$.  A concrete description of the two links can be given as follows
(first noted in \cite{LR}):  Let $A \in I_{V_1}$ be any homogeneous polynomial
having no component in common with $F_2$.  Then link $V_1$ to some
intermediate scheme $Y$ using the complete intersection $(A,F_2)$, and link
$Y$ to $Z$ using the complete intersection $(AF_1,F_2)$.

One can also check, using various methods, that the $E$-type resolutions of
$V_1$ and $Z$ are related as follows.  If ${\mathcal I}_{V_1}$ has an $E$-type
resolution
\[
0 \rightarrow {\mathcal E} \rightarrow \bigoplus_{i=1}^m {\mathcal
O}_{\proj{n}} (-a_i) \rightarrow {\mathcal I}_{V_1} \rightarrow 0,
\]
where $H^1_*({\mathcal E}) = 0$, then ${\mathcal I}_Z$ has an $E$-type
resolution
\begin{equation} \label{E basic double link}
0 \rightarrow {\mathcal E}(-d_1) \oplus {\mathcal O}_{\proj{n}} (-d_1-d_2)
\rightarrow \bigoplus_{i=1}^m {\mathcal O}_{\proj{n}} (-d_1-a_i) \oplus
{\mathcal O}_{\proj{n}} (-d_2) \rightarrow {\mathcal I}_{Z} \rightarrow 0,
\end{equation}
Note that the stable equivalence of $\mathcal E$ and ${\mathcal E}(-d_1)
\oplus {\mathcal O}_{\proj{n}} (-d_1-d_2)$ is obvious.

The Lazarsfeld-Rao property says, basically, that in an even liaison class,
all the minimal elements look alike and that the entire class can be built up
from an arbitrary minimal element using Basic Double Linkage and deformation.
More precisely, we have the following statement.

\begin{theorem} \label{LR property} {\bf (Lazarsfeld-Rao property)}
Let $\mathcal L$ be an even liaison class of codimension two subschemes of
$\proj{n}$.

\begin{itemize}

\item[(a)] If $V_1,V_2 \in {\mathcal L}^0$ then there is a flat deformation
from one
to the other through subschemes all in ${\mathcal L}^0$.

\item[(b)] If $V_0 \in {\mathcal L}^0$ and $V \in {\mathcal L}^h$ ($h \geq 1$)
then there is a sequence of subschemes $V_0,V_1,\dots,V_t$ such that for all
$i$, $1 \leq i \leq t$, $V_i$ is a basic double link of $V_{i-1}$, and $V$ is a
deformation of $V_t$ through subschemes all in ${\mathcal L}^h$.
\end{itemize}
\end{theorem}

We stress that the deformations mentioned in Theorem \ref{LR property} are
carried out entirely within the even liaison class $\mathcal L$.  They preserve
cohomology, not only dimensionally but even structurally.

Theorem \ref{LR property} was first proved for codimension two locally
Cohen-Macaulay subschemes of $\proj{n}$ in \cite{BBM}.  At approximately the
same time, it was proved (as part of a much broader theory) for curves in
$\proj{3}$ in \cite{MP}.  It was proved for codimension two subschemes of a
smooth
arithmetically Gorenstein subscheme in \cite{BM6}. Finally, in codimension two
it was later extended to arbitrary unmixed ideals in \cite{N-gorliaison} and
\cite{nollet}.  We now give the general idea of the proof of \cite{BBM}, and
refer the reader to that paper for the details, as well as to \cite{MP},
\cite{N-gorliaison} and \cite{nollet}.

\begin{proof} (Sketch)
There are three basic components of the proof.

\begin{itemize}
\item[1.] (Bolondi, \cite{bolondi}) If $V_1,V_2 \in {\mathcal L}^h$ (in
particular they have the same deficiency modules) and if they have the same
Hilbert function then the desired deformation can be found.  So it is
reduced to
a question of Hilbert functions.

\item[2.] If $V_1,V_2 \in {\mathcal L}^h$ and if they do {\em not} have the
same
Hilbert function then by studying locally free $N$-type resolutions one can
show
that there is a ``smaller" $V'$ in the even liaison class  (i.e.\ $V' \in
{\mathcal L}^{h'}$ for some $h' < h$).  Combined with the first part, this
proves
that the minimal elements all lie in the same flat family.

\item[3.] Given $V_0 \in {\mathcal L}^0$ and $V \in {\mathcal L}^h$, by
studying
(\ref{E basic double link}) and knowing that $V_0$ and $V$ are linked in an
even
number of steps, it is possible to ``predict'' what basic double links are
needed to start with $V_0$ and arrive at a scheme $V_t$ with $E$-type
resolution
which agrees (except for the maps) with that of $V$, up to trivially adding
free
summands to both modules in the resolution.  This means that $V$ and $V_t$ have
the same Hilbert function and deficiency modules, so we again apply the first
part.
\end{itemize}
\end{proof}

\begin{remark}
(i) From the name ``Lazarsfeld-Rao property'' one would naturally
expect that the
paper \cite{LR} of Lazarsfeld and Rao was important in the development of the
above theorem.  In fact, it really inspired it (although many people doubted
that something so general would hold).  We can state the main result of
\cite{LR} in the following way.  For a curve $C \subset \proj{3}$, let
\[
e(C) := \max \{ t | h^2({\mathcal I}_C(t)) \neq 0 \} =
\max\{ t | h^1({\mathcal O}_C (t)) \neq 0 \} = \max\{ t | h^0({\omega_C}(-t))
\neq 0 \}.
\]
   Then
\begin{itemize}
\item[a.] If $C$ lies on no surface of degree $e(C) +3$ then ${\mathcal L}_C$
has the Lazarsfeld-Rao property and $C \in {\mathcal L}_C^0$.

\item[b.] If $C$ lies on no surface of degree $e(C)+4$ then furthermore $C$ is
the {\em only} element of~${\mathcal L}_C^0$.
\end{itemize}

For example, suppose that $C \subset \proj{3}$ is a set of $\geq 2$ skew
lines.
Then $e(C) =-2 $.  Thus since $C$ cannot lie on a surface of degree 1,
part a.\
gives that $C \in {\mathcal L}_C^0$.  If $C$ furthermore does not lie on a
quadric surface then $C$ is the only minimal element of its even liaison
class.

Similarly, one can apply it to rational curves, where $e(C) = -1$, and get
analogous statements: a rational curve lying on a quadric surface is not
minimal (it is linked to a set of skew lines), one lying on a cubic surface is
minimal but not unique (it moves in a linear system) and one not lying on a
cubic surface is the unique minimal curve.  A different, more geometric
approach
to the minimality of skew lines and rational curves (not using \cite{LR}), and
other related questions, can be found in \cite{Mi-invar}.

\medskip

(ii) Let us recall the concept of elementary CI-biliaison (in the case of
curves). Let $C \subset \PP^n$ be a curve which is an effective
divisor on a complete intersection surface $S \subset \PP^n$. Let
$F_1$
be a hypersurface meeting $S$ transversally such that $C \subset S
\cap F_1$. Let $C'$ be the curve linked to $C$ by $S \cap F_1$. Choose
a hypersurface $F_2$ such that it meets $S$ transversally and $C'
\subset S \cap F_2$. Denote by $C''$ the curve linked to $C'$ by $S \cap F_2$.
Then it is said that $C''$ is obtained from $C$ by an {\it elementary
   CI-biliaison} on $S$. It is called {\it ascending} if $\deg F_2 - \deg F_1
 \geq 
0$, otherwise {\it decsending}. As already indicated in Remark
\ref{rem-Hartsh} $C''$ is obtained from $C$ by an  elementary
   CI-biliaison on $S$ if and only if $C'' \sim C + h H$. 
Observe that elementary CI-biliaison is a
generalization of basic double CI-linkage (cf.\ Remark \ref{basic
   double CI-link}). Recently, in \cite{strano-biliai} R.\ Strano has
obtained the following
variant of the Lazarsfeld-Rao property: Let $C \subset \PP^3$ be a
curve which is not \aCM. Then $C$ can be obtained from a minimal curve in
its even liaison class by finitely many ascending elementary
CI-biliaisons. Thus, using the more general elementary biliaison instead
of basic double links we can avoid the possible final deformation
which is allowed in Theorem \ref{LR property}.
\end{remark}

\begin{remark} \label{all possible d,g}
If one knows the Hilbert function of a curve $C$ in $\proj{3}$ (or of a
codimension
two subscheme in general) then one can write the Hilbert function of all
possible
basic double links from $C$.  Hence the Lazarsfeld-Rao property can be used to
give a complete list of all possible $(d,g) = \hbox{(degree,genus)}$
combinations
that occur in an even liaison class ($g$ is the {\em arithmetic}
genus), if one
only knows it for a minimal element.

For example consider curves in $\proj{3}$ that are \aBM\ but not
arithmetically Cohen-Macaulay.  In Example
\ref{construction of aBM curves} we saw one way to construct them, such
that the
result has its leftmost component in degree $2N-2$ (where $N = \dim_k
H^1_*({\mathcal I}_C)$), which according to Proposition \ref{buchs leftmost
shift} makes it a minimal element of its even liaison class.  Hence its degree,
genus, and even its Hilbert function, are uniquely determined, thanks to the
Lazarsfeld-Rao property.  The following, from \cite{migbook}, is a complete
list
of the possible $(d,g)$ that can occur for \aBM\ curves in $\proj{3}$ when $d
\leq 10$.  It includes two curves for which $N = 2$: one with $(d,g) = (8,5)$
and $H^1_*({\mathcal I}_C)$ concentrated in degree 2, and one with $(d,g) =
(10,
10)$ and $\dim_k H^1_*({\mathcal I}_C)_2 = \dim_k H^1_*({\mathcal I}_C)_3 =
1$.
The rest have $\dim_k H^1_*({\mathcal I}_C) = 1$.

\medskip

\begin{center}
\begin{tabular}{|c|c|c|c|c|c|c|c|c|c|ccccccccccccccccccc}
\hline
degree &
2  & 3 & 4 & 5 & 6 & 7 & 8 & 9 & 10 \\
\hline
genus & $-1$ & \hbox{d.n.e.} & 0 & 1 & 3 & 4,6 & 5,6,8 10 & 8, 9, 15 & 10,
11, 13, 15, 21 \\
\hline
\end{tabular}
\end{center}
Note that there is no such curve of degree 3.
\end{remark}

\bigskip

As was the case with the necessary and sufficient conditions for G-linkage, the
biggest open problem is to find a way to extend these results to higher
codimension.  One intermediate situation was studied in \cite{BM6} (cf.\
also \cite{N-gorliaison}), where
liaison was studied not in projective space but rather on a smooth \aG\
subvariety $X$ of projective space.  It was shown that codimension two liaison
here  behaves almost identically to that in $\proj{n}$, even though of
course the
objects being linked have codimension greater than two in $\proj{n}$. These
results have been further generalized in \cite{N-gorliaison} to codimension
two subschemes of an arbitrary integral \aG\ subscheme.

One interesting difference concerns \acm\ subvarieties.  Here we mean \acm\ in
projective space (i.e.\ the deficiency modules vanish), but such a subvariety
need not have a finite resolution over the Gorenstein coordinate ring $R/I_X$.
It was shown that the notion of minimality still makes sense, viewed not in
terms of the shift of the modules (which are zero) but rather in terms of
$N$-type resolutions.  Then it was shown that the Lazarsfeld-Rao property holds
in such a situation on $X$.

Note that the linkage on $X$ is by complete intersections on $X$, which however
are only \aG\ as subschemes of $\proj{n}$. But if we turn to Gorenstein liaison
in $\proj{n}$ with no such restriction, the situation becomes much
less optimistic.

First, we can see right away that there is no hope for a statement which is
identical to that for codimension two.  The following example was taken from
\cite{migbook}.  Consider the non-degenerate curve in $\proj{4}$ in the
following configuration:


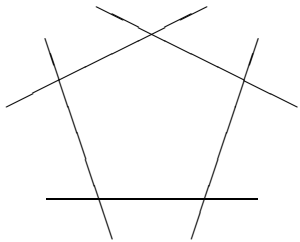
\begin{figure}[ht]\label{deg-5-gor-curve}
  \begin{picture}(370,100)
  \put (145,15){\line (1,0){80}}
\put (170,0){\line (-1,3){25.3}}
\put (200,0){\line (1,3){25.3}}
\put (130,50){\line (2,1){75.9}}
\put (240,50){\line (-2,1){75.9}}
  \end{picture}
\end{figure}

\noindent This curve is \aG.  As such, it links two skew lines to a curve of
degree 3 consisting of the disjoint union of a line and two lines meeting in a
point.  One checks that both of these curves have Rao module which is one
dimensional, occurring in degree 0.  Since this is the minimal shift, it is
clear that the elements of ${\mathcal L}^0$ do not all have the same degree,
hence are not in a flat family.

So if there is a nice structure for an even liaison class under Gorenstein
liaison, what should the statement be?  The next natural guess, due to
Hartshorne \cite{Hart-co-3} is that perhaps the elements of ${\mathcal L}^0$
satisfy the property that while there may be curves of different degrees, those
curves of the same degree at least lie in a flat family.  He showed this for
the liaison class of two skew lines.  However, it was shown to be false in
general by Lesperance \cite{lesperance1}, who gave an example of two sets of
curves ``usually'' in the same even liaison class which are in the minimal
shift
and have the same degree and even arithmetic genus, but which do not lie in the
same flat family.  His example was extended somewhat by Casanellas
\cite{Cas-thesis}, who looked at the same kind of curves but in $\proj{5}$.
(Lesperance was not able to show that all of his curves are in the same even
liaison class, even though they do have the same Rao module.  Casanellas showed
that this obstacle disappears in $\proj{5}$.)

So at the moment no one has a good idea of how to find an analog to the
Lazarsfeld-Rao property for Gorenstein liaison of subschemes of $\proj{n}$ of
codimension $\geq 3$.  A first problem seems  to be to find a good concept of a
minimal element of an even G-liaison class. In the even liaison class $\cL$
of a non-\aCM\ subscheme of codimension two, the minimal elements are the
elements
of smallest degree in $\cL$  and all these elements have the same Hilbert
function. In particular, a non-\aCM\ curve of degree two in $\PP^3$ must be
minimal in its even liaison class. In higher codimension the situation is
very  different. It is still true that two curves of degree two in $\PP^n$
are in the same even liaison class if and only if their Hartshorne-Rao
modules are isomorphic according to  \cite{NNS1}, but such curves can have
different
genera.

A naive idea would be to define the minimal elements in an even
G-liaison class as the ones achieving the minimal shift and having
minimal Hilbert polynomial. Consider the curves of degree two in
$\PP^n$ whose Hartshorne-Rao module  is isomorphic to  the ground
field $K$.  Such a curve can have every arithmetic genus $g$
satisfying $- \frac{n-1}{2} \leq g \leq -1$, but it is
non-degenerate if and only if $- \frac{n-1}{2} \leq g \leq 2-n$
(\cite{NNS1}). Thus, for $n \geq 4$ minimal curves in the sense
just discussed were  degenerate.

It should be remarked that the authors wonder if the
Lazarsfeld-Rao property, even as it is stated in codimension two, might
hold for
CI-liaison in higher codimension. There are some encouraging result in
\cite{HU-Annals}.


\section{Remarks on the different liaison concepts} \label{compare section}

We have already seen that for subschemes whose codimension is at least
three, G-linkage and CI-linkage generate very different equivalence
classes. In this section we want to discuss these differences a
bit more systematically. Finally, we compare briefly the
equivalence classes generated by (algebraic) CI-linkage and
geometric CI-linkage.

As we have mentioned in Section \ref{sec-necess},  Rao's correspondence
gives  the only known method for distinguishing  between G-liaison
classes. The situation is different for CI-liaison. There are
various invariants, numerical (\cite{Huneke-liai}) as well as structural
(\cite{Buch-Ul}, \cite{HU-Annals},  \cite{KMMNP}),  which allow one to
distinguish between
CI-liaison classes of \aCM\ subschemes. In order to give the
flavour of such invariants, we state a particularly clean result
which has been shown in \cite{Buch-Ul} by algebraic means, whereas
a more geometric proof has been given in \cite{KMMNP}, Proposition
6.8.

\begin{theorem} \label{thm-ci-inv}
Suppose  $V, W \subset \PP^n$, $n \geq 4$, are \aCM\ subschemes of
codimension $3$. If $V$ and $W$ belong to the same CI-liaison
class then there are isomorphisms of graded $R$-modules
$$
\HH^i(K_V \otimes_R I_V) \cong \HH^i(K_W \otimes_R I_W) \fora i =
1,\ldots,n-3.
$$
\end{theorem}

In other words, the modules $\HH^i(K_V \otimes_R I_V)$ are
invariants of the CI-liaison class of $V$. They must vanish if $V$
is licci.

\begin{corollary} \label{cor-licci}
Let $V \subset \PP^n$, $n \geq 4$, be an \aCM\ subscheme of
codimension $3$. If $V$ is licci then $\HH^i(K_V \otimes_R I_V) =
0$ for all $i = 1,\ldots,n-3$.
\end{corollary}

\begin{proof}
Let $I \subset R$ be a complete intersection of codimension three.
Then we have the following isomorphisms (ignoring degree shifts)
$$
K_{R/I} \cong R/I;
$$
thus
$$
K_{R/I} \otimes_R I \cong I/I^2 \cong (R/I)^3.
$$
Since $\HH^i(R/I) = 0$ for $i \leq n-3$ because $R/I$ is \CM,
Theorem \ref{thm-ci-inv} proves the claim.
\end{proof}

For example, this result can used to reprove that the rational
normal curve in $\PP^4$ is not licci (cf.\ Example
\ref{ex-rat-n-curve}).

In \cite{KMMNP},  the previous theorem has been used to investigate
CI-liaison classes of curves on a Castelnuovo surface.

\begin{example}
Let $S \subset \PP^4$ be a general Castelnuovo surface, i.e.\ the
blow-up of a set of $8$ general points in $\PP^3$ embedded into
$\PP^4$ by the linear system $\mid 4 E_0 - 2 E_1 - E_2 - \ldots -
E_8 \mid$. Note that $S$ is an \aCM\ surface of degree $5$ which
contains a rational normal curve $C$ of $\PP^4$. Denote by $H_S$
the general hyperplane section of $S$. Furthermore, denote by
$C_j$ any curve in the linear system $\mid C + j H_S \mid$. Then
we have (cf.\ \cite{KMMNP}, Example 7.9)
\begin{itemize}
\item[(a)] The curve $C_j$ is not licci if $j \geq 0$.
\item[(b)] The curves $C_i$ and $C_j$ belong to different
CI-liaison classes whenever $1 \leq i < j$ and $j \geq 3$.
\end{itemize}
Since we know that all \aCM\ curves on $S$ are glicci (cf.\
Remark \ref{rem-acm-glicci}) we obtain that the G-liaison class of $C$
contains infinitely many CI-liaison classes.
\end{example}

So far CI-liaison invariants beyond the G-liaison invariants given
by Rao's correspondence are known only for \aCM\ subschemes. It
seems plausible to expect such additional invariants also for
non-\aCM\ subschemes. The problem of finding them deserves further
investigation.  Here is possibly the simplest situation.  In \cite{migliore-p4}
the following conjecture was made.

\begin{conjecture} \label{2 skew lines in p4}
If $C$ is a set of two skew lines in $\proj{4}$, spanning a hyperplane $H$, and
if $C'$ is another set of two skew lines in $\proj{4}$, spanning a hyperplane
$H'$, then $C$ is in the CI-liaison class of $C'$ if and only if $H = H'$.
\end{conjecture}

This conjecture would say that somehow the hyperplane
$H$ is a geometric invariant of the CI-liaison class of $C$, so
there must be some other algebraic invariant in addition to the
Rao module.  We have seen above that Hartshorne and Lesperance
independently showed that $C$ and $C'$ are in the same G-liaison
class, so this invariant would not hold for G-liaison.
\smallskip

We have seen that liaison in codimension two has two natural
generalizations in higher codimension: CI-liaison and G-liaison.
The former can be understood as a theory about divisors on
complete intersections while G-liaison is a theory about divisors
on \aCM\ schemes with property $G_1$. Thus, G-liaison is a much
coarser equivalence relation than CI-liaison. It has the advantage
that it is well suited for studying linear systems. The even
CI-liaison classes are rather small. In fact, it seems very
difficult to find enough invariants which would completely characterize an
even CI-liaison class.

It is also worth mentioning two disadvantages of G-liaison. The
first is related to our thin knowledge of \aG\ subschemes. Given a
subscheme $V$,  it is difficult to find ``good'' G-links of $V$;
i.e.\  ``good'' \aG\ subschemes $X$ containing $V$, where ``good''
often means small. For example, it is not too difficult to
determine the smallest degree of a complete intersection
containing $V$, while it is not known how to find an \aG\ subscheme
of smallest degree containing~$V$.

The second concerns lifting the information on hyperplane sections.
If $V, W \subset \PP^n$ are \aCM\ subschemes and $H \subset \PP^n$
is a general hyperplane such that $V \cap H$ and $W \cap H$ are
linked by the complete intersection $\bar{X} \subset H$ then there is a
complete intersection $X \subset \PP^n$ linking $V$ to $W$ such
that $\bar{X} = X \cap H$. The corresponding conclusion fails if
we replace ``complete intersection'' by ``\aG'' (cf.\ \cite{KMMNP},
Example 2.12).
\smallskip

In Section \ref{sec-relations} we defined geometric CI-linkage. It is also a
symmetric relation, thus its transitive closure is an equivalence relation
which is
essentially the same as CI-liaison. However, we have to be a little bit
careful what
we mean here. If $V$ is not a generic complete intersection then clearly it
does not
participate in a geometric CI-link. Thus, we make the following definition.

\begin{definition} \label{def-gci}
Let $H(c,n)$ denote the set of all equidimensional generic complete
intersections
  of $\PP^n$ of codimension $c$.
  \end{definition}

Note that this differs from the corresponding definition of Rao \cite{rao1}
not only
in allowing arbitrary codimension, but also in removing his assumption that the
schemes are \lCM.

Geometric CI-liaison is an equivalence relation on $H(c, n)$ while
CI-liaison is an
equivalence relation of the set of all equidimensional subschemes of
$\PP^n$ having
codimension $c$. But if we restrict the latter to $H(c, n)$ we get the
following:

\begin{theorem} \label{thm-geom-CI}
Algebraic and geometric CI-linkage generate the same equivalence relation
on $H(c,
n)$. That is, if $V, W \in H(c, n)$ are two generic complete intersections such
that there is a sequence of (algebraic) CI-links
$$
V \sim V_1 \sim \ldots V_s \sim W
$$
with all $V_i \in H(c, n)$ then there is a sequence of geometric CI-links
from $V$
to $W$.
\end{theorem}

For the proof we refer to \cite{KMMNP}, Theorem 4.14. The result
generalizes Rao's
Theorem 1.7 in \cite{rao1} which deals with the case $c=2$.
\smallskip

The last result leaves open the following problem. Suppose there are $V, W
\in H(c, n)$
such that there is a sequence of algebraic CI-links
$$
V \sim V_1 \sim \ldots V_s \sim W
$$
where some of the $V_i$ are {\it not} generic complete intersections. Is
there still
a sequence of geometric CI-intersections from $V$ to $W$?

The answer is known in codimension two. It uses the observation of Rao
(\cite{rao1},
Remark 1.5) that for a given $f \in I_V$ where $V \in H(2, n)$ there is
always a
form $g \in I_V$ of sufficiently large degree such that the complete
intersection
defined by $(f, g)$ links $V$ to a scheme $V'$ which is also a generic complete
intersection. Combining this fact with an analysis of the arguments which
establish
injectivity of Rao's correspondence in codimension two one gets the following:

\begin{theorem} \label{thm-geom-c-2}
Let $V, W \in H(2, n)$ be two subschemes such
that there is a sequence of (algebraic) CI-links
$$
V \sim V_1 \sim \ldots V_s \sim W.
$$
Then there is a sequence of geometric CI-links from $V$ to $W$.
\end{theorem}

For details of the proof we refer to \cite{KMMNP},Theorem 4.16.
\smallskip

It is an open question if the analogue of Theorem
\ref{thm-geom-c-2} is also true for subschemes of codimension $c
\geq 3$.


\section{Applications of Liaison} \label{application section}

In this section we mention some applications of liaison that have been
made in the literature.  It is not at all intended to be a complete list.

\subsection{Construction of \aG\ schemes with nice properties}

Here we describe in somewhat more detail the result of \cite{MN3} mentioned on
page \pageref{MN3 ref}.  It represents one of the few applications so far of
Gorenstein liaison as opposed to complete intersection liaison.

It is an open question to determine what Hilbert
functions are possible for Artinian Gorenstein graded $K$-algebras.  Indeed,
this seems to be intractable at the moment.  However, it was shown by Harima
\cite{harima} that the Hilbert functions of the Artinian Gorenstein graded
$K$-algebras {\em with the Weak Lefschetz property} (cf.\ Definition \ref{wlp
definition}) are precisely the SI-sequences (see page \pageref{def of
SI-sequence} for the definition).  Another open question is to determine the
possible Hilbert functions of reduced, \aG\ subschemes of $\proj{n}$ of any
fixed codimension.  Again, it is not clear if this problem can be solved or
not, but in the same way as the Artinian case, we have a partial result.  That
is, in \cite{MN3} it was shown that every SI-sequence gives rise to a reduced
union of linear varieties which is \aG\ and whose general Artinian reduction
has the Weak Lefschetz property.

\begin{remark} It would be very nice to show that every reduced \aG\ subscheme
has the property that its general Artinian reduction has the Weak Lefschetz
property.  If this were the case, then the result of \cite{MN3} would give a
classification of the Hilbert functions of reduced \aG\ subschemes of
$\proj{n}$, namely they would be those functions whose appropriate difference
is an SI-sequence.
\end{remark}

The construction given in \cite{MN3} is somewhat technical, and we give only
the main ideas.  One of the interesting points of this construction is that
it works in completely the opposite direction from the usual application of
liaison.  That is, instead of starting with a scheme $V$ and finding a suitable
\aG\ scheme $X$ containing it, we start with a (very reducible) \aG\ scheme $X$
and find a suitable subscheme $V$ to link using $X$.  Here are the main
steps of the
proof.

\begin{itemize}

\item[(a)] Suppose that we have a {\em geometric} link $V_1 \stackrel{X}{\sim}
V_2$, where $V_1$ (and hence also $V_2$) are \acm, and $X$ is \aG\ (not
necessarily a complete intersection).  Suppose you know the Hilbert
function of 
$V_1$ and of $X$.  Then using Corollary \ref{cor-dir-linked-id} we can write
the Hilbert function of $V_2$ (see also Example \ref{ex-linkage} (iii)).  From
the exact sequence
\[
0 \rightarrow I_X \rightarrow I_{V_1} \oplus I_{V_2} \rightarrow I_{V_1} +
I_{V_2} \rightarrow 0
\]
we also can get the Hilbert function of $R/(I_{V_1} + I_{V_2})$.

\item[(b)] Using induction on the codimension, we construct our \aG\ schemes
$X$ which are not complete intersections in general.  They have the following
properties.

\begin{itemize}

\item[(i)] They are {\em generalized stick figures}.  This means that they are
the reduced union of linear varieties of codimension $c$ (say), and no three
components meet in a linear variety of codimension $c+1$.  In the case of
curves, this is precisely the notion of a {\em stick figure}.

There are
several advantages to using generalized stick figures for $X$.  First, there
are many possible subconfigurations that we can link using $X$, if we can just
devise a way to find the ``right'' ones.  Second, any such link is guaranteed
to be geometric, since $X$ is reduced.  Third, after making such a link and
finding the sum of the linked ideals, the result is guaranteed to be reduced,
thanks to the fact that it is a generalized stick figure!  (This idea was used
earlier in \cite{GM5} for the case of CI-linked stick figure curves in
$\proj{3}$.)

\item[(ii)] Their Hilbert functions are ``maximal'' with a flat part in the
middle.  They are constructed inductively as a sum of G-linked ideals, by
finding a suitable subset with ``big'' Hilbert function, which in turn is
constructed by Basic Double G-Linkage.  For example, here are the
$h$-vectors of
the \aG\ schemes in low codimension:
\[
\begin{array}{ccccccccccccccccccccc}
\hbox{codim 2: \ \ } 1 & 2 & 3 & \dots & t-1 \\ \\
\hbox{codim 3: \ \ } 1 & 3 & 6 & \dots & \binom{t}{2} \\ \\
\hbox{codim 4: \ \ } 1 & 4 & 10 & \dots & \binom{t+1}{3}
\end{array}
\overbrace{
\begin{array}{cccccccccccccccccccc}
t & t & \dots & t \\ \\
\binom{t+1}{2} & \binom{t+1}{2} & \dots & \binom{t+1}{2} \\ \\
\binom{t+2}{3} & \binom{t+2}{3} & \dots & \binom{t+2}{3}
\end{array}
}^{\hbox{flat}}
\begin{array}{ccccccccccccccccc}
t-1 & \dots & 3 & 2 & 1 \\ \\
\binom{t}{2} & \dots & 6 & 3 & 1 \\ \\
\binom{t+1}{3} & \dots & 10 & 4 & 1
\end{array}
\]
\end{itemize}

\item[(c)] The schemes $X$ obtained in (b) will be used to link.  We will
assume that $\codim X = c-1$ and construct our schemes in codimension $c$.
Suppose that a desired SI-sequence $\underline{h}$ is given.  We use the
formula
of part (a) to work backwards, to determine the Hilbert function of an \acm\
subconfiguration $V_1 \subset X$ that would be needed to produce
$\underline{h}$ as a sum of linked ideals.

\item[(d)] We use our knowledge of the schemes $X$ to prove that an \acm\
scheme $V_1 \subset X$, as described in (c), in fact does exist.  This is the
most technical part of the proof.
\end{itemize}


\subsection{Smooth curves in $\proj{3}$}

\begin{itemize}

\item[A.] A long-standing problem, with many subtle variations, was to
determine
the possible pairs $(d,g)$ of degree and genus of smooth curves in $\proj{3}$
(or $\proj{n}$).  This was solved by Gruson and Peskine \cite{GP} for curves in
$\proj{3}$ and by Rathmann \cite{rathmann} for curves in $\proj{4}$ and
$\proj{5}$.  Substantial progress has been made by Chiantini, Ciliberto and Di
Gennaro \cite{CCD} in higher projective spaces.

One variation of this problem is to determine a bound for the
(arithmetic) genus of a non-degenerate, integral, degree $d$ curve $C \subset
\proj{3}$ lying on an irreducible surface $S$ of degree $k$, and to describe
the extremal curves.  This problem was solved by Harris \cite{harris}, who gave
a specific bound.  Furthermore, he showed that the curves which are extremal
with respect to this bound are precisely the curves residual to a plane curve
via certain complete intersections.  Note that they are thus \acm.  (A deeper
problem is to bound the genus of a smooth curve in $\proj{3}$ not lying on any
surface of degree $<k$.  There is much progress on this problem, beginning with
work of Hartshorne and Hirschowitz \cite{HH}.)

\item[B.] \label{UPP questions} Harris' work mentioned above used the Hilbert
function of the general hyperplane section of the curve $C$.  He showed that
the general hyperplane section must have the Uniform Position Property (see
Definition \ref{def of CB and UPP}).  (Note that Harris' proof of the uniform
position property for a general hyperplane section required characteristic
zero.  It has been proved in characteristic $p$ for $\proj{n}$, $n \geq 4$, by
Rathmann \cite{rathmann}.)   This led to natural questions:

\begin{itemize}
\item[Q1.] What are all the possible Hilbert functions for the general
hyperplane section of an integral curve in $\proj{3}$?  (Same question for
$\proj{n}$.)

\item[Q2.] What are all the possible Hilbert functions for the general
hyperplane section of an integral \acm\ curve in $\proj{3}$?  (Same question
for $\proj{n}$.)

\item[Q3.] What are all the possible Hilbert function of sets of points in
$\proj{2}$ with the Uniform Position Property?  (Same question for
$\proj{n-1}$.)

\item[Q4.] Do the questions above (for fixed $n$) have the same answer?
\end{itemize}

The answer to these questions is known for $n=3$, but open otherwise (see also
Section \ref{hf question section}).  The answer to Q4 is ``yes'' when
$n=3$, and
the Hilbert functions that arise are those of so-called {\em decreasing type}.
This means  the following.  Let
$Z$ be the set of points (either the hyperplane section of an integral curve or
a set of points with the Uniform Position Property).  Then the Hilbert
function of the Artinian reduction, $A$, of $R/I_Z$ looks as follows.  Let
$d_1$ be the degree of the first minimal generator of $I_Z$, and $d_2$ the
degree of the second. Note that $d_1 \leq d_2$. Let $r$ be the
Castelnuovo-Mumford regularity of ${\mathcal I}_Z$. Then
\[
h_A(t) =
\left \{
\begin{array}{ll}
t+1 & \hbox{if $t < d_1$} \\
d_1 & \hbox{if $d_1 \leq t \leq d_2 -1$} \\
\hbox{(strictly decreasing)} & \hbox{if $d_2 -1 \leq t \leq r$} \\
0 & \hbox{if $t \geq r$}
\end{array}
\right.
\]
Work on this problem was carried out in \cite{GP}, \cite{MR}, \cite{sauer}.
The interesting part is to construct an integral \acm\ curve with the desired
$h$-vector, and this was done in \cite{MR} by a nice application of liaison.  A
completely different approach, using lifting techniques, was carried out in
\cite{CO}.
\end{itemize}


\subsection{Smooth surfaces in $\proj{4}$, smooth threefolds in $\proj{5}$}

In the classification of smooth codimension two subvarieties (and by
Hartshorne's conjecture, we stop with threefolds in $\proj{5}$), it has
typically been the case that adjunction theory or other methods have been used
to narrow down the possibilities (see for instance \cite{BS-book}), and then
liaison has been used to construct examples.

We give an illustration of this idea by sketching a result of Mir\'o-Roig from
\cite{rosa3}.  A natural question is to determine the degrees $d$ for which
there exists a smooth, non-\acm\ threefold in $\proj{5}$.  It had been
shown by
  B\u anic\u a \cite{banica} that such threefolds exist for any odd $d \geq 7$
and for any even $d = 2k > 8$ with $k = 5s +1$, $5s+2$, $5s+3$ or $5s+4$.  It
had been shown by Beltrametti, Schneider and Sommese \cite{BSS} that
any smooth
threefold in $\proj{5}$ of degree 10 is \acm.

It remained to consider the case where $d = 10n$, $n \geq 2$.  Mir\'o-Roig
proved the existence of such threefolds using liaison.  Her idea was to begin
with well-known non-\acm\ threefolds in $\proj{5}$ and use the fact that the
property of being \acm\ is preserved under liaison.  In addition, she used the
following result of Peskine and Szpiro \cite{PS} to guarantee smoothness:

\begin{theorem}\label{smoothresid}
Let $X \subset \proj{n}$, $n \leq 5$, be a local
complete intersection of codimension two.  Let $m$ be a twist such that
${\mathcal I}_X (m)$ is globally generated.  Then for every pair $d_1 ,d_2 \geq
m$ there exist forms $F_i \in H^0 ({\mathcal I}_X (d_i )), \ i = 1,2$, such
that
the corresponding hypersurfaces $V_1$ and $V_2$ intersect properly and link $X$
to a variety $X'$.  Furthermore, $X'$ is a local complete intersection with no
component in common with $X$, and $X'$ is nonsingular outside a set of positive
codimension in Sing $X$.
\end{theorem}

\noindent (This special case of the theorem is quoted from
\cite{decker-popescu}, Theorem 2.1.)  Mir\'o-Roig  considered an \aBM\
threefold
$Y$ with locally free resolution
\[
0 \rightarrow {\mathcal O}_{\proj{5}} \oplus {\mathcal O}_{\proj{5}}(1)^3
\rightarrow \Omega^1 (3) \rightarrow {\mathcal I}_Y(6) \rightarrow 0
\]
(see also Example \ref{ex-rao-corr}).  Since ${\mathcal I}_Y (6)$ is globally
generated, Theorem \ref{smoothresid} applies.  Linking by two general
hypersurfaces of degrees 6 and 7, respectively, she obtains a smooth residual
threefold $X$ of degree 30, and using the mapping cone construction she obtains
the locally free resolution of ${\mathcal I}_X$.  Playing the same kind of
game, she is able to obtain from $X$ smooth threefolds of degrees $10n$, $n
\geq 5$, by linking $X$ using hypersurfaces of degree 10 and $n+3$.  The
remaining cases, degrees 20 and 40, are obtained by similar methods, starting
with different $Y$.


\subsection{Hilbert function questions} \label{hf question section}

We have seen above that liaison is useful for showing the existence of
interesting objects.  In this section we will see that liaison can sometimes be
used to prove non-existence results, as well as results which reduce the
possibilities.  For instance, we consider the question of describing the
possible Hilbert functions of sets of points in
$\proj{3}$ with the Uniform Position Property.

\begin{example}
Does there exist a set of points in $\proj{3}$ with the Uniform Position
Property and $h$-vector
\[
1 \ \ 3 \ \ 6 \ \ 5 \ \ 6,
\]
and if so, what can we say about it?  Suppose that such a set, $Z$, does
exist.  Note that the growth in the $h$-vector from degree 3 to degree 4 is
maximal, according to Macaulay's growth condition \cite{mr.macaulay}.  This
implies, thanks to \cite{BGM} Proposition 2.7, that the components
$[I_Z]_3$ and
$[I_Z]_4$ both have a GCD of degree 1, defining a plane $H$.  It also follows
using the same argument as \cite{BGM} Example 2.11 that $Z$ consists of either
14 or 15 points on $H$, plus 6 or 7 points not on $H$ (of which 4 or 5 are on a
line).  Such a $Z$ clearly does not have the Uniform Position Property!
\end{example}

\begin{example}
Does there exist a set of points in $\proj{3}$ with the Uniform Position
Property and $h$-vector
\[
1 \ \ 3 \ \ 6 \ \ 5 \ \ 5,
\]
and if so, what can we say about it?  Let $Z$ be such a set.
In this case we do not have maximal growth from degree 3 to degree 4, but we
again consider the component in degree~3.  This time we will not have a GCD,
but we can consider the base locus of the linear system $|[{I}_Z]_3|$.  Suppose
that this base locus is zero-dimensional.  Then three
general elements of $[I_Z]_3$ give a complete intersection, $I_X =
(F_1,F_2,F_3)$.  This means that $Z$ is linked by $X$ to a zeroscheme $W$, and
we can make a Hilbert function ($h$-vector) calculation (cf.\ Corollary
\ref{cor-dir-linked-id} and Example \ref{ex-linkage} (c)):
  \medskip
\begin{center}
\begin{tabular}{c|ccccccccc}
degree & 0 & 1 & 2 & 3 & 4 & 5 & 6 & 7  \\
\cline{1-9}
$R/I_X$ &  1 & 3 & 6 & 7 & 6 & 3 & 1 & 0 \\[2pt]
$R/I_Z$ &  1 & 3 & 6 & 5 & 5 & 0 & 0 & 0 \\ [2pt]
$R/I_W$ &  0 & 0 & 0 & 2 & 1 & 3 & 1 & 0
\end{tabular}
\end{center}

\medskip
This means that the residual, $W$, has $h$-vector $1 \ 3 \ 1 \ 2$, which is
impossible (it violates Macaulay's growth condition).

Thus we are naturally led to look for an example consisting of a set of 20
general points, $Z$, on an irreducible curve $C$ of degree 5.  (We do not
justify this, although similar considerations can be found in the proof of
Theorem 4.7 of \cite{BGM}, but we hope that it is clear that this is the
natural place to look, even if it is not clear that it is the {\em only} place
to look.)  The Hilbert function of $Z$ has to agree with that of $Z$ up to
degree 4.  One can check that a general curve $C$ of degree 5 and genus 1 will
do the trick (and no other will).  Hence the desired set of points does exist.
\end{example}


\subsection{Arithmetically Buchsbaum curves specialize to stick figures}

We have seen how to use Liaison Addition to construct minimal \aBM\ curves
(Example \ref{construction of aBM curves}) and how to use the Lazarsfeld-Rao
property to give all the possible $(d,g)$ combinations possible for \aBM\
curves (Remark \ref{all possible d,g}).  Now we sketch how these ideas were
refined in \cite{BM5} and applied to show that every \aBM\ specializes to a
stick figure.  This is a special case of the Zeuthen problem, a long-standing
problem that was solved a few years ago by Hartshorne \cite{hart-zeuthen}.  The
general question is whether every smooth curve in $\proj{3}$ specializes to a
stick figure, and Hartshorne showed that the answer is ``no."  This makes it
more interesting that the answer is ``yes'' for \aBM\ curves.

Let $C$ be an \aBM\ curve.  The basic idea here is that the Lazarsfeld-Rao
property provides the desired deformation, if we can produce a stick figure
using basic double links which is cohomologically the same as $C$.  So there
are two parts to the story.  First we have to produce a minimal element which
is a stick figure, and second we have to study basic double links and show that
we can always keep producing stick figures.

For the first part, it is a refinement of the construction given in Example
\ref{construction of aBM curves}.  Skipping details, we merely note here that
if $C_1$ and $C_2$ are both pairs of skew lines chosen generically, then $F_1$
and $F_2$ can be chosen to be unions of planes, and for a sufficiently general
choice, the curve $C$ constructed by Liaison Addition will be a stick figure.
To see that this procedure can give a minimal element for {\em any} Buchsbaum
even liaison class is somewhat more technical, but is an extension of this
idea.

For the second part, recall that a basic double link is obtained by starting
with a curve $C$ and a surface $F$ containing $C$, and taking the union $Y$ of
$C$ and a general hyperplane section of $F$.  If $C$ is a union of lines and
$F$ is a union of planes then clearly $Y$ will also be a union of lines.  The
first problem is to show that we can always arrange that there exists a surface
$F$ which is a union of planes.  For instance, if $C$ is a union of $\geq 3$
skew lines on a quadric surface (this is not \aBM, but gives the idea), and if
we want $\deg F = 2$, then $F$ clearly cannot be chosen to be a union of
planes.  So we have to show that a union of planes can always be obtained in
our case.  But there is a more subtle problem.

For example, suppose that $C$ is a set of two skew lines, and suppose that we
make a sequence of three basic double links using $F_1$, $F_2$ and $F_3$ of
degrees 20, 15 and 4 respectively, obtaining curves $Y_1$, $Y_2$ and $Y_3$ of
degrees 22, 37 and 41 respectively.  A little thought shows that one cannot
avoid that $Y_3$ have a triple point!  (The key is that $\deg F_1 > \deg F_2 >
\deg F_3$.)  Thus this sequence of basic double links {\em cannot} yield a
stick
figure.

The solution to this dilemma is to show that there is a {\em
cohomologically equivalent} sequence of basic double links using surfaces $G_1,
G_2, G_3$ with $\deg G_1 \leq \deg G_2 \leq \deg G_3$.  Then the type of
problem described in the last paragraph does not occur.  Again, the details are
technical, and we refer the reader to \cite{BM4} and \cite{BM5}.


\subsection{The minimal free resolution of generic forms}

An important problem, variations of which have been studied by many people, is
to describe the Hilbert function or minimal free resolution of an ideal $I
\subset R = K[x_1,\dots,x_n]$ generated by a general set of forms of fixed
degrees (not necessarily all the same).  The answer to the Hilbert function
problem has been conjectured by Fr\"oberg and we will not describe it
here.  It
is known to hold when $n \leq 3$ and when the number of generators is $n+1$.

For the minimal free resolution, the answer has been conjectured by Iarrobino.
At the heart of this is the idea that if the forms are general then there
should be no ``ghost terms'' in the minimal free resolution, i.e.\ there should
be no summand $R(-t)$ that appears in consecutive free modules in the
resolution.  One can see immediately that this is too optimistic, however. For
instance, if $I$ has two generators of degree 2 and one of degree 4 then there
is a term $R(-4)$ corresponding to a first syzygy and a term $R(-4)$
corresponding to a generator.  So the natural conjecture is that apart from
such terms which are forced by Koszul relations, there should be no ghost
terms.

This was proved to be false in \cite{MMR2}.  A simple counterexample is the
case of four generators in $K[x_1,x_2,x_3]$ of degrees 4,4,4 and 8.  The
minimal free resolution turns out to be
\[
0 \rightarrow
\left (
\begin{array}{c}
R(-10) \\
\oplus \\
R(-11)^2
\end{array}
\right )
\rightarrow
\left (
\begin{array}{c}
R(-8)^3 \\
\oplus \\
R(-9)^2 \\
\oplus \\
R(-10)
\end{array}
\right )
\rightarrow
\left (
\begin{array}{c}
R(-4)^3 \\
\oplus \\
R(-8)
\end{array}
\right )
\rightarrow R \rightarrow R/I \rightarrow 0
\]
The term $R(-8)$ that does not split arises from Koszul relations, as above,
but the summand $R(-10)$ shared by the second and third modules also does
not split and this does {\em not} arise from Koszul relations.

The paper \cite{MMR2} made a general study of the minimal free resolution of
$n+1$ general forms in $R$.  The minimal free resolution was obtained in many
cases (depending on the degrees of the generators) and the main tools were
liaison and a technical lemma from \cite{MN3} giving a bound on the graded
Betti  numbers for Gorenstein rings.  The key to this work is Corollary
\ref{cor-aci-gor} above, which says that our ideal $I$ can always be directly
linked to a Gorenstein ideal.

Here is the basic idea.  Knowing the Hilbert function for the $n+1$ general
forms leads to the Hilbert function of the linked Gorenstein ideal.  The
technical lemma of \cite{MN3} then gives good bounds for the graded Betti
numbers of the linked Gorenstein ideal, and in fact these bounds can often be
shown to be sharp.  Then the mapping cone obtained from the first sequence of
Lemma \ref{lem-st-ex-seq} can be used to give a free resolution of
$R/I$.  One can then determine to what extent this resolution is minimal.  In
particular, ghost terms in the minimal free resolution of the Gorenstein ideal
translate to ghost terms in the minimal free resolution of $I$.  Especially
when $n=3$, we can often arrange ghost terms for the Gorenstein ideal (thanks
to the Buchsbaum-Eisenbud structure theorem \cite{BuEi} and the work of Diesel
\cite{diesel}).


  \section{Open Problems} \label{open prob sect}

In this section we collect the open questions that were mentioned in the
preceding sections, and add some more.

\begin{enumerate}

\item Describe the Hilbert functions for general hyperplane sections of
integral curves in $\proj{n}$ ($n \geq 4$) and for sets of points in
$\proj{n-1}$ with the Uniform Position Property.  (See the discussion starting
on page \pageref{UPP questions}.)

\item Find a description of all the possible Hilbert functions of
Artinian Gorenstein graded $K$-algebras.  Find a description of all the
possible Hilbert functions of reduced \aG\ subschemes of $\proj{n}$.  (Is the
answer to this last question precisely the SI-sequences?)

\item  Classify  the possible graded Betti
numbers for Gorenstein algebras in codimension $\geq 4$.  See Questions
\ref{possible artin gor hf} and \ref{possible reduced gor hf} and the
discussion following them.

\item \label{stci question} It is an old problem (see e.g.\
\cite{hartshorne-book} Exer.\ 2.17 (d)) whether every irreducible curve $C
\subset \proj{3}$ is a set-theoretic complete intersection.  It is not true
that a curve which is a set-theoretic complete intersection must be \ACM\ (see
e.g.\ \cite{rao3}).  However, the first
author has conjectured that such a curve must be linearly normal.  Some
progress
in this direction was achieved by Jaffe \cite{jaffe}.  In the first draft of
these notes we made the comment here that we were not aware even of a curve
which is a set-theoretic complete intersection but is not self-linked.
However, R.\ Hartshorne has provided us with an example, which we have recorded
in Example \ref{hartshorne ex} below.

\item Find conditions that are necessary and sufficient for two schemes in
codimension $\geq 3$ to be evenly CI-linked or evenly G-linked.

\item \label{gor g-linked?} In particular, is it true that two \ACM\ schemes
of the same codimension are G-linked in finitely  many steps?  As an important
first case, is it true that two arithmetically Gorenstein subschemes of the
same
codimension are G-linked in finitely many steps?

\item Extend the known CI-liaison invariants for \aCM\ subschemes
(cf., e.g., Theorem \ref{thm-ci-inv})
to non-\aCM\ subschemes which allow one to distinguish
CI-liaison classes within an even G-liaison class of a non-\aCM\
subscheme.

\item Compare the equivalence relations generated by geometric G-liaison and
(algebraic) G-liaison on the set of subschemes of $\PP^n$ having
codimension $c$ and
being generically Gorenstein. (cf.\ Section \ref{compare section} for
results in the
case of CI-liaison.)

\item Find a structure theorem similar to the LR-property that holds for
G-liaison or for CI-liaison in higher codimension.

\item Establish upper and lower bounds for the various  minimal shifts
   attached to an equidimensional scheme
   (cf.\ Remark  \ref{rem-shifts}).

\item Find a good concept for minimal elements in an even G-liaison class
   (cf.\ Section \ref{str-even-liaison}).

\item Find conditions like the theorem of Peskine and Szpiro \cite{PS}
(cf.\ Theorem \ref{smoothresid}) which
guarantee that a G-linked residual scheme is smooth (in the right
codimension). Find applications of this to the classification of smooth
codimension 3 subschemes. See \cite{Mig-Pet} for more on this idea.

\end{enumerate}

\begin{example} \label{hartshorne ex}
In an earlier draft of these notes we asked if there is any smooth
curve in $\proj{3}$ which is a set-theoretic complete intersection but not
self-linked.  We believed that there should be such a curve, but were not aware
of one.  This example is due to Robin Hartshorne, who has kindly allowed us to
reproduce it here.

A curve is self-linked if it is a set-theoretic complete intersection of
multiplicity  2.  So here we will construct, for every integer $d >0$, a
smooth
curve in $\proj{3}$ that is set-theoretically the complete intersection of
multiplicity $d$, but of no  lower multiplicity.

Start with a smooth plane curve of degree $d$, having a $d$-fold
inflectional tangent at a point $P$. Let $X$ be the cone over that
curve in $\proj{3}$. Let $L$ be the cone over $P$. Then $L$ is a line on $X$,
$dL$ is a complete intersection on $X$, and no lower multiple of $L$ is a
complete intersection of $X$  with another surface. Now let $C$ be a smooth
curve in the  linear system $|L + mH|$ on $X$, for $m \gg 0$.

Note that $dC$ is linearly equivalent to $dL + mdH = (md+1)H$. Therefore $dC$
is the intersection of $X$ with another surface in $\proj{3}$, and so
$C$ is a set-theoretic complete intersection of multiplicity $d$. Note that no
smaller multilple of $C$  is the complete intersection of $X$ with anything
else, because $eC$ for $e < d$ is  not a Cartier divisor on $X$. But could $eC$
be an intersection  of two other surfaces? Since $C$ has degree $md+1$, if $F$
is any  other surface containing $C$, then the degree of $X \cdot F$ is $d
\cdot \deg F$, so  $\deg F > m$. So if $C$ is the set-theoretic complete
intersection of $F$ and $G$, then $\deg F \cdot G$  is $> m^2$, and the
multiplicity of the structure on $C$ is  $> m^2/(md+1)$ , which for $m \gg
0$ is
$> d$. (In fact, to obtain $m^2/(md+1) > d$, i.e.\ $m(m-d^2) > d$, it is enough
to take $m > d^2$.)
\end{example}

\end{document}